\documentclass[10pt]{article}
\pdfoutput=1 

\usepackage{amsmath, amssymb}
\usepackage{enumerate, caption, graphicx}
\usepackage[letterpaper,text={6.75in,9.5in},centering]{geometry}
\usepackage[numbers,comma,square,sort&compress]{natbib}

\newtheorem{Lemma}{Lemma}[section]
\newtheorem{Corollary}[Lemma]{Corollary}
\newtheorem{Definition}[Lemma]{Definition}
\newtheorem{Hypothesis}{Hypothesis}
\newtheorem{Proposition}[Lemma]{Proposition}
\newtheorem{Remark}[Lemma]{Remark}
\newtheorem{Theorem}{Theorem}

\newcommand{\qed}{\hspace*{\fill}$\rule{0.3\baselineskip}{0.35\baselineskip}$}
\newenvironment{Proof}[1][.]{\begin{trivlist}\item[]\textbf{Proof#1 }}{\qed\end{trivlist}}
\newenvironment{Acknowledgment}{\begin{trivlist}\item[]\textbf{Acknowledgments }}{\end{trivlist}}

\setcaptionmargin{0.25in}

\makeatletter\@addtoreset{equation}{section}\makeatother

\def\Re{\mathop\mathrm{Re}\nolimits}
\def\Im{\mathop\mathrm{Im}\nolimits}
\def\Rg{\mathop\mathrm{Rg}\nolimits}
\def\errfn{\mathop\mathrm{errfn}\nolimits}
\def\sgn{\mathop\mathrm{sgn}\nolimits}
\def\Span{\mathop\mathrm{span}\nolimits}
\def\graph{\mathop\mathrm{graph}\nolimits}

\newcommand{\rmd}{\mathrm{d}}
\newcommand{\rme}{\mathrm{e}}
\newcommand{\rmi}{\mathrm{i}}

\begin{document}

\title{Nonlinear stability of time-periodic viscous shocks}

\author{Margaret Beck\footnotemark[1]\\
Division of Applied Mathematics\\
Brown University\\
Providence, RI~02912, USA
\and
Bj\"orn Sandstede\footnotemark[1]\\
Division of Applied Mathematics\\
Brown University\\
Providence, RI~02912, USA
\and
Kevin Zumbrun\\
Department of Mathematics\\
Indiana University\\
Bloomington, IN 47405, USA
}

\renewcommand{\thefootnote}{\fnsymbol{footnote}}
\footnotetext[1]{The majority of this work was done while MB and BS were affiliated with the Department of Mathematics, University of Surrey, Guildford, GU2~7XH, UK}
\renewcommand{\thefootnote}{\arabic{footnote}}

\date{\today}
\maketitle

\begin{abstract}
In order to understand the nonlinear stability of many types of time-periodic travelling waves on unbounded domains, one must overcome two main difficulties: the presence of embedded neutral eigenvalues and the time-dependence of the associated linear operator. This problem is studied in the context of time-periodic Lax shocks in systems of viscous conservation laws. Using spatial dynamics and a decomposition into separate Floquet eigenmodes, it is shown that the linear evolution for the time-dependent operator can be represented using a contour integral similar to that of the standard time-independent case. By decomposing the resulting Green's distribution, the leading order behavior associated with the embedded eigenvalues is extracted. Sharp pointwise bounds are then obtained, which are used to prove that time-periodic Lax shocks are linearly and nonlinearly stable under the necessary conditions of spectral stability and minimal multiplicity of the translational eigenvalues. The latter conditions hold, for example, for small-oscillation time-periodic waves that emerge through a supercritical Hopf bifurcation from a family of time-independent Lax shocks of possibly large amplitude.
\end{abstract}

\begin{small}
\setcounter{tocdepth}{1}
\tableofcontents
\end{small}
\newpage


\section{Introduction}\label{introduction}

We are interested in solutions to partial differential equations on unbounded domains that are time-periodic in an appropriate comoving frame and spatially asymptotic to constants. Such solutions, which we refer to as time-periodic waves, are therefore of the form
\begin{equation} \label{profile}
u(x,t) = \bar{u}(x-ct,t) \quad\mbox{with}\quad
\lim_{y\to\pm\infty} \bar{u}(y,t) = u_\pm, \quad 
\bar{u}(y,t+T) = \bar{u}(y,t),
\end{equation}
where $c$ is the speed of the wave, and $T>0$ is the time period. We refer to travelling waves $\bar{u}(x-ct)$ as time-independent waves because they become stationary in the comoving frame. Time-periodic solutions to partial differential equations arise in many applications. Examples are pulsating, spinning, and cellular instabilities in detonation waves, which have been investigated in \cite{TexierZumbrun06_gas,KasimovStewart02}, and interfaces that connect wave trains in reaction-diffusion systems, which have been analyzed, for instance, in \cite{SandstedeScheel04}. The purpose of this paper is to prove that spectral stability of time-periodic waves in viscous conservation laws implies their nonlinear stability. In the examples given above, two difficulties arise that prevent us from applying standard theory: the linearization about the wave is time-periodic, and its Floquet spectrum contains essential spectrum up to the origin.

We will focus our analysis on time-periodic Lax shocks, which arise in systems of viscous conservation laws. This setting is useful not only because of its relation to real physical models that exhibit periodic phenomena, but also because the underlying structure within these systems provides intuition and aids in their analysis. As discussed in \cite{TexierZumbrun05,TexierZumbrun06_framework,TexierZumbrun06_gas, SandstedeScheel06}, time-periodic viscous Lax shocks may arise through Hopf bifurcations from stationary shock waves. Spectral stability of the bifurcating waves has been treated in \cite{SandstedeScheel06}, where it was shown that they are spectrally stable if the bifurcation is supercritical and spectrally unstable if the bifurcation is subcritical. Our goal is to prove nonlinear stability of arbitrary time-periodic Lax shocks, possibly far away from any stationary shocks. In the stationary case, pointwise Green's function estimates have proved to be very useful in establishing nonlinear stability; see, for instance, \cite{HowardZumbrun06, MasciaZumbrun03}. We shall use ideas from spatial dynamics, and in particular the exponential-dichotomy theory for ill-posed elliptic problems developed in \cite{PeterhofSandstedeScheel97, SandstedeScheel01, SandstedeScheel06}, to extend this approach to the time-periodic case.

\begin{figure}[t]
\centering\includegraphics{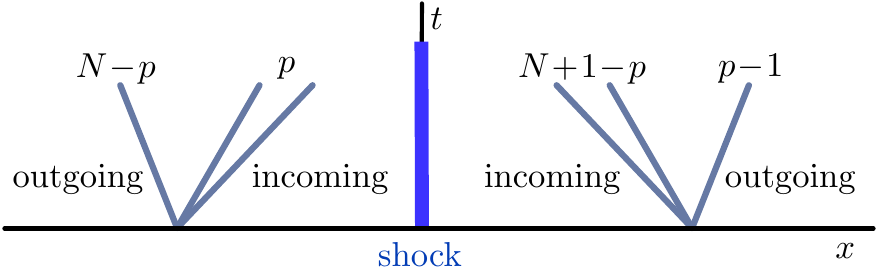}
\caption{The $N+1$ incoming and $N-1$ outgoing far-field characteristics $x=a_k^\pm t$ of the linearization about a Lax $p$-shock for a system of $N$ viscous conservation laws are illustrated for $p=2$ and $N=3$.}
\label{F:p-shocks}
\end{figure}

We now state our hypotheses and results in detail. Consider a time-periodic viscous shock profile $u(x,t)=\bar{u}(x-ct,t)$, as in (\ref{profile}), of a parabolic system of conservation laws
\begin{equation} \label{main}
u_t + f(u)_x = u_{xx}, \qquad x\in\mathbb{R},\quad u\in\mathbb{R}^N,
\end{equation}
where $f\in C^3$. Working in a coordinate system that moves along with the shock and appropriately rescaling time, we may, without loss of generality, consider a standing profile $\bar{u}(x,t)$ with minimal temporal period $2\pi$. Accordingly, we take $c=0$ and $T=2\pi$ from now on. We shall assume that the profile $\bar{u}(x,t)$ is a Lax $p$-shock:

\begin{Hypothesis}\label{H1}
The ordered eigenvalues $a_1^\pm<\ldots<a_N^\pm$ of $f_u(u_\pm)$ are real, distinct, and nonzero, and there is a number $p\in\{1,\dots,N\}$ so that
$a_{N-p}^-<0<a_{N-p+1}^-$ and $a_{N-p+1}^+<0<a_{N-p+2}^+$.
\end{Hypothesis}
The eigenvalues $a_k^\pm$ of $f_u(u_\pm)$ determine the characteristics $x=a_k^\pm t$ of the linear system
\begin{equation*}
v_t + f_u(u_\pm) v_x = 0,
\end{equation*}
whose general solution is a linear combination of solutions of the form $v(x-a_k^\pm t)$. We say that the characteristics corresponding to $a_k^-<0$ and $a_k^+>0$ are \emph{outgoing}, while the characteristics associated with $a_k^->0$ and $a_k^+<0$ are called \emph{incoming}, and refer the reader to Figure~\ref{F:p-shocks} for an illustration. For future use, we denote by $l_k^\pm$ and $r_k^\pm$ the left and right eigenvectors of $f_u(u_\pm)$ associated with the characteristic speeds $a_k^\pm$ and normalize these vectors so that $\langle l_j^-,r^-_k\rangle=\langle l_j^+,r^+_k\rangle=\delta_{jk}$.

\begin{Definition}\label{notation}
Throughout this paper, we will often use the notation $a_\mathrm{out}^-$ and $a_\mathrm{out}^+$ to denote the outgoing characteristics $a_k^-<0$ and $a_k^+>0$ of $u_-$ and $u_+$, respectively. Similarly, $a_\mathrm{in}^-$ and $a_\mathrm{in}^+$ denote the incoming characteristics $a_k^->0$ and $a_k^+<0$. We often write
\[
\sum_{a_\mathrm{out}^-} := \sum_{a_k^-<0} \qquad\qquad
\sum_{a_\mathrm{in}^-}  := \sum_{a_k^->0} \qquad\qquad
\sum_{a_\mathrm{out}^+} := \sum_{a_k^+>0} \qquad\qquad
\sum_{a_\mathrm{in}^+}  := \sum_{a_k^+<0}
\]
to emphasize the interpretation of these sums in terms of characteristics. Similar notation will also be used for the associated left and right eigenvectors.
\end{Definition}

In order to state our spectral assumptions, we linearize \eqref{main} about the time-periodic shock $\bar{u}(x,t)$ and obtain
\begin{equation}\label{linearized}
u_t = u_{xx} - [f_u(\bar{u}(x,t))u]_x,
\end{equation}
where $f_u(\bar{u}(x,t))$ is $2\pi$-periodic in $t$. Because of the time-periodicity of the linear operator on the right-hand side, the appropriate characterization of spectral stability is in terms of its Floquet spectrum. Let $\Phi_{2\pi}$ denote the time-$2\pi$ map associated with the semiflow of (\ref{linearized}) posed on $L^2(\mathbb{R},\mathbb{R}^N)$. The Floquet spectrum is defined as
\begin{equation}\label{e:floq}
\Sigma = \{ \sigma\in\mathbb{C}:\; \rme^{2\pi\sigma}\in\mbox{ spectrum of } \Phi_{2\pi} \}.
\end{equation}
It follows from Hypothesis~\ref{H1} that $\bar{u}_x(x,t)$ and $\bar{u}_t(x,t)$ decay exponentially to zero as $|x|\to\infty$, uniformly in time; see \cite{SandstedeScheel06}. Thus, the origin $\sigma=0$ lies in the spectrum $\Sigma$ and belonging to it are two exponentially decaying eigenfunctions, namely $\bar{u}_x$ and $\bar{u}_t$, which reflect the translation-invariance of the underlying equations with respect to space and time. We will assume that the origin is the only part of the spectrum that lies in the closed right half-plane other than the spectrum at $\rmi\mathbb{Z}$ that results from the fact that Floquet exponents are not unique.

\begin{Definition}\label{spectralstab}
The profile $\bar{u}$ is said to be \emph{spectrally stable} if 
\begin{enumerate}[\bf (S1)]
\item The time-$2\pi$ map $\Phi_{2\pi}$ associated with (\ref{linearized}) on $L^2(\mathbb{R}, \mathbb{R}^N)$ has no eigenvalues in $\{\Re\sigma\geq0\}\setminus\rmi\mathbb{Z}$.
\item Equation (\ref{linearized}) has exactly two linearly independent localized solutions with period $2\pi$ in $t$.
\item The outgoing characteristics and the jump $[\bar{u}]:=u_+-u_-$ across the shock $\bar{u}$ are linearly independent:
\begin{equation}\label{majdacondition}
\det(r_1^-,\dots,r_{N-p}^-,[\bar{u}],r_{N+2-p}^+,\dots,r_N^+) \neq 0.
\end{equation}
\item The formal adjoint equation
\[
w_t = w_{xx} + f_u^T(\bar{u}(x,t)) w_x
\]
has a non-constant, bounded, time $2\pi$-periodic solution $w(x,t)$. Setting $\psi_2(x,t):=w(x,-t)$, we assume that
\begin{equation}
\int_\mathbb{R}\int_0^{2\pi} \langle\psi_2(x,t),\bar{u}_t(x,t)\rangle\,\rmd t\,\rmd x \neq 0
\label{E:transverse_connection}
\end{equation}
and $\lim_{x\to\pm\infty}\psi_2(x,t)=:\psi_\pm\in L_\mathrm{in}^\pm$, where $L_\mathrm{in}^\pm$ denote the subspaces of $\mathbb{R}^N$ spanned by the left eigenvectors $l_\mathrm{in}^\pm$ of $f_u(u_\pm)$ associated with the incoming characteristics $a_\mathrm{in}^\pm$.
\end{enumerate}
\end{Definition}

Hypothesis~(S2) guarantees that the geometric multiplicity of the zero eigenvalue is minimal. Condition~(S3) is the same as for stationary Lax shocks, where it is known as the Liu--Majda condition; see \cite{ZumbrunHoward98} for more details. Furthermore, (S3) implies that there is a unique, up to scalar multiples, nonzero $\psi_1\in\mathbb{R}^N$ so that $\psi_1$ is perpendicular to the outgoing characteristic eigenvectors that appear in the determinant (\ref{majdacondition}). This constant $\psi_1$ is the adjoint eigenfunction associated with $\bar{u}_x$, and (S3) implies that
\begin{equation}
\int_0^{2\pi}\int_\mathbb{R} \langle\psi_1,\bar{u}_x(x,t)\rangle\,\rmd x\,\rmd t
= \langle \psi_1,[\bar{u}] \rangle \neq 0.
\label{E:transverse_connection_2}
\end{equation}
Given (S2), conditions (S3)-(S4) furthermore imply that zero is a root of order two of an appropriate Evans function. The integrals in (\ref{E:transverse_connection}) and (\ref{E:transverse_connection_2}) can be interpreted as the Melnikov integrals associated with speed and frequency for the spatial-dynamics existence problem for $\bar{u}(x,t)$; see \cite{SandstedeScheel01,SandstedeScheel06}.

Our first result, paralleling the theory of the time-independent case \cite{ZumbrunHoward98,MasciaZumbrun03}, is the equivalence of linearized and spectral stability.

\begin{Theorem}\label{stabcrit}
Assume \ref{H1} and pick $\rho\geq0$, then spectral stability is equivalent to linearized stability in $L^1\cap H^\rho$ (that is, each solution of (\ref{linearized}) with initial data in $L^1\cap H^\rho$ converges in this space to $\Span\{\bar{u}_x,\bar{u}_t\}$ as $t\to\infty$).
\end{Theorem}

Our next theorem is the main result of this paper: it asserts that spectral stability implies nonlinear stability and gives detailed pointwise estimates for how perturbations decay as $t\to\infty$. Before we can state this result, we need additional notation. Let
\begin{eqnarray*}
\theta_\mathrm{gauss}(x,t) & := &
\sum_{a_\mathrm{out}^-} \frac{\rme^{-|x-a_\mathrm{out}^-t|^2/M t}}{\sqrt{1+t}} +
\sum_{a_\mathrm{out}^+} \frac{\rme^{-|x-a_\mathrm{out}^+t|^2/M t}}{\sqrt{1+t}} \\
\theta_\mathrm{inner}(x,t) & := &
\frac{1}{\sqrt{1+|x|+t}} \left(
\sum_{a_\mathrm{out}^-} \frac{1}{\sqrt{1+|x-a_\mathrm{out}^-t|}} +
\sum_{a_\mathrm{out}^+} \frac{1}{\sqrt{1+|x-a_\mathrm{out}^+t|}}
\right)\\
\theta_\mathrm{outer}(x,t) & := &
\frac{1}{(1+|x-a_1^-t|+\sqrt{t})^{\frac32}} +
\frac{1}{(1+|x-a_N^+t|+\sqrt{t})^{\frac32}},
\end{eqnarray*}
where $M>0$ is a sufficiently large constant, and define
\[
\chi(x,t) := \left\{ \begin{array}{lcl}
1 & \quad & x\in[a_1^-t,a_N^+t] \\
0 & & x\notin[a_1^-t,a_N^+t]
\end{array}\right.
\]
to be the characteristic function of the characteristic cone $[a_1^-t,a_N^+t]$.

\begin{Theorem}\label{nonlin}
Define the weighted norm $\|v\|_{H^3_w}:=\|(1+x^2)^{\frac34}v\|_{H^3}$. If \ref{H1} and spectral stability (S1)--(S4) hold, then the profile $\bar{u}$ is nonlinearly stable with respect to initial perturbations $v_0$ for which $\|v_0\|_{H^3_w}$ is sufficiently small. More precisely, there exist constants $C>0$ and $\delta>0$ such that, for each $v_0$ with $\|v_0\|_{H^3_w}<\delta$, there exist functions $(q,\tau)(t)$ and constants $(q_*,\tau_*)$ so that, for all $x\in\mathbb{R}$ and $t\geq0$, we have
\begin{eqnarray}\label{pointwise}
|u(x,t)-\bar{u}(x-q_*-q(t),t-\tau_*-\tau(t))| & \leq &
C \|v_0\|_{H^3_w} \left[ \theta_\mathrm{gauss}+\chi\theta_\mathrm{inner}+(1-\chi)\theta_\mathrm{outer} \right](x,t) \\ \nonumber
\|u(\cdot,t)-\bar{u}(\cdot-q_*-q(t),t-\tau_*-\tau(t))\|_{H^3} & \leq & C \|v_0\|_{H^3_w}
\end{eqnarray}
and
\[
|(q_*,\tau_*)| + (1+t)^{\frac12}|(q,\tau)(t)| + (1+t)|(\dot{q},\dot\tau)(t)| \leq C \| v_0\|_{H^3_w},
\]
where $u(x,t)$ denotes the solution of \eqref{main} with initial data $u_0(x)=\bar{u}(x,0)+v_0(x)$.
\end{Theorem}

The pointwise bound \eqref{pointwise} yields as a corollary the sharp $L^p$ decay rate
\[
\|u(\cdot,t)-\bar{u}(\cdot-q_*-q(t),t-\tau_*-\tau(t))\|_{L^p} \leq C \| v_0\|_{H^3_w}
(1+t)^{-\frac{1}{2}(1-\frac{1}{p})}, \qquad 1\leq p\leq\infty.
\]
The statement of Theorem~\ref{nonlin} can be understood as follows. First, for initial data sufficiently close to the underlying time-periodic wave, solutions to the full nonlinear equation (\ref{main}) will converge to an appropriate space and time translate of the wave. In addition, the perturbation to the solution profile will be of a certain form, given by the functions $\theta_\mathrm{gauss}$, $\theta_\mathrm{inner}$, and $\theta_\mathrm{outer}$, as it decays to zero. The function $\theta_\mathrm{gauss}$ consists of Gaussians that move along the outgoing characteristics at speeds $a_j^-<0$ and $a_j^+>0$. The function $\theta_\mathrm{inner}$ accounts for the nonlinear interactions that occur within the characteristic cone $[a_1^-t,a_N^+t]$ encoded in the characteristic function $\chi$. The algebraically decaying tail of the initial data, outside the characteristic cone, is captured by the function $\theta_\mathrm{outer}$.

Theorem~\ref{nonlin} is obtained as a consequence of detailed estimates of the solution operator of the linearization
\begin{equation}\label{linearov}
u_t = u_{xx} - [f_u(\bar{u}(x,t))u]_x.
\end{equation}
More precisely, we shall derive pointwise bounds of the \emph{Green's distribution} $\mathcal{G}(x,t;y,s)$ of (\ref{linearov}), which is its fundamental solution, defined as the solution at $(x,t)$ of (\ref{linearov}) with initial data $\delta(x-y)$ at time $s$, where $\delta(x-y)$ is the delta-function centered at $y$. Recall that $\bar{u}(x,t)$ converges to the constant vectors $u_\pm$ exponentially fast as $x\to\pm\infty$, and that $a_\mathrm{out}^\pm$ and $a_\mathrm{in}^\pm$ are, respectively, the outgoing and incoming characteristics associated with $u_\pm$. We denote by $\errfn(z):=\frac{1}{\sqrt{\pi}}\int_{-\infty}^z\rme^{-\xi^2}\,\rmd\xi$ the error function.

\begin{Theorem}\label{greenbounds}
Under the assumptions of Theorem~\ref{nonlin}, the Green's distribution $\mathcal{G}(x,t;y,s)$ associated with the linearized system \eqref{linearov} can be written as $\mathcal{G}=\mathcal{E}_1+\mathcal{E}_2+\tilde{\mathcal{G}}$ so that the following is true: There are positive constants $\eta$, $C$, and $M$ so that, for $y\leq0$ and $t\geq s$, we have
\begin{eqnarray}\label{E1}
\mathcal{E}_1(x,t;y,s) & = & \bar{u}_x(x,t) \pi_1(y,s,t) \\ \label{e_1}
\pi_1(y,s,t) & = & \sum_{a_\mathrm{in}^-} \left(\errfn\left(\frac{y+a_\mathrm{in}^-(t-s)}{\sqrt{4(t-s+1)}}\right)-\errfn\left(\frac{y-a_\mathrm{in}^-(t-s)}{\sqrt{4(t-s+1)}}\right)\right)l_{1,\mathrm{in}}^{-}(y,s)^T \\ \label{E2}
\mathcal{E}_2(x,t;y,s) & = & \bar{u}_t(x,t) \pi_2(y,s,t) \\ \label{e_2}
\pi_2(y,s,t) & = & \sum_{a_\mathrm{in}^-} \left(\errfn\left(\frac{y+a_\mathrm{in}^-(t-s)}{\sqrt{4(t-s+1)}}\right)-\errfn\left(\frac{y-a_\mathrm{in}^-(t-s)}{\sqrt{4(t-s+1)}}\right)\right)l_{2,\mathrm{in}}^{-}(y,s)^T
\end{eqnarray}
for appropriate functions $l_{j,\mathrm{in}}(y,s)$ that are $2\pi$-periodic in $s$ and satisfy
\begin{equation}\label{ljkbounds}
|\partial_s^\alpha l_{1,\mathrm{in}}^-(y,s)| + |\partial_s^\alpha l_{2,\mathrm{in}}^-(y,s)| \leq C, \qquad
|\partial_s^\alpha \partial_y^\beta l_{1,\mathrm{in}}^-(y,s)| + |\partial_s^\alpha \partial_y^\beta l_{2,\mathrm{in}}^-(y,s)| \leq C \rme^{-\eta|y|}
\end{equation}
for $0\leq|\alpha|\leq1$ and $1\leq\beta\leq2$, and
\begin{eqnarray*}
|\partial_{y}^\alpha \tilde{\mathcal{G}}(x,t;y,s)| & \leq &
C\rme^{-\eta(|x-y|+(t-s))}
\\ & & \nonumber
+ C \left( (t-s)^{-\frac{|\alpha|}{2}}+|\alpha| \rme^{-\eta|y|} \right)
\left(
\sum_{a^-} (t-s)^{-\frac12} \rme^{-(x-y- a^-(t-s))^2/M(t-s)} \rme^{-\eta x^+}
\right. \\ & & \nonumber
+ \sum_{a_\mathrm{out}^-,\,a_\mathrm{in}^-} \chi_{\{ |a_\mathrm{out}^- (t-s)|\ge |y| \}}
(t-s)^{-\frac12} \rme^{-(x-a_\mathrm{in}^-((t-s)-|y/a_\mathrm{out}^-|))^2/M(t-s)} \rme^{-\eta x^+}
\\ & & \nonumber \left.
+ \sum_{a_\mathrm{in}^-,\, a_\mathrm{out}^+} \chi_{\{ |a_\mathrm{in}^- (t-s)|\ge |y| \}}
(t-s)^{-\frac12} \rme^{-(x-a_\mathrm{out}^+ ((t-s)-|y/a_\mathrm{in}^-|))^2/M(t-s)} \rme^{-\eta x^-} \right)
\end{eqnarray*}
for $0\leq\alpha\leq1$. In the above, $x^\pm:=\max\{0,\pm x\}\geq0$ denotes the positive/negative part of $x$, and $\chi_J(x)$ is the indicator function of the interval $J$. Symmetric bounds hold for $y\geq0$ and $t\geq s$. Furthermore, if we replace $\bar{u}(x,t)$ with $\bar{u}(x-q_*,t-\tau_*)$, the estimates above remain true  uniformly for $(q_*,\tau_*)$ in any compact set.
\end{Theorem}

The terms $\mathcal{E}_1$ and $\mathcal{E}_2$ correspond to the translational eigenmodes in space and time, respectively. The error functions in (\ref{e_1}) and (\ref{e_2}) record the effects of information that gets transported along the incoming characteristics, which determine the ultimate space and time translate of the wave to which the perturbed solution converges. The term $\tilde{\mathcal{G}}$ encapsulates higher order terms arising from the parts of the spectrum to the left of the imaginary axis, including the continuous spectrum near zero, and from higher order terms associated with the translational eigenmodes. Overall, the description of the Green's distribution in Theorem~\ref{greenbounds} is exactly analogous to that in \cite{HowardZumbrun06} for the time-independent case, with the addition of the new term $\mathcal{E}_2$. Thus, in Theorem~\ref{greenbounds}, we have effectively performed a ``time-asymptotic conjugation'' to the time-independent case, analogous to the usual Floquet conjugation for finite-dimensional differential equations, but carried out only on the low-frequency modes important for time-asymptotic behavior.

We remark that we have made no effort to state minimal regularity assumptions;  instead, we chose to maintain an argument structure that generalizes easily to the case of quasilinear and partially parabolic viscosity $(B(u)u_x)_x$ as treated in \cite{RaoofiZumbrun07,TexierZumbrun06_gas}. In the constant-viscosity case treated here, $H^3$ could be replaced by $L^\infty$ in Theorem~\ref{nonlin}, and the short-time theory of \S\ref{auxests} and \S\ref{proof}, which is based on energy estimates, could be replaced by simpler $L^\infty$ theory based on integral equations and Picard iteration as in \cite[\S5]{HowardZumbrun06}. We expect that our results will extend to the quasilinear, partially parabolic case as in the closely related analyses of \cite{RaoofiZumbrun07,TexierZumbrun06_gas}.

We have focused here on pure Lax shocks, which are the most common type of viscous shocks, and the only type arising in standard gas dynamics. Other types of stationary viscous shocks include undercompressive, overcompressive, and mixed-type profiles: pure over- and undercompressive profiles arise in magnetohydrodynamics (MHD) and phase-transitional models, but mixed under--overcompressive profiles are also possible, as described in \cite{LiuZumbrun95, ZumbrunHoward98}. The nonlinear stability of stationary non-Lax viscous shocks has recently been addressed in \cite{HowardZumbrun06}. One key difference is that, in the pure Lax and the undercompressive case, the neutral eigenvalues result only from space translates of the profile $\bar{u}$, whereas in other cases they also involve deformations of $\bar{u}$; see \cite{ZumbrunHoward98, MasciaZumbrun03, MasciaZumbrun04_largeamp, Zumbrun01, Zumbrun04, Zumbrun04_lecturenotes} for further discussion. Our results for the case of time-periodic viscous shocks carry over to nonclassical over-, under-, and mixed over--undercompressive waves, with virtually no changes in the analysis, provided we substitute for the stability condition (S2) the condition that there exist $\ell$ eigenmodes at $\sigma=0$, where $\ell$ is the dimension of the manifold of time-periodic travelling-wave connections between $u_-$ and $u_+$ and for (S3)--(S4) the more fundamental condition
\begin{equation*}
\det\left[\int_\mathbb{R}\int_0^{2\pi}\langle \psi_j(x,t),\phi_k(x,t)\rangle\,\rmd x\,\rmd t\right]\ne 0
\end{equation*}
for bases $\psi_j$ and $\phi_k$ of left and right genuine eigenmodes of $\sigma=0$, where $j,k=1,\dots,\ell$. With this change, the analysis goes through essentially unchanged; see \cite{HowardZumbrun06} for the time-independent case.

The plan for the remainder of the paper is as follows. In \S\ref{Explanation_of_Method}, we give a detailed explanation of the method by which we will prove the main results of the paper, focusing on the intuition behind the ideas, rather than the technical details. 
Section~\ref{S:greens_functions} contains the linear theory we develop for equation (\ref{linearov}) using spatial dynamics, including the contour integral representation of the linear evolution, which provides the framework for the decomposition of the Green's distribution in \S\ref{S:res_decomp_est}. Using this decomposition, Theorems~\ref{stabcrit}  and~\ref{greenbounds} are proved in \S\ref{lin}, while \S\ref{nonlinsec} contains the nonlinear analysis and the proof of Theorem~\ref{nonlin}. We end this paper in \S\ref{s:conclusions} with a brief summary and a discussion of open problems.


\section{Outline of the method}\label{Explanation_of_Method}

We now give a detailed outline of the method used to obtain the results in \S\ref{introduction}, focusing on the intuition and main ideas involved. Rigorous justification can be found in \S\ref{S:greens_functions}-\S\ref{nonlinsec} below.


\subsection{Nonlinear stability of stationary shocks}

We first recall the techniques that can be used to study the nonlinear stability of stationary Lax shocks $\bar{u}(x)$ of (\ref{main}). Taking a solution of the form $u(x,t)=\bar{u}(x)+v(x,t)$, the resulting equation for the perturbation $v$ is
\[
v_t = \mathcal{L}v + Q(v,v_x)_x, \qquad
v(\cdot,0) = v_0,
\]
where
\[
\mathcal{L}v = v_{xx} - [f_u(\bar{u}(x))v]_x, \qquad
Q(v,v_x) = -f(\bar{u}+v)+f(\bar{u})+f_u(\bar{u})v.
\]
One way to understand the linearized equation
\[
v_t = \mathcal{L}v
\]
is to take its Laplace transform and obtain
\[
\lambda \hat{v} - v_0 = \mathcal{L}\hat{v}.
\]
This equation can be solved using the resolvent operator, and taking the inverse Laplace transform of the result leads to the standard representation
\[
\rme^{t \mathcal{L}} = \frac{1}{2\pi\rmi}
\int_\Gamma \rme^{t\lambda} [\lambda-\mathcal{L}]^{-1}\,\rmd\lambda
\]
of the linear semigroup, where $\Gamma$ is a curve in the complex plane that lies entirely in the resolvent set of $\mathcal{L}$. For viscous shocks, the essential spectrum of $\mathcal{L}$ contains the origin. Hence, we cannot easily derive good decay estimates for the linear semigroup, because we cannot move the contour $\Gamma$ into the origin, where analyticity of the resolvent $[\lambda-\mathcal{L}]^{-1}$ breaks down. Instead, we exploit that $\mathcal{L}$ is a differential operator and construct its Green's distribution. The Green's distribution $\mathcal{G}(x,t,y)$ is given by the semigroup acting on the Dirac delta function centered at $y$:
\begin{equation}\label{e:gd}
\mathcal{G}(x,t,y)
= \left(\rme^{t \mathcal{L}}\delta(\cdot-y)\right)(x)
= \frac{1}{2\pi\rmi}\int_\Gamma \rme^{t\lambda} \left( [\lambda-\mathcal{L}]^{-1} \delta(\cdot-y) \right)(x)\,\rmd\lambda
=: \frac{1}{2\pi\rmi}\int_\Gamma \rme^{t\lambda} G(x,y,\lambda)\,\rmd\lambda.
\end{equation}
Using variation of constants, solutions to the nonlinear equation can then be written as
\begin{equation}\label{e:vc}
v(x,t) = \int_{\mathbb{R}} \mathcal{G}(x,t,y) v_0(y) \,\rmd y - \int_0^t \int_{\mathbb{R}} \mathcal{G}_y(x,t-s,y) Q(v,v_x)(y,s)\,\rmd y \,\rmd s.
\end{equation}
To derive decay estimates for the Green's distribution, and hence for the integrals in (\ref{e:vc}), we can now deform the contour $\Gamma$ in the rightmost integral in (\ref{e:gd}) pointwise for each $(x,y)$. The key is that the resolvent kernel $G(x,y,\lambda)$ satisfies the ordinary differential equation
\begin{equation}\label{E:res_kern_eqn}
\lambda G - \delta(x-y) = G_{xx} - [f_u(\bar{u}(x))G]_x,
\end{equation}
so that ODE techniques such as the Gap Lemma \cite{GardnerZumbrun98,KapitulaSandstede98} can be used to extend the resolvent kernel meromorphically across $\lambda=0$ for each fixed $(x,y)$. It is then possible, see \cite{HowardZumbrun06, ZumbrunHoward98}, to move the contour $\Gamma$ for each $(x,y)$ to extract the leading order behavior of the Green's distribution $\mathcal{G}(x,t,y)$.

We now illustrate the outcome of extending the resolvent kernel meromorphically across the origin $\lambda=0$. First, $\lambda=0$ is a simple embedded eigenvalue with the exponentially localized eigenfunction $\bar{u}_x$. The associated adjoint eigenfunction is a constant $\psi_1\in\mathbb{R}^N$. Eigenvalues correspond to poles of the resolvent, and we therefore expect to obtain the term $\lambda^{-1}\bar{u}_x(x)\langle\psi_1,\cdot\rangle$, which involves the spectral projection, when we extend the resolvent kernel across the origin. A second contribution should appear due to the essential spectrum: To determine its contribution, we consider the constant-coefficient equation
\begin{equation}\label{e:rk}
\lambda G - G_{xx} + A(x) G_x = \delta(x-y), \qquad
A(x) = \left\{\begin{array}{lcl}
	f_u(u_-) & \quad & x<0 \\
	f_u(u_+) & \quad & x>0 \\
\end{array} \right.
\end{equation}
for $y\leq0$, say, where we neglected the shock profile, which we already accounted for through the spectral projection. To construct bounded solutions of this equation, we need to find the roots of the  characteristic equation $\det(\lambda-\nu^2+f_u(u_\pm)\nu)=0$, which are given by $\nu\approx a_j^\pm$ and $\nu\approx-\lambda/a_j^\pm$ for $\lambda$ near zero. The latter roots are the dangerous ones as they give small exponential rates for the associated solutions of (\ref{e:rk}), and thus they will determine the next term, after the spectral projection, in the expansion of the resolvent kernel. In addition, we need to select those roots that give exponential decay of the resolvent kernel in $(x,y)$ when $\lambda>0$. Using this information and the notation $a^- = \{ a^-_\mathrm{in}, a^-_\mathrm{out}\}$, the expansion
\[
G(x,y,\lambda) \approx
\frac{1}{\lambda} \bar{u}_x(x)\langle\psi_1,\cdot\rangle
+ \left\{\begin{array}{lclcl}\displaystyle
\sum_{a_\mathrm{out}^-,\,a^-} c_1 \rme^{-\lambda x/a_\mathrm{out}^+ + \lambda y/a^-}
& \quad & \mbox{ for } x\leq y\leq0 & \quad & (1) \\ \displaystyle
\sum_{a_\mathrm{in}^-,\,a^-} c_2 \rme^{-\lambda x/a^- + \lambda y/a_\mathrm{in}^-} 
& \quad & \mbox{ for } y\leq x\leq0 & \quad & (2) \\ \displaystyle
\sum_{a_\mathrm{in}^-,\,a_\mathrm{out}^+} c_3 \rme^{-\lambda x/a_\mathrm{out}^+ + \lambda y/a_\mathrm{in}^-}
& \quad & \mbox{ for } y\leq0\leq x & \quad & (3)
\end{array}\right.
\]
for the resolvent kernel can be derived when $y\leq0$ for appropriate coefficients $c_j$ that depend on the summation indices. The different cases in the equation above account for the transport of the perturbation along different characteristics, which we will discuss in more detail below. The inverse Laplace transform of the resolvent kernel is given by
\begin{eqnarray*}
\mathcal{G}(x,t,y) & \approx &
\bar{u}_x(x) \langle\psi_1,\cdot\rangle
+ \sum_{a_\mathrm{out}^-} \frac{c_1}{\sqrt{4\pi t}} \rme^{-\frac{(x-y-a_\mathrm{out}^-t)^2}{4t}}
+ \sum_{a_\mathrm{in}^-,\,a_\mathrm{out}^-} \frac{c_2}{\sqrt{4\pi t}}
	\rme^{-\frac{(x-a_\mathrm{out}^-(t-|y/a_\mathrm{in}^-|))^2}{4t}}
\\ &&
+ \sum_{a_\mathrm{in}^-,a_\mathrm{out}^+} \frac{c_3}{\sqrt{4\pi t}}
	\rme^{-\frac{(x-a_\mathrm{out}^+(t-|y/a_\mathrm{in}^-|))^2}{4t}}.
\end{eqnarray*}
Thus, the spectral projection onto the eigenfunction $\bar{u}_x$ appears, while the essential spectrum leads to Gaussians that move along the characteristics as they decay. It turns out that it is advantageous to use the term
\[
\frac{1}{\lambda} \bar{u}_x(x) \sum_{a_\mathrm{in}^-} \rme^{-\lambda y/a_\mathrm{in}^-} \langle c_\mathrm{in}^-,\cdot\rangle, \qquad y\leq0
\]
in place of $\bar{u}_x \langle \psi_1, \cdot \rangle$, which uses an expansion of the adjoint eigenfunction in terms of the weak spatial eigenvalues $\nu=-\lambda/a_\mathrm{in}^\pm$. This leads to the decomposition $\mathcal{G}(x,t,y)=\mathcal{E}(x,t,y)+\tilde{\mathcal{G}}(x,t,y)$ of the Green's distribution with
\begin{eqnarray*}
\mathcal{E}(x,t,y) & = & \bar{u}_x(x) \pi(y,t) \\
\pi(y,t) & \approx &
\sum_{a^-_\mathrm{in}}
\left( \errfn\left(\frac{y+a^-_\mathrm{in} t}{\sqrt{4(t+1)}}\right)
- \errfn\left(\frac{y-a_\mathrm{in}^- t}{\sqrt{4(t+1)}}\right) \right)
\langle c_\mathrm{in}^-,\cdot\rangle \\
\tilde{\mathcal{G}}(x,t,y) & \approx &
\sum_{a^-_\mathrm{out}} \frac{c_1}{\sqrt{4 \pi t}}
\rme^{-\frac{(x-y-a_\mathrm{out}^- t)^2}{4t}}
+ \sum_{a_\mathrm{in}^-,\,a_\mathrm{out}^-} \frac{c_2}{\sqrt{4 \pi t}}
\rme^{-\frac{(x-a_\mathrm{out}^-(t - |y/a^-_\mathrm{in}|))^2}{4t}}
+ \sum_{a^-_\mathrm{in},\,a^+_\mathrm{out}} \frac{c_3}{\sqrt{4 \pi t}}
\rme^{-\frac{(x-a_\mathrm{out}^+(t - |y/a^-_\mathrm{in}|))^2}{4t}}
\end{eqnarray*}
for $y\leq0$, where we again use the notation $a^\pm_\mathrm{out}$ and $a^\pm_\mathrm{in}$ to denote outgoing and incoming characteristics to aid in the following intuitive explanation of the above representation.

\begin{figure}
\centering\includegraphics{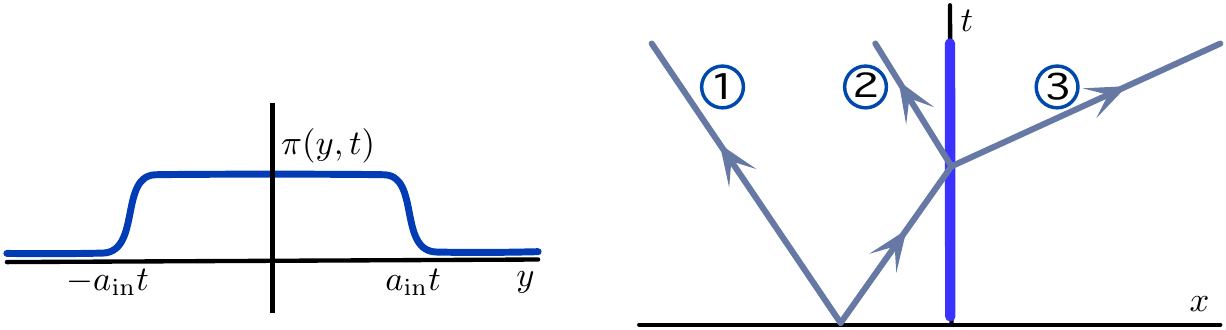}
\caption{On the left, a sketch of the expansion $\pi(y,t)$ of the adjoint eigenfunction is shown. On the right, the correspondence between the three terms in $\tilde{\mathcal{G}}$ and transport along the characteristics is illustrated.}
\label{F:greens_diagram}
\end{figure}

The term $\mathcal{E}$ consists of the eigenfunction $\bar{u}_x$ and the inverse Laplace transform $\pi(y,t)$ of the expansion of its adjoint eigenfunction. Hence, $\mathcal{E}$ is an expansion of the spectral projection. The sum in $\pi(y,t)$ is taken over the incoming directions, and it therefore tracks the initial data and nonlinear interactions resulting from the perturbation as they move in towards the shock and records their effect on the shock location. A sketch of the function $\pi$ is given in Figure~\ref{F:greens_diagram}: notice that $\pi(y,t)\to1$ pointwise as $t\to\infty$, and so $\bar{u}_x \pi \to \bar{u}_x\langle \psi_1, \cdot \rangle$. The three pieces of $\tilde{\mathcal{G}}$ record, respectively, how the perturbation moves along the different characteristics that transport it in three distinct ways, as illustrated in Figure~\ref{F:greens_diagram}: outwards away from the shock (1); inwards towards the shock and then reflected back out again (2); and in towards the shock, through it, and outwards on the other side of the shock (3).

We can now outline the nonlinear stability argument, which uses the above decomposition in the integral equation (\ref{e:vc}). More precisely, to track the shock location and remove the neutral direction $\bar{u}_x$ along the shock translates from (\ref{e:vc}), we write solutions as
\[
u(x+q(t),t) = \bar{u}(x) + v(x,t).
\]
Exploiting the decomposition of the Green's distribution, we define the phase shift $q(t)$ implicitly using
\begin{equation}\label{nl:q}
q(t) =
- \int_\mathbb{R} \pi(y,t)v_0(y) \,\rmd y
+ \int_0^t \int_\mathbb{R} \pi_y(y,t-s) \left[ Q(v, v_x) + \dot{q}v \right](y,s) \,\rmd y\,\rmd s
\end{equation}
and find after some algebra that the perturbation $v$ satisfies the integral equation
\begin{equation}\label{nl:v}
v(x,t) =
\int_\mathbb{R} \tilde{\mathcal{G}}(x,t;y) v_0(y) \,\rmd y
- \int_0^t \int_\mathbb{R} \tilde{\mathcal{G}}_y(x,t-s;y)
\left[ Q(v, v_x) + \dot{q}v \right](y,s) \,\rmd y \,\rmd s.
\end{equation}
Note that the evolution of $v$ is now governed only by the decaying part $\tilde{\mathcal{G}}$ of the Green's distribution. In \cite{HowardZumbrun06, ZumbrunHoward98}, it was shown that pointwise bounds, which result from the above formulas for $\mathcal{E}$ and $\tilde{\mathcal{G}}$, can now be used to establish existence of solutions to these integral equations and to prove that $q(t)$ converges to a limit $q_*$, while $v$ decays algebraically to zero.


\subsection{Nonlinear stability of time-periodic shocks}

Given the success of pointwise estimates in establishing stability of stationary viscous shocks, we would like to use the same technique in the case of time-periodic shocks. For this approach to work, we must show that the resolvent kernel and the Green's distribution in equation (\ref{e:gd}) all have well-defined counterparts for a time-periodic linear operator.

For time-periodic linear operators, the appropriate notion of the spectrum is the Floquet spectrum, and Floquet exponents $\sigma$ and the associated Floquet eigenfunctions $u(x,t)$ are found as solutions to the linearized equation
\begin{equation}\label{e:fe}
\sigma u + u_t = u_{xx} - [f_u(\bar{u}(x,t))u]_x,
\end{equation}
where any eigenfunction must be localized in space and satisfy $u(x,t+2\pi)=u(x,t)$ for all $t$. Due to the nonuniqueness of Floquet exponents, it suffices to consider only $\sigma$ with $-\frac12<\Im\sigma\leq\frac12$ (see \S\ref{s:g*} below). Based on (\ref{e:fe}) and comparing with (\ref{E:res_kern_eqn}), the resolvent kernel in the time-periodic setting should satisfy
\begin{equation}\label{e:fe_delta}
\sigma u + u_t - \delta(x-y)\delta(t-s) = u_{xx} - [f_u(\bar{u}(x,t))u]_x,
\end{equation}
where the additional temporal delta function sets the initial time $t=s$. In order to solve (\ref{e:fe_delta}), we first focus on (\ref{e:fe}) and write it as the first-order system
\begin{equation*}
U_x =
\begin{pmatrix}
0 & 1 \\ \partial_t+\sigma+f_{uu}(\bar{u})[\bar{u}_x,\cdot] & f_u(\bar{u})
\end{pmatrix} U
=: \mathcal{A}(x,\sigma) U
\end{equation*}
in the evolution variable $x$, posed on the space $H^{m+\frac12}(S^1)\times H^m(S^1)$ for some $m\geq0$. When this spatial dynamical system has an exponential dichotomy, given by solution operators $\Phi^\mathrm{s}(x,y,\sigma)$ for $x\geq y$ and $\Phi^\mathrm{u}(x,y,\sigma)$ for $x\leq y$ with $\Phi^\mathrm{s}(x,x,\sigma)+\Phi^\mathrm{u}(x,x,\sigma)=\mathrm{id}$, then the inhomogeneous equation
\begin{equation}\label{e:sd_H}
U_x = \mathcal{A}(x,\sigma)U + \mathcal{H}(x)
\end{equation}
can be solved uniquely, and the solution is given by
\begin{equation}\label{e:vc1}
U(x) = \int_{-\infty}^x \Phi^\mathrm{s}(x,z,\sigma) \mathcal{H}(z)\,\rmd z - \int_x^\infty \Phi^\mathrm{u}(x,z,\sigma) \mathcal{H}(z)\,\rmd z.
\end{equation}
Thus, solving equation (\ref{e:fe_delta}) is equivalent to using $\mathcal{H}(x)=(0,-\delta(x-y)\delta(\cdot-s))^T$, which leads to
\[
G(x,y,\sigma,s) = \left\{\begin{array}{rcl}
-P_1 \Phi^\mathrm{s}(x,y,\sigma)
\left(\begin{array}{c}0\\ \delta(\cdot-s)\end{array}\right) & \quad & x>y \\
 P_1 \Phi^\mathrm{u}(x,y,\sigma)
\left(\begin{array}{c}0\\ \delta(\cdot-s)\end{array}\right) & \quad & y>x
\end{array}\right.
\]
as the resolvent kernel of (\ref{e:fe}), where $P_1$ projects onto the first component. Taking the inverse Laplace transform, we find the Green's distribution via
\[
\mathcal{G}(x,t,y,s) = \frac{1}{2\pi\rmi} \int_{-\frac\rmi2}^{\frac\rmi2} \rme^{\sigma t} [G(x,y,\sigma,s)](t) \,\rmd\sigma.
\]
There are various issues that need to be addressed to make this approach work. Foremost among these issues is the regularity of the resolvent kernel $G(x,y,\sigma,s)$, since it is obtained by taking $[\mathcal{H}(x)](t)=(0,-\delta(x-y)\delta(t-s))^T$, for which we cannot solve (\ref{e:sd_H}) in $H^{m+\frac12}\times H^m$.


\section{Construction of the Green's distribution}\label{S:greens_functions}

Consider the linearization
\begin{equation}\label{e:lin}
u_t = u_{xx} - [f_u(\bar{u}(x,t))u]_x
\end{equation}
about the shock profile $\bar{u}(x,t)$. We say that $\mathcal{G}(x,t;y,s)$ is the Green's distribution of (\ref{e:lin}) if, for each given $u_0\in L^1(\mathbb{R})$, the function
\[
u(x,t) = \int_\mathbb{R} \mathcal{G}(x,t;y,s) u_0(y)\,\rmd y
\]
is a classical solution of (\ref{e:lin}) for $t>s$, and we have $u(x,t)\to u_0(x)$ as $t\searrow s$ for almost every $x\in\mathbb{R}$.

In this section, we shall show that (\ref{e:lin}) has a Green's distribution $\mathcal{G}(x,t;y,s)$. However, knowing its mere existence is not sufficient: to prove linear or nonlinear stability of time-periodic viscous shocks, we need to establish appropriate pointwise bounds for the Green's distribution. Thus, we shall construct the Green's distribution in a way that allows us to derive such bounds. Our strategy for finding the Green's distribution is as follows. Starting with the Green's function $\mathcal{G}_0(x,t)$ of the damped heat equation, which satisfies
\[
[\partial_t-\partial_x^2+1] \mathcal{G}_0 = 0, \qquad
\mathcal{G}_0|_{t=s} = \delta(x-y),
\]
we will iteratively construct a sequence $\mathcal{G}_j$ of functions that satisfy
\[
[\partial_t-\partial_x^2+1] \mathcal{G}_j = [1-\partial_x \cdot f_u(\bar{u})] \mathcal{G}_{j-1}, \qquad
\mathcal{G}_j|_{t=s}=0.
\]
We shall see that $\mathcal{G}_j$ will be become more regular as $j$ increases and, in addition, the difference $\mathcal{G}_*$ of the sum $\check{\mathcal{G}}=\sum_{j=0}^\ell\mathcal{G}_j$ of these functions and the desired Green's distribution $\mathcal{G}$ will become smoother as well. For a sufficiently large $\ell$, we can then construct the difference $\mathcal{G}_*$ as the inverse Laplace transform of a function $G_*$, which, in turn, is given via the expression (\ref{e:vc1}) for some regular $\mathcal{H}$. Pointwise bounds for $\check{\mathcal{G}}$ can now be derived immediately, upon exploiting their explicit construction. As we shall see in \S\ref{S:res_decomp_est}, we can also obtain pointwise bounds for $G_*$ by using spatial dynamics. The following theorem summarizes the existence result for $\mathcal{G}$ and the pointwise bounds for $\check{\mathcal{G}}$.

\begin{Theorem}\label{t:gd}
Equation (\ref{e:lin}) has a Green's distribution $\mathcal{G}(x,t;y,s)$, which lies in $C^2_x\cap C^1_t$ for $t>s$ and is bounded uniformly in $y$ for $x\in\mathbb{R}$ and $t>s$. The Green's distribution $\mathcal{G}$ can be written as $\mathcal{G}=\check{\mathcal{G}}+\mathcal{G}_*$ so that the following is true: For each $a\in\mathbb{R}$, there is a constant $M\geq1$ so that the function $\check{\mathcal{G}}$ obeys the pointwise estimate
\begin{equation}\label{est:gj}
\left|\partial_y^\alpha \check{\mathcal{G}}(x,t;y,s)\right| \leq
C (t-s)^{\frac{-1-|\alpha|}{2}} \rme^{-\frac{|x-y-a(t-s)|^2}{M(t-s)}},\qquad
x,y\in\mathbb{R},\quad
t>s
\end{equation}
for $0\leq\alpha\leq2$. Furthermore, the function $\mathcal{G}_*$ is given as the contour integral
\begin{equation}\label{n:ci}
\mathcal{G}_*(x,t;y,s) = \frac{1}{2\pi\rmi} \int_{\mu-\frac\rmi2}^{\mu+\frac\rmi2} \rme^{\sigma t} [G_*(x,y,\sigma;s)](t) \,\rmd\sigma, \qquad \mbox{for fixed }\mu>0
\end{equation}
for a function $G_*$ that is defined pointwise in $(x,t;y,s)$ and is analytic in $\sigma$ for $\sigma$ to the right of the Floquet spectrum $\Sigma$.
\end{Theorem}

In particular, the estimate (\ref{est:gj}) shows that, in Theorem~\ref{greenbounds}, the term $\check{\mathcal{G}}$ can be subsumed into the remainder term $\tilde{\mathcal{G}}$. Thus, once the preceding result is proved, it remains to establish pointwise bounds for $\mathcal{G}_*$ to complete the proof of Theorem~\ref{greenbounds}: this will be accomplished in \S\ref{S:res_decomp_est}-\S\ref{lin}. In the remainder of this section, we will prove Theorem~\ref{t:gd} and establish additional properties of $\mathcal{G}_j$ and $\mathcal{G}_*$.


\subsection{Construction of $\mathcal{G}_j$}\label{s:gj}

Let $\mathcal{G}_0(x,t;y,s)$ be the Green's distribution
\begin{equation}\label{n:g0}
\mathcal{G}_0(x,t;y,s) = \frac{1}{\sqrt{4\pi(t-s)}} \rme^{-\frac{(x-y)^2}{4(t-s)}-(t-s)}
\end{equation}
associated with the heat equation
\[
[\partial_t-\partial_x^2+1] \mathcal{G}_0(x,t;y,s) = 0, \qquad
\mathcal{G}_0(x,s;y,s) = \delta(x-y),
\]
where $x,y\in\mathbb{R}$ and $t>s$. We use $\mathcal{G}_0$ to define the functions $\mathcal{G}_j(x,t;y,s)$ recursively for $j\geq1$ as solutions of
\[
[\partial_t-\partial_x^2+1] \mathcal{G}_j = [1-\partial_x \cdot f_u(\bar{u}(x,t))] \mathcal{G}_{j-1}, \qquad
\mathcal{G}_j(x,s;y,s) = 0.
\]
Note that these functions are, by Duhamel's formula, given explicitly by
\begin{eqnarray}\label{n:gj}
\mathcal{G}_j(x,t;y,s) & = &
\int_s^t \int_\mathbb{R} \mathcal{G}_0(x,t;\tilde{y},\tilde{s})
[1-\partial_{\tilde{y}} \cdot f_u(\bar{u}(\tilde{y},\tilde{s}))]
\mathcal{G}_{j-1}(\tilde{y},\tilde{s};y,s)\,\rmd\tilde{y}\,\rmd\tilde{s}
\\ \nonumber & = &
\int_s^t \int_\mathbb{R} \mathcal{G}_0(x,t;\tilde{y},\tilde{s})
\mathcal{G}_{j-1}(\tilde{y},\tilde{s};y,s)\,\rmd\tilde{y}\,\rmd\tilde{s}
+ \int_s^t \int_\mathbb{R} \partial_{\tilde y}\mathcal{G}_0(x,t;\tilde{y},\tilde{s}) f_u(\bar{u}(\tilde{y},\tilde{s})) \mathcal{G}_{j-1}(\tilde{y},\tilde{s};y,s)\,\rmd\tilde{y}\,\rmd\tilde{s}.
\end{eqnarray}
We are interested in finding the Green's distribution $\mathcal{G}$ of $u_t=u_{xx}-[f_u(\bar{u})u]_x$, which satisfies
\[
[\partial_t-\partial_x^2+\partial_x \cdot f_u(\bar{u})] \mathcal{G} = 0, \qquad
\mathcal{G}(x,s;y,s) = \delta(x-y).
\]
Thus, seeking $\mathcal{G}$ in the form $\mathcal{G}=\mathcal{G}_*+\sum_{j=0}^\ell\mathcal{G}_j$, we find that $\mathcal{G}$ is the desired Green's distribution if and only if $\mathcal{G}_*(x,t;y,s)$ satisfies
\begin{equation}\label{n:g*}
[\partial_t-\partial_x^2+\partial_x \cdot f_u(\bar{u})] \mathcal{G}_* = [1-\partial_x \cdot f_u(\bar{u})] \mathcal{G}_\ell, \qquad
\mathcal{G}_*(x,s;y,s) = 0.
\end{equation}
We postpone the discussion of $\mathcal{G}_*$ to the next section and focus in the remainder of this section on pointwise bounds for the contributions $\mathcal{G}_j$.

\begin{Lemma}\label{parametrixlem}
Assume that $f\in C^k$, then, for each $j\geq0$, there is a constant $C\geq1$ such that
\begin{equation}\label{Gjbounds}
\left|\partial_{x,y}^\alpha\partial_{t,s}^\beta \mathcal{G}_j(x,t;y,s)\right| \leq C (t-s)^{\frac{j-1-|\alpha|}{2}-|\beta|} \rme^{-\frac{|x-y|^2}{C(t-s)}-\frac{t-s}{C}}
\end{equation}
for multi-indices $(\alpha,\beta)$ with $0\leq\max\{|\alpha_x|+|\beta_t|,|\alpha_y|+|\beta_s|\}\leq k-1$, and
\begin{equation}\label{2Gjbounds}
\left|(\partial_x+\partial_y)^\gamma (\partial_t+\partial_s)^\delta \partial_{x,y}^\alpha \partial_{t,s}^\beta \mathcal{G}_j(x,t;y,s)\right| \leq
C (t-s)^{\frac{j-1-|\alpha|}{2}-|\beta|} \rme^{-\frac{|x-y|^2}{C(t-s)}-\frac{t-s}{C}}
\end{equation}
for $0\leq\max\{|\alpha_x|+|\beta_t|,|\alpha_y|+|\beta_s|\}+|\gamma|+|\delta| \leq k-1$.
\end{Lemma}

\begin{Proof}
Equation (\ref{n:g0}) can be used to obtain \eqref{Gjbounds}-\eqref{2Gjbounds} for $j=0$ by direct computation. We can therefore proceed by induction to obtain estimates for $j\geq1$. The semigroup property for the damped heat equation $Cu_t=u_{xx}-u$ implies that
\[
\int_\mathbb{R} (t-\tilde{s})^{-\frac12} \rme^{-\frac{|x-\tilde{y}|^2}{C(t-\tilde{s})}-\frac{t-\tilde{s}}{C}}
(\tilde{s}-s)^{-\frac12} \rme^{-\frac{|\tilde{y}-y|^2}{C(\tilde{s}-s)}-\frac{\tilde{s}-s}{C}} \,\rmd\tilde{y} =
\tilde{C} (t-s)^{-\frac12} \rme^{-\frac{|x-y|^2}{C(t-s)}-\frac{t-s}{C}}.
\]
Using the induction hypothesis and \eqref{n:gj}, we can therefore estimate
\begin{eqnarray*}
|\mathcal{G}_j(x,t;y,s)| & \leq & 
\int_s^t \int_\mathbb{R} \Big( |\mathcal{G}_0(x,t;\tilde{y},\tilde{s})| |\mathcal{G}_{j-1}(\tilde{y},\tilde{s};y,s)| + |\partial_{\tilde y}\mathcal{G}_0(x,t;\tilde{y},\tilde{s})| |f_u(\bar{u}(\tilde{y},\tilde{s}))\mathcal{G}_{j-1}(\tilde{y},\tilde{s};y,s)| \Big)\,\rmd\tilde{y}\,\rmd\tilde{s}
\\ & \leq &
C (t-s)^{-\frac12} \rme^{-\frac{|x-y|^2}{C(t-s)}-\frac{t-s}{C}}
\int_s^t (\tilde{s}-s)^{\frac{j-1}{2}} \left( 1+(t-\tilde{s})^{-\frac{1}{2}} \right) \rmd\tilde{s}
\\ & \leq &
\tilde{C} (t-s)^{\frac{j-1}{2}} \rme^{-\frac{|x-y|^2}{C(t-s)}-\frac{t-s}{C}},
\end{eqnarray*}
which verifies the estimate \eqref{Gjbounds} for $\alpha=\beta=0$.

To obtain \eqref{2Gjbounds} for $\alpha=\beta=0$, we observe that $(\partial_x+\partial_y)\mathcal{G}_0(x,t;y,s)=(\partial_t+\partial_s)\mathcal{G}_0(x,t;y,s)=0$ since $\mathcal{G}_0$ depends on its arguments through $x-y$ and $t-s$ only. Thus, after integration by parts in $(\tilde{y},\tilde{s})$, we obtain
\begin{eqnarray*}
(\partial_x+\partial_y) \mathcal{G}_j(x,t;y,s) & = &
\int_s^t \int_\mathbb{R} \mathcal{G}_0(x,t;\tilde{y},\tilde{s}) (\partial_{\tilde{y}}+\partial_y) \mathcal{G}_{j-1}(\tilde{y},\tilde{s};y,s)\,\rmd\tilde{y}\,\rmd\tilde{s}
\\ &&
+ \int_s^t \int_\mathbb{R} \partial_{\tilde y} \mathcal{G}_0(x,t;\tilde{y},\tilde{s}) f_u(\bar{u}(\tilde{y},\tilde{s})) (\partial_{\tilde{y}}+\partial_y)
\mathcal{G}_{j-1}(\tilde{y},\tilde{s};y,s)\,\rmd\tilde{y}\,\rmd\tilde{s}
\\ &&
+ \int_s^t \int_\mathbb{R} \partial_{\tilde y} \mathcal{G}_0(x,t;\tilde{y},\tilde{s}) [\partial_{\tilde y}f_u(\bar{u}(\tilde{y},\tilde{s}))] \mathcal{G}_{j-1}(\tilde{y},\tilde{s};y,s)\,\rmd\tilde{y}\,\rmd\tilde{s},
\end{eqnarray*}
and similarly for $(\partial_t+\partial_s)$. Using the induction hypothesis, we can now estimate these integrals as above to obtain \eqref{2Gjbounds} for $\alpha=\beta=0$. Notice that we lose one degree of regularity for $f_u\in C^{k-1}$ for each application of $(\partial_x+\partial_y)$ or $(\partial_t+\partial_s)$, hence $|\gamma|+|\delta|$ in total.

Finally, using \eqref{2Gjbounds} to shift $(y,s)$ to $(\tilde{y},\tilde{s})$ derivatives where needed, and arguing by induction on $j$, $|\alpha|$,
$|\beta|$, and $|\gamma|$, $|\delta|$, we may estimate 
\begin{eqnarray*}
\lefteqn{ |\partial_{x,y}^\alpha \partial_{t,s}^\beta \mathcal{G}_j(x,t;y,s)| \leq }
\\ &&
\int_s^{\frac{s+t}{2}} \int_\mathbb{R} \left( \left|\partial_x^{|\alpha|}\partial_t^{|\beta|} \mathcal{G}_0(x,t;\tilde{y},\tilde{s})\right| \left|\mathcal{G}_{j-1} (\tilde{y},\tilde{s};y,s) \right|
+ \left|\partial_x^{|\alpha|}\partial_t^{|\beta|} \partial_{\tilde y}\mathcal{G}_0(x,t;\tilde{y},\tilde{s})\right| \left|f_u(\bar{u}(\tilde{y},\tilde{s}))\mathcal{G}_{j-1}(\tilde{y},\tilde{s};y,s)\right|\right) \rmd\tilde{y}\,\rmd\tilde{s}
\\ &&
+ \int_{\frac{s+t}{2}}^t \int_\mathbb{R} \left( \left|\mathcal{G}_0(x,t;\tilde{y},\tilde{s})\right| \left|\partial_y^{|\alpha|}\partial_s^{|\beta|} \mathcal{G}_{j-1}(\tilde{y},\tilde{s};y,s)\right|
+ \left|\partial_{\tilde y}\mathcal{G}_0(x,t;\tilde{y},\tilde{s})\right| \left|f_u(\bar{u}(\tilde{y},\tilde{s}))\partial_y^{|\alpha|}\partial_s^{|\beta|} \mathcal{G}_{j-1}(\tilde{y},\tilde{s};y,s)\right| \right)\,\rmd\tilde{y}\,\rmd\tilde{s}
\end{eqnarray*}
plus terms that contain lower-order derivatives of $\mathcal{G}_j$, $(\partial_x+\partial_y)\mathcal{G}_j$, or $(\partial_s+\partial_t)\mathcal{G}_j$, yielding
\begin{eqnarray*}
|\partial_{x,y}^\alpha \partial_{t,s}^\beta \mathcal{G}_j(x,t;y,s)| & \leq &
C (t-s)^{-\frac12} \rme^{-\frac{|x-y|^2}{C(t-s)}-\frac{t-s}{C}}
\int_s^{\frac{s+t}{2}} \int_\mathbb{R} (t-\tilde{s})^{\frac{-|\alpha|-2|\beta|}{2}} \left(1+ (t-\tilde{s})^{-\frac{1}{2}}\right) (s-\tilde s)^{\frac{j-1}{2}}
\rmd\tilde{s} \\ &&
+ C (t-s)^{-\frac12} \rme^{-\frac{|x-y|^2}{C(t-s)}-\frac{t-s}{C}}
\int_{\frac{s+t}{2}}^t \int_\mathbb{R} (\tilde{s}-s)^{\frac{j-1-|\alpha|-2|\beta|}{2}} \left(1+ (t-\tilde{s})^{-\frac{1}{2}}\right) \rmd\tilde{s}
\\ & \leq & 
\tilde{C} (t-s)^{\frac{j-1-|\alpha|-2|\beta|}{2}} \rme^{-\frac{|x-y|^2}{C(t-s)}-\frac{t-s}{C}}
\end{eqnarray*}
as claimed, which verifies \eqref{Gjbounds}. A similar argument yields \eqref{2Gjbounds}. Noting that each shift of $x,t$ or $y,s$ derivative costs one degree of regularity for $f_u\in C^{k-1}$, and that we must shift $|\alpha_y|+|\alpha_s|$ derivatives in the first integral and $|\alpha_x|+|\alpha_t|$ in the second, we obtain the stated range of indices.
\end{Proof}

The computations in the proof of Lemma~\ref{parametrixlem} may be recognized as parametrix estimates as in standard short-time parabolic theory \cite{Friedman,GilbargTrudinger,Ladyzhenskaya}.

The bounds for $\check{\mathcal{G}}=\sum_{j=0}^\ell\mathcal{G}_j$ asserted in Theorem~\ref{t:gd} are now a consequence of Lemma~\ref{parametrixlem}. Indeed, for each fixed $a\in\mathbb{R}$ and $C>0$, there is a constant $M\geq1$, chosen sufficiently large, so that
\[
\rme^{-\frac{(x-y)^2}{C(t-s)}-\frac{(t-s)}{C}} \leq M \rme^{-\frac{(x-y-a(t-s))^2}{M(t-s)}}
\]
for all $x,y\in\mathbb{R}$ and $t>s$. To complete the proof of Theorem~\ref{t:gd}, it therefore remains to verify the assertions about $\mathcal{G}_*$.


\subsection{Construction of $\mathcal{G}_*$}\label{s:g*}

In this section, we verify the assertions made in Theorem~\ref{t:gd} about the contribution $\mathcal{G}_*$ to the Green's distribution $\mathcal{G}$ and show that it is given by (\ref{n:ci}) for an appropriate function $G_*$, which we shall refer to as the resolvent kernel.

Recall that $\mathcal{G}_*$ needs to satisfy equation (\ref{n:g*}), given by
\begin{equation}\label{n:g*1}
[\partial_t-\partial_x^2+\partial_x f_u(\bar{u}(x,t))] \mathcal{G}_*(x,t;y,s) =
[1-\partial_x f_u(\bar{u}(x,t))] \mathcal{G}_\ell(x,t;y,s), \qquad
\mathcal{G}_*(x,s;y,s) = 0.
\end{equation}
We shall use the Laplace transform to solve (\ref{n:g*1}), since this approach will allow us to reformulate (\ref{n:g*1}) as a spatial dynamical system, which facilitates the verification of pointwise bounds.

It is convenient to shift the time-dependent coefficients, and we therefore consider the equivalent system
\begin{equation}\label{e:le}
[\partial_t-\partial_x^2+\partial_x f_u(\bar{u}(x,t+s))] \mathcal{G}_*(x,t) = g_\ell(x,t;y,s),\qquad
\mathcal{G}_*(x,0) = 0,
\end{equation}
where we omit the arguments $(y,s)$ in the notation for $\mathcal{G}_*(x,t)$ and use the notation
\[
g_\ell(x,t;y,s) := [1-\partial_x f_u(\bar{u}(x,t+s))] \mathcal{G}_\ell(x,t+s;y,s).
\]
The term $f_u(\bar{u}(x,t+s))$ is smooth and $2\pi$-periodic in $t$ and can therefore be represented by its Fourier series
\[
f_u(\bar{u}(x,t+s)) = \sum_{k\in\mathbb{Z}} f_k(x) \rme^{\rmi k(t+s)}.
\]
Recall that the Laplace transform is defined by
\begin{equation*}
\hat{v}(x,\lambda) = \int_0^\infty \rme^{-\lambda t} v(x,t) \,\rmd t 
\end{equation*}
for all $\lambda\in\mathbb{C}$ with $\Re\lambda>\mu$, where $\mu$ is chosen so that the integral is convergent. Taking the Laplace transform of equation (\ref{e:le}), we obtain
\begin{equation*}
\lambda \hat{\mathcal{G}}_*(x,\lambda) = \partial_x^2 \hat{\mathcal{G}}_*(x,\lambda)
- \partial_x \left( \sum_{k\in\mathbb{Z}} \rme^{\rmi ks} f_k(x)
\hat{\mathcal{G}}_*(x,\lambda-\rmi k) \right) + \hat{g}_\ell(x,\lambda;y,s).
\end{equation*}
The effect of the time-periodic coefficients is that the equations for different values of $\lambda$ couple. Note, however, that the equations for $\hat{\mathcal{G}}_*$ at $\lambda=\lambda_1$ and $\lambda=\lambda_2$ couple only if $\lambda_1-\lambda_2\in\rmi\mathbb{Z}$. To exploit this fact, we define
\begin{equation*}
\lambda = \sigma+\rmi n, \qquad
-\frac{1}{2}<\Im\sigma\leq\frac{1}{2}, \qquad
n\in\mathbb{Z},
\end{equation*}
which is motivated by \cite{CostinCostinLebowitz04} and similar to a Fourier--Bloch wave decomposition of periodic functions. For each $\sigma$, we can view the above equation as a system of infinitely many second-order ODEs in $x$. In particular, if we define $\hat{\mathcal{G}}_*^n(x,\sigma):=\hat{\mathcal{G}}_*(x,\sigma+\rmi n)$, and similarly for $\hat{g}_\ell$, we arrive at the system
\begin{equation}\label{E:lt_linearov_sigma}
(\sigma+\rmi n) \hat{\mathcal{G}}_*^n = \partial_x^2 \hat{\mathcal{G}}_*^n - \partial_x \left( \sum_{k\in\mathbb{Z}} \rme^{\rmi ks} f_k(x) \hat{\mathcal{G}}_*^{n-k}\right) + \hat{g}_\ell^n(x,\sigma;y,s).
\end{equation}
The following result, which we will prove in \S\ref{s:rk} below, implies the existence of solutions to (\ref{E:lt_linearov_sigma}). Recall the definition (\ref{e:floq}) of the Floquet spectrum $\Sigma$ of (\ref{e:le}).

\begin{Proposition}\label{p:rk}
Fix $\ell\geq3$. For each $\sigma$ to the right of the Floquet spectrum $\Sigma$, the system
\begin{equation}\label{e:sd1}
\sigma G_* + \partial_t G_* = \partial_x^2 G_* - [f_u(\bar{u}(x,t+s)) G_*]_x + \sum_{n\in\mathbb{Z}} \rme^{\rmi nt}\hat{g}_\ell^n(x,\sigma;y,s), \qquad
G_*(x,0;y,s,\sigma) = 0
\end{equation}
has a unique solution $G_*(x,t;y,s,\sigma)$ in $C^2_x\cap C^1_t$ that is $2\pi$-periodic in $t$. Furthermore, $G_*(x,t;y,s,\sigma)$ is analytic in $\sigma$ and lies in $C^2_y\cap C^1_s$ with respect to $(y,s)$.
\end{Proposition}

Writing $G_*(x,t;y,s,\sigma)$ as the Fourier series
\[
G_*(x,t;y,s,\sigma) = \sum_{n\in\mathbb{Z}} \rme^{\rmi nt} \hat{\mathcal{G}}_*^n(x,\sigma;y,s),
\]
we see that the Fourier coefficients $\hat{\mathcal{G}}_*^n(x,\sigma;y,s)$ satisfy (\ref{E:lt_linearov_sigma}). Furthermore, the inverse Laplace transform
\begin{equation}\label{E:int_greens}
\mathcal{G}_*(x,t;y,s)
= \frac{1}{2\pi\rmi} \int_{\mu-\rmi\infty}^{\mu+\rmi\infty} \rme^{\lambda t} \hat{\mathcal{G}}_*(x,\lambda) \,\rmd\lambda
= \frac{1}{2\pi\rmi} \int_{\mu-\frac\rmi2}^{\mu+\frac\rmi2} \rme^{\sigma t} G_*(x,t;y,s,\sigma) \,\rmd\sigma
\end{equation}
of $\hat{\mathcal{G}}_*(x,\lambda)$ is well defined for each fixed $\mu>0$ and satisfies equation (\ref{e:le}). This completes the proof of Theorem~\ref{t:gd}, subject to the proof of Proposition~\ref{p:rk}, which we will give in \S\ref{s:rk} below.

From independence of (\ref{E:int_greens}) with respect to $\mu>0$ together with analyticity of $G_*$ with respect to $\sigma$ on $\{\Re\sigma>0\}$, we may conclude from Cauchy's integral theorem that
\begin{equation}\label{keyper}
\rme^{(\sigma+\rmi/2)t} G_*(x,t;y,s,\sigma +\rmi/2) =
\rme^{(\sigma-\rmi/2)t} G_*(x,t;y,s,\sigma -\rmi/2)
\end{equation}
for all $t$. This is evident at a formal level, as is the more general property that $\rme^{\sigma t}G_*$ is periodic in $\sigma$ with period $\rmi$ or, equivalently, $\rme^{-\rmi k t}G_*(x,t;y,s,\sigma)=G_*(x,t;y,s,\sigma+\rmi k)$ for all $k\in\mathbb{Z}$, which follows since the left-hand side satisfies (\ref{e:sd1}) with $\sigma$ replaced by $\sigma-\rmi k$.


\subsection{Construction of the resolvent kernel $G_*$}\label{s:rk}

In this section, we prove Proposition~\ref{p:rk}. Throughout this section, $\sigma$ will lie in the set $\Omega$, which denotes the connected component of $\mathbb{C}\setminus\Sigma$ that contains $\sigma=\infty$. We need to construct a solution $G_*(x,t,\sigma)$ of
\[
\sigma G_* + \partial_t G_* = \partial_x^2 G_* - [f_u(\bar{u}(x,t+s)) G_*]_x + \tilde{g}_\ell(x,t,\sigma;y,s), \qquad
G_*(x,0,\sigma) = 0
\]
that is analytic in $\sigma\in\Omega$ and $2\pi$-periodic in $t$, where $\tilde{g}_\ell(x,t,\sigma;y,s)$ denotes the Fourier series of $\{\hat{g}_\ell^n(x,\sigma;y,s)\}$. We rewrite this equation as the spatial dynamical system
\begin{equation}\label{E:sd_evalue_inhom}
U_x = \mathcal{A}(x,\sigma) U + \Delta_\ell(x) :=
\begin{pmatrix} 0 & 1 \\ \partial_t+\sigma+f_{uu}(\bar{u})[\bar{u}_x,\cdot] & f_u(\bar{u}) \end{pmatrix} U + \begin{pmatrix} 0 \\ -\tilde{g}_\ell(x,t,\sigma;y,s) \end{pmatrix}
\end{equation}
in the evolution variable $x$, with $\bar{u}=\bar{u}(x,t+s)$, where $U(x)=(u,u_x)^T$ is a $2\pi$-periodic function in time for each fixed $x$.

To construct solutions of (\ref{E:sd_evalue_inhom}), we first focus on the associated homogeneous system
\begin{equation}\label{E:sd_evalue}
U_x = \mathcal{A}(x,\sigma) U
= \begin{pmatrix} 0 & 1 \\
\partial_t+\sigma+f_{uu}(\bar{u})[\bar{u}_x,\cdot] & f_u(\bar{u})
\end{pmatrix} U.
\end{equation}
We consider the spatial dynamical system (\ref{E:sd_evalue}) on the Hilbert space $Y_m=H^{m+\frac12}(S^1)\times H^{m}(S^1)$, where $m\geq0$ will be chosen later. It can be shown that the operator on the right-hand side of (\ref{E:sd_evalue}) is closed and densely defined with domain $Y_{m+\frac12}$; see \cite{SandstedeScheel01}. Equation (\ref{E:sd_evalue}) is ill-posed, in the sense that solutions to arbitrary initial data in $Y_m$ may not exist: indeed, the leading-order operator
\[
\begin{pmatrix} 0 & 1 \\ \partial_t  & 0 \end{pmatrix}
\]
has spectrum given by $\{\pm\sqrt{\rmi k}:\;k\in\mathbb{Z}\}$ and therefore does not generate a semigroup on $Y_m$. Nevertheless, equation (\ref{E:sd_evalue}) provides a useful framework for analyzing the PDE (\ref{e:le}), since it admits exponential dichotomies whose properties can be related to spectral properties of (\ref{e:le}):

\begin{Definition}[{\cite[\S2.1]{PeterhofSandstedeScheel97}}]\label{d:ed}
Let $J=\mathbb{R}^+$, $\mathbb{R}^-$ or $\mathbb{R}$. Equation (\ref{E:sd_evalue}) is said to have an exponential dichotomy on $J$ if there exist positive constants $K$ and $\kappa^\mathrm{s}<0<\kappa^\mathrm{u}$ and two strongly continuous families of bounded operators $\Phi^\mathrm{s}(x,z)$ and $\Phi^\mathrm{u}(x,z)$ on $Y_m$, defined respectively for $x\geq z$ and $x\leq z$, such that
\[
\sup_{x\geq z,\, x,z\in J} \rme^{-\kappa^\mathrm{s}(x-z)} \|\Phi^\mathrm{s}(x,z)\|_{L(Y_m)} + \sup_{x\leq z,\, x,z\in J} \rme^{-\kappa^\mathrm{u}(x-z)} \|\Phi^\mathrm{u}(x,z)\|_{L(Y_m)} \leq K,
\]
the operators $P^\mathrm{s}(x):=\Phi^\mathrm{s}(x,x)$ and $P^\mathrm{u}(x):=\Phi^\mathrm{u}(x,x)$ are complementary projections for all $x\in J$, and the functions $\Phi^\mathrm{s}(x,z)U_0$ and $\Phi^\mathrm{u}(x,z)U_0$ satisfy (\ref{E:sd_evalue}) for $x>z$ and $x<z$, respectively, with values in $Y_m$ for each fixed $U_0\in Y_m$.
\end{Definition}

It follows from \cite[Remark~2.5 and Theorem~2.6]{SandstedeScheel01} that (\ref{E:sd_evalue}) has an exponential dichotomy on $\mathbb{R}$ for each $\sigma\in\Omega$. Furthermore, \cite[Proof of Theorem~1]{PeterhofSandstedeScheel97} implies that the operators $\Phi^\mathrm{s}(x,z,\sigma)$ and $\Phi^\mathrm{u}(x,z,\sigma)$ are analytic in $\sigma$ for $\sigma\in\Omega$ as functions into $L(Y_m)$. As a consequence of \cite[\S6.1]{SandstedeScheel01}, we can then solve the inhomogeneous equation
\[
U_x = \mathcal{A}(x,\sigma) U + \mathcal{H}(x)
\]
uniquely for each $\mathcal{H}\in L^2(\mathbb{R},Y_m)$ via the variation-of-constants formula
\begin{equation}\label{e:sf}
U(x) = \int_{-\infty}^x \Phi^\mathrm{s}(x,z,\sigma) \mathcal{H}(z)\,\rmd z
+ \int_\infty^x \Phi^\mathrm{u}(x,z,\sigma) \mathcal{H}(z)\,\rmd z,
\end{equation}
and the solution satisfies
\[
U \in H^1(\mathbb{R},Y_m) \cap L^2(\mathbb{R},Y_{m+\frac12}).
\]
Our goal is to apply these results to equation (\ref{E:sd_evalue_inhom}),
\[
U_x = \mathcal{A}(x,\sigma) U + \Delta_\ell(x), \qquad
\Delta_\ell(x) := \begin{pmatrix} 0 \\ -\tilde{g}_\ell(x,t,\sigma;y,s) \end{pmatrix},
\]
for $\sigma\in\Omega$. To establish the regularity of the right-hand side $\Delta_\ell(x)$, we repeat the iterative construction of the components $\mathcal{G}_j$ from \S\ref{s:gj} for their Laplace--Fourier transforms, which is akin to the bootstrapping arguments carried out in \cite[\S5.3]{SandstedeScheel01} and \cite{CostinCostinLebowitz04}. Thus, we are led to consider the equation
\[
V_x = \mathcal{A}_0(\sigma) V, \qquad
\mathcal{A}_0(\sigma) = \begin{pmatrix} 0 & 1 \\ \partial_t+\sigma+1 & 0 \end{pmatrix},
\]
which corresponds to the heat equation $u_t=u_{xx}-u$ that we utilized in \S\ref{s:gj}. Writing this equation in terms of its Fourier modes $\hat{V}^n$, we obtain the system
\[
\partial_x \hat{V}^n = \begin{pmatrix} 0 & 1 \\ \rmi n+\sigma+1 & 0 \end{pmatrix} \hat{V}^n,
\]
which admits the exponential dichotomy
\begin{eqnarray}\label{e:ed0}
\Phi^{\mathrm{u},n}_0(x,y,\sigma) & = & \frac{1}{2\sqrt{\sigma+1+\rmi n}}
\begin{pmatrix} \sqrt{\sigma+1+\rmi n} & 1 \\
\sigma+1+\rmi n & \sqrt{\sigma+1+ \rmi n} \end{pmatrix}
\rme^{-\sqrt{\sigma+1+\rmi n}|x-y|} \\ \nonumber
\Phi^{\mathrm{s},n}_0(x,y,\sigma) & = & \frac{1}{2\sqrt{\sigma+1+\rmi n}} \begin{pmatrix} \sqrt{\sigma+1+\rmi n} & -1 \\
-(\sigma+1+\rmi n) & \sqrt{\sigma+1+\rmi n} \end{pmatrix}
\rme^{-\sqrt{\sigma+1+\rmi n}|x-y|}.
\end{eqnarray}
In line with the definition of $\mathcal{G}_0$ as the Green's function of $u_t=u_{xx}-u$, we consider the equation
\begin{equation}\label{E:G0_eqn}
\partial_x V_0 = \mathcal{A}_0(\sigma) V_0 + \Delta_0(x), \qquad
\Delta_0(x) = \begin{pmatrix} 0 \\ -\delta(x-y)\delta(t) \end{pmatrix}.
\end{equation}
Using (\ref{e:sf}) and (\ref{e:ed0}), the Fourier modes $\hat{V}_0^n$ of the solution $V_0(x)$ are then given by
\begin{eqnarray*}
\hat{V}^n_0(x,y,\sigma) & = &
\int_{-\infty}^x \Phi^{\mathrm{s},n}_0(x,z,\sigma)
\begin{pmatrix} 0 \\ -\delta(z-y) \end{pmatrix}  \,\rmd z
+ \int_\infty^x \Phi^{\mathrm{u},n}_0(x,z,\sigma)
\begin{pmatrix} 0 \\ -\delta(z-y) \end{pmatrix} \,\rmd z
\nonumber \\ & = &
\frac{1}{2\sqrt{\sigma+1+\rmi n}}
\begin{pmatrix} 1 \\ -\sgn(x-y) \sqrt{\sigma+1+\rmi n} \end{pmatrix}
\rme^{-\sqrt{\sigma+1+\rmi n}|x-y|}.
\end{eqnarray*}
The function $V_1$ that corresponds to $\mathcal{G}_1$ can be found by solving the equation
\[
\partial_x V_1 = \mathcal{A}_0(\sigma) V_1 + \mathcal{B}(x) V_0,
\]
where
\begin{equation}\label{e:B}
\mathcal{B}(x) := \mathcal{A}(x,\sigma)-\mathcal{A}_0(\sigma)
= \begin{pmatrix} 0 & 0 \\ f_{uu}(\bar{u})[\bar{u}_x,\cdot]-1 & f_u(\bar{u})\end{pmatrix}.
\end{equation}
The equation for $V_1$ can again be solved explicitly using the exponential dichotomy (\ref{e:ed0}) for the Fourier modes. To do so, let $\{a_k\}$ and $\{b_k\}$ denote the Fourier components of $f_u(\bar{u})$ and $f_{uu}(\bar{u})[\bar{u}_x,\cdot]$, respectively, (suppressing the dependence on $s$) and note that
\[
[\widehat{\mathcal{B}V_0}]^n(x,y,\sigma) =
\sum_{k\in\mathbb{Z}} \rme^{-\sqrt{\sigma+1+\rmi k}|x-y|}
\begin{pmatrix} 0 \\ \displaystyle \frac{b_{n-k}(x)}{2\sqrt{\sigma+1+\rmi k}} - \sgn(x-y) \frac{a_{n-k}(x)}{2} \end{pmatrix},
\]
where we set $\tilde{b}_0:=b_0-1$ and dropped the tilde. Applying the exponential dichotomy (\ref{e:ed0}), we obtain
\begin{eqnarray}\label{E:G1}
\lefteqn{\hat{V}^n_1(x,y,\sigma) \;=\;} \\ \nonumber &&
\int_{\mathbb{R}} \rme^{-\sqrt{\sigma+1+\rmi n}|x-z|}
\sum_{k\in\mathbb{Z}} \frac{\rme^{-\sqrt{\sigma+1+\rmi k}|z-y|}}{4} 
\begin{pmatrix} \displaystyle
- \frac{b_{n-k}(z)}{\sqrt{\sigma+1+\rmi n}\sqrt{\sigma+1+\rmi k}}
+ \frac{\sgn(x-y) a_{n-k}(z)}{\sqrt{\sigma +1+ \rmi n}} \\ \displaystyle
- \frac{b_{n-k}(z)}{\sqrt{\sigma+1+\rmi k}} + \sgn(x-y) a_{n-k}(z) \end{pmatrix}\rmd z.
\end{eqnarray}
We record the following estimate on $V_1$, which will be used below. By estimating the ``worst" term in $\hat{V}_1^n$, we obtain
\begin{eqnarray*}
\|V_1\|_{L^2(\mathbb{R}, L^2(S^1)^2)}^2 & = &
\sum_n \int_\mathbb{R} |\hat{V}^n_1(x,y,\sigma)|^2 \,\rmd x
\\ & \leq &
C \sum_n \left\| \rme^{-\sqrt{\sigma+1+\rmi n}|\cdot|} * \sum_k a_{n-k}(\cdot) \rme^{-\sqrt{\sigma+1+\rmi k}|\cdot -y|}\right\|^2_{L^2(\mathbb{R})}
\\ & \leq &
C \sum_n \left\|\rme^{-\sqrt{\sigma+1+\rmi n}|\cdot|} \right\|_{L^2(\mathbb{R})}^2 
\left\| \sum_k a_{n-k}(\cdot) \rme^{-\sqrt{\sigma+1+\rmi k}|\cdot -y|}\right\|_{L^1(\mathbb{R})}^2
\\ & \leq &
C \sum_n \left\|\rme^{-\sqrt{\sigma+1+\rmi n}|\cdot|} \right\|_{L^2(\mathbb{R})}^2 \left( \sum_k \|a_{n-k}(\cdot) 
\rme^{-\sqrt{\sigma+1+\rmi k}|\cdot -y|}\|_{L^1(\mathbb{R})}\right)^2
\\ & \leq &
C \sum_n \left\|\rme^{-\sqrt{\sigma+1+\rmi n}|\cdot|} \right\|_{L^2(\mathbb{R})}^2 \left( \sum_k \|a_{n-k}\|_{L^\infty(\mathbb{R})} \|\rme^{-\sqrt{\sigma+1+\rmi k}|\cdot|}\|_{L^1(\mathbb{R})}\right)^2
\\ & \leq &
C \sum_n \frac{1}{|\sigma+1+\rmi n|^{\frac12}} \left( \sum_k \|a_{n-k}\|_{L^\infty} \frac{1}{|\sigma+1+\rmi k |^{\frac12} }\right)^2.
\end{eqnarray*}
We now reorder the second sum by defining $j:=n-k$ and use the fact that
\[
(1+|n|^\alpha) \leq C(1+|n-j|^\alpha)(1+|j|^\alpha)
\]
to obtain
\begin{eqnarray*}
\|V_1\|_{L^2(\mathbb{R}, L^2(S^1)^2)}^2
& \leq &
C \sum_n \frac{1}{|\sigma+1+\rmi n|^{\frac12}} \left( \sum_j \|a_j\|_{L^\infty} 
\frac{1}{|\sigma+1+\rmi(n-j) |^{\frac12} }\right)^2
\\ & \leq &
C \sum_n \frac{1}{|\sigma+1+\rmi n|^{\frac12}} \left( \sum_j \|a_j\|_{L^\infty} \frac{(|\sigma+1|^{\frac12}+|j|^{\frac12})}{|\sigma+1|^{\frac12}(|\sigma+1|^{\frac12}+|n|^{\frac12})} \right)^2
\\ & \leq &
C \sum_n \frac{1}{1+|n|^{\frac32}},
\end{eqnarray*}
which is sufficient for convergence and guarantees that $V_1\in L^2(\mathbb{R},(L^2(S^1))^2)$. In the above derivation, we used the assumptions that $f$ and $\bar{u}$ are sufficiently smooth so that $\{a_k\}\in\ell^{1,\alpha}(L^\infty_x(\mathbb{R}))$ with $\alpha=\frac12$, where $\ell^{1,\alpha}(L^\infty_x(\mathbb{R}))=\{\{u_k\}:\;\sum_k |k|^{\alpha}\|u_k\|_{L^\infty}<\infty\}$. Note that $V_1$ is analytic in $\sigma$ for $\Re\sigma>-1$.

From this point onwards, we can obtain the Fourier--Laplace transforms $V_j$ of $\mathcal{G}_j$ inductively by solving
\[
\partial_x V_j = \mathcal{A}_0(\sigma) V_j + \mathcal{B}(x) V_{j-1}.
\]
In fact, the following lemma will allow us to repeat the bootstrapping procedure until we obtain any degree of smoothness in the inhomogeneity that we like, subject to restrictions only from the smoothness of $f(u)$.

\begin{Lemma}\label{l:bs1}
Assume that $\{c_k\}_{k\in\mathbb{Z}}\in\ell^{1,\alpha}(L^\infty_x(\mathbb{R}))$ for some $\alpha\geq\beta>0$. If $\tilde{V}$ satisfies
\[
[\tilde{V}(x,y)](t) = \sum_{n\in\mathbb{Z}} \rme^{\rmi nt} \tilde{V}_n(x,y), \qquad
\|\tilde{V}_n(\cdot,y)\|_{L^1} \leq \frac{C}{1+|n|^\beta}
\]
uniformly in $y$, and $V$ is defined by
\[
[V(x,y)](t) = \sum_{n\in\mathbb{Z}} \rme^{\rmi nt} V_n(x,y), \qquad
V_n(x,y) := \int_\mathbb{R} \rme^{-\sqrt{\sigma+1+\rmi n}|x-z|} \sum_{k\in\mathbb{Z}} c_{n-k}(z) \tilde{V}_k(z,y)\, \rmd z,
\]
then the functions $V_n$ satisfy
\[
\|V_n(\cdot,y)\|^2_{L^2} \leq \frac{\tilde{C}}{1+|n|^{\frac12+2\beta}}
\qquad \mbox{and} \qquad
\|V_n(\cdot,y)\|_{L^1}^2 \leq \frac{\tilde{C}}{1+|n|^{1+2\beta}},
\]
where the constant $\tilde{C}$ is independent of $\sigma$ and $n$.
\end{Lemma}

\begin{Proof}
The assertion follows from a calculation similar to the one used above to bound $\hat{V}^n_1$. We have
\begin{eqnarray*}
\|V_n(\cdot,y)\|^2_{L^2} & = &
\left\| \rme^{-\sqrt{\sigma+1+\rmi n}|\cdot|} * \left( \sum_{k\in\mathbb{Z}} c_{n-k}(\cdot)\tilde{V}_k(\cdot,y) \right) \right\|_{L^2}^2 \\ & \leq &
C \left\| \rme^{-\sqrt{\sigma+1+\rmi n}|\cdot|} \right\|_{L^2}^2
\left\| \sum_{k\in\mathbb{Z}} c_{n-k}(\cdot) \tilde{V}(\cdot,y) \right\|_{L^1}^2 \\ & \leq &
\frac{C}{1+|n|^\frac{1}{2}} \left( \sum_{k\in\mathbb{Z}} \|c_{n-k}\|_{L^\infty} \| \tilde{V}(\cdot,y)_k \|_{L^1} \right)^2 \\ & \leq &
\frac{C}{1+|n|^\frac{1}{2}} \left( \|c\|_{\ell^{1,\beta}(L^\infty(\mathbb{R}))} \frac{1}{1+|n|^\beta} \right)^2 \\ & \leq &
\frac{C}{1+|n|^{\frac{1}{2} + 2\beta}}.
\end{eqnarray*}
An analogous estimate holds for the $L^1$ norm of $V_n(\cdot,y)$.
\end{Proof}

The following lemma illustrates how the regularity of the functions $V_j$ increases after each iteration of the bootstrapping procedure.

\begin{Lemma}\label{l:bs2}
Pick any sufficiently small $\epsilon>0$, then, for each $j\geq0$, there is an $\eta>0$ so that the function $V_j$ is analytic for $\Re\sigma>-\frac12$ and
\begin{eqnarray*}
\lefteqn{ \rme^{\eta|\cdot-y|} V_j, \;
\rme^{\eta|\cdot-y|}\rme^{\eta|y|}(\partial_y V_j+\partial_x V_j) } \\ & \in &
L^2_x(\mathbb{R},H_t^{j+\frac12-\epsilon}\times H_t^{j-\frac12-\epsilon}) \cap
L^1_x(\mathbb{R},H_t^{j+1-\epsilon}\times       H_t^{j-\epsilon}) \cap
H^1_x(\mathbb{R},H_t^{j-\frac12-\epsilon}\times H_t^{j-\frac32-\epsilon}),
\end{eqnarray*}
where $H^k_t:=H^k(S^1)$ for all $k$.
\end{Lemma}

\begin{Proof}
The claims about $V_0$ and $V_1$ follows from their explicit formulas and estimates of the type given above: If we restrict $\sigma$ to $\Re\sigma>-\frac12$, we can also extract a factor $\rme^{-\eta|x-y|}$ from the convolution integral (\ref{E:G1}). Furthermore, applying $\partial_y$ and integrating by parts proves the claim about $(\partial_x+\partial_y)V_1$ as we have $\|\partial_z\mathcal{B}(z)\|_{L(Y_m)}\leq\rme^{-\theta|z|}$.

Next, note that the coefficients $\{ a_j \}$ and $\{ b_j \}$ that denote the Fourier components of $f_u(\bar{u})$ and $f_{uu}(\bar{u})[\bar{u}_x,\cdot]$, respectively, satisfy the hypothesis on $\{c_j\}$ in Lemma~\ref{l:bs1} for a value of $\alpha$ that is determined by the smoothness of the nonlinearity $f(u)$. To estimate $V_j$, we use the fact that
\[
V_{j}(x,y) =
\int_{-\infty}^x \Phi^{\mathrm{s}}_0(x,z,\sigma) \mathcal{B}(z) V_{j-1}(z,y) \,\rmd z
+ \int_\infty^x \Phi^{\mathrm{u}}_0(x,z,\sigma) \mathcal{B}(z) V_{j-1}(z,y) \,\rmd z
\]
and denote by $\tilde{V}_j^n$ the least well-behaved of the two components of each Fourier mode $\hat{V}_j^n$. If $V_{j-1}$ satisfies
\[
V_{j-1} \in
L^2(\mathbb{R},H^{j+\frac12-\epsilon}(S^1) \times H^{j-\frac12-\epsilon}(S^1)) 
\cap
L^1(\mathbb{R},H^{j+1-\epsilon}(S^1) \times H^{j-\epsilon}(S^1)),
\]
then Lemma~\ref{l:bs1} implies that
\begin{eqnarray*}
\|V_j\|^2_{L^1(\mathbb{R},L^2(S^1)\times L^2(S^1))} & \leq &
C \sum_{n\in\mathbb{Z}} \left\| \int_\mathbb{R} \rme^{-\sqrt{\sigma+1+in}|x-z|}
\begin{pmatrix} \frac{1}{\sqrt{\sigma+1+in}} \\ 1 \end{pmatrix}
\sum_{k\in\mathbb{Z}} a_{n-k}(z) \tilde{V}_j^n(z,y)\,\rmd z \right\|_{L^1(\mathbb{R})}^2
\\ & \leq &
\sum_{n\in\mathbb{Z}} \begin{pmatrix} \frac{1}{1+|n|} \frac{1}{1+|n|^{j+2}} \\
\frac{1}{1+|n|^{j+2}} \end{pmatrix}
\\ & \leq &
C \sum_{n\in\mathbb{Z}} \begin{pmatrix} \frac{1}{1+|n|^{j+3}} \\
\frac{1}{1+|n|^{j+2}} \end{pmatrix}
\end{eqnarray*}
and
\begin{eqnarray*}
\|V_j\|^2_{L^2(\mathbb{R},L^2(S^1)\times L^2(S^1))} & \leq &
C \sum_{n\in\mathbb{Z}} \left\| \int_\mathbb{R} \rme^{-\sqrt{\sigma+1+in}|x-z|}
\begin{pmatrix} \frac{1}{\sqrt{\sigma+1+in}} \\ 1 \end{pmatrix}
\sum_{k\in\mathbb{Z}} a_{n-k}(z) \tilde{V}_j^n(z,y)\,\rmd z \right\|_{L^2(\mathbb{R})}^2
\\ & \leq &
C \sum_{n\in\mathbb{Z}} \begin{pmatrix} \frac{1}{1+|n|^{j+\frac52}} \\
\frac{1}{1+|n|^{j+\frac32}} \end{pmatrix},
\end{eqnarray*}
which gives the $L^1_x$ and $L^2_x$ bounds on $V_j$. Note that we can again extract a factor $\rme^{-\eta|x-y|}$ from the convolution integrals. To address the $H^1_x$ bound, note that the derivative $\partial_x$ falls only on the exponential $\rme^{-\sqrt{\sigma+1+\rmi n}|x-z|}$, and the resulting extra factor of $|n|^{\frac12}$ leads to the loss of smoothness in time stated in the lemma.
\end{Proof}

We now turn to equation (\ref{E:sd_evalue_inhom}), which is equivalent to
\begin{equation}\label{n:sdg*}
\partial_x U = \mathcal{A} U + \mathcal{B}(x) V_\ell(x).
\end{equation}
If we choose $\ell=3$, then
\[
P_2\mathcal{B}V_3 \in L^2_x(\mathbb{R},H^{2+\gamma}(S^1)) \cap H^1_x(\mathbb{R},H^{1+\gamma}(S^1))
\]
for each fixed $0<\gamma\ll1$, where $P_2$ denotes the projection onto the second component. Thus, if we set $m=2+\gamma$ in the definition of the underlying space $Y_m$ so that $Y_m=H^{\frac52+\gamma}(S^1)\times H^{2+\gamma}(S^1)$, then $\mathcal{B}V_3\in L^2(\mathbb{R},Y_m)$. As discussed above, (\ref{n:sdg*}) then has a unique solution $U_*=(G_*,\partial_x G_*)$, which depends analytically on $\sigma$ and lies in $H^1(\mathbb{R},Y_m)$ so that
\[
\begin{pmatrix} G_* \\ \partial_x G_* \end{pmatrix} \in
H^1_x(\mathbb{R},H^{\frac52+\gamma}(S^1) \times H^{2+\gamma}(S^1)).
\]
In particular, $[G_*(x,\sigma)](t)$ and $[\partial_x G_*(x,\sigma)](t)$ are continuous functions. Inspecting (\ref{n:sdg*}) and using the regularity of $U_*$, we find that $\partial_x^2 G_* \in H^1_x(\mathbb{R},H^{1+\gamma}(S^1))$ so that $\partial_x^2 G_*$ is also continuous in $(x,t)$.

It remains to discuss regularity with respect to $(y,s)$. Differentiability with respect to $y$ is a consequence of the iteration scheme together with Lemma~\ref{l:bs2}. Differentiability with respect to $s$ follows similarly: since we replaced $t$ by $t+s$ at the beginning of our analysis, the initial inhomogeneity $\Delta_0(x)$ from (\ref{E:G0_eqn}) does not depend on $s$, and the dependence of $V_0$ on $s$ is only through $\mathcal{B}(x)$, that is, through $f_u(\bar{u}(x,t+s))$ and its derivatives. In particular, the recursive construction of $V_j$ shows that $\partial_s V_j$ lies in the same space as $V_j$. Once we shift back to the original time variable, $s$-derivatives of $V_j$ become equivalent to $t$-derivatives and therefore lose one degree of regularity in time.

This completes the proof of Proposition~\ref{p:rk}. We record as a corollary the representation of $G_*$ given through the variation-of-constants formula (\ref{e:sf}):

\begin{Corollary}\label{c:g*}
The resolvent kernel $G_*$ can be represented as
\begin{eqnarray}\label{E:G_dich_reln}
[G_*(x,y,\sigma;s)](t) & = &
P_1\left[ \int_{-\infty}^x \Phi^\mathrm{s}(x,z,\sigma)
\mathcal{B}(z) V_3(z,y,\sigma;s) \,\rmd z \right](t)
\\ \nonumber & &
+ P_1\left[ \int_\infty^x \Phi^\mathrm{u}(x,z,\sigma)
\mathcal{B}(z) V_3(z,y,\sigma;s) \,\rmd z \right](t),
\end{eqnarray}
where $\Phi^\mathrm{s,u}$ is the exponential dichotomy associated with the operator $\mathcal{A}$ on the space $Y_m$, and $P_1$ is the projection onto the first component.
\end{Corollary}

In \S\ref{S:res_decomp_est}, we will show that the exponential dichotomies $\Phi^\mathrm{s}(x,y,\sigma)$ and $\Phi^\mathrm{u}(x,y,\sigma)$ can be extended meromorphically by separating out the translational and essential eigenmodes. The most natural way to transfer this result to the full resolvent kernel $G(x,y,\sigma,s)$ is via the formal expression
\begin{equation}\label{E:G_natural}
G(x,y,\sigma;s) = 
\left\{\begin{array}{rcl}
 P_1 \Phi^\mathrm{s}(x,y,\sigma)\delta(\cdot) & \quad & x>y \\
-P_1 \Phi^\mathrm{u}(x,y,\sigma)\delta(\cdot) & \quad & x<y,
\end{array}\right.
\end{equation}
where $\cdot$ denotes the argument $t$. However, as discussed above, we cannot apply the dichotomy directly to $\delta(t)$. If one could argue, possibly using test functions and the uniqueness of strong solutions to equation (\ref{E:sd_evalue_inhom}), that the Green's distributions given through $\mathcal{G}=\mathcal{G}_*+\sum_j\mathcal{G}_j$ and via (\ref{E:G_natural}) must be equivalent in a distributional sense, then one could instead work directly with the dichotomies via (\ref{E:G_natural}). However, we do not know how to make such an argument work because test functions would be smooth in $t$. In order to solve the initial value problem associated with (\ref{e:lin}), we must work with initial data of the form $\delta(t)u_0(x)$, which are not smooth in time. The function $\delta(t)$ represents the spreading of the initial data amongst all Fourier modes and is a key aspect of the dynamics. This appears to be an important distinction between the time-periodic and time-independent problems. 


\section{Meromorphic extension and bounds for the resolvent kernel}\label{S:res_decomp_est}

The goal of this section is to extend the resolvent kernel $G_*(x,y,\sigma;s)$ that we constructed in the last section meromorphically across the essential spectrum near $\sigma=0$ and to derive sharp pointwise bounds for $G_*$ with respect to $(x,y)$. To state our result, we define the spatial eigenvalues $\nu_j^\pm(\sigma)$ as the $N$ solutions of the characteristic equation
\[
\det(\nu^2 - f_u(u_\pm) \nu - \sigma) = 0
\]
that are close to zero when $\sigma$ is close to zero. These eigenvalues are analytic in $\sigma$ and have the expansion
\begin{equation*}
\nu_j^\pm(\sigma) = -\frac{\sigma}{a_j^\pm} + \frac{2\sigma^2}{[a_j^\pm]^3} + \mathrm{O}(|\sigma|^3), \qquad
j=1,\ldots,N.
\end{equation*}
We shall also use the notation $\nu_\mathrm{out}^\pm$ and $\nu_\mathrm{in}^\pm$ to denote the spatial eigenvalues associated with the outgoing and incoming characteristics $a_\mathrm{out}^\pm$ and $a_\mathrm{in}^\pm$, respectively; see Definition~\ref{notation}. The following theorem is the main result of this section.

\begin{Theorem}\label{resbounds}
Assume that Hypothesis~\ref{H1} is met and that the shock profile $\bar{u}(x,t)$ is spectrally stable so that (S1)--(S4) in Definition~\ref{spectralstab} are met, then there exist positive constants $C$, $\eta$ and $\epsilon$ so that the following is true: The resolvent kernel $[G_*(x,y,\sigma;s)](t)$ has a meromorphic extension in $\sigma$ into $\{\sigma\in\mathbb{C}:\;\Re\sigma\geq-\epsilon\}$ and can be written as
\[
G_*(x,y,\sigma;s) = E_1(x,y,\sigma;s) + E_2(x,y,\sigma;s) + R(x,y,\sigma;s),
\]
where the terms $E_j$ have a pole at $\sigma=0$, while $R$ is analytic in $\sigma$. For $y\leq0$, we have
\begin{eqnarray}\label{E1lambda}
E_1(x,y,\sigma,t;s) & = & \frac{1}{\sigma} \sum_{\nu_\mathrm{in}^-} \bar{u}_x(x,t) l_{1,\mathrm{in}}^{-}(y,s)^T \rme^{-\nu_\mathrm{in}^-(\sigma)y} \\ \nonumber
E_2(x,y,\sigma,t;s) & = & \frac{1}{\sigma} \sum_{\nu_\mathrm{in}^-} \bar{u}_t(x,t) l_{2,\mathrm{in}}^{-}(y,s)^T \rme^{-\nu_\mathrm{in}^-(\sigma)y}
\end{eqnarray}
for appropriate functions $l_{j,\mathrm{in}}(y,s)$ that are $2\pi$-periodic in $s$, and a symmetric representation holds for $y\geq0$. The remainder term $R$ satisfies the following pointwise bounds, where $\alpha$ is any multi-index with $0\leq|\alpha_x|+|\alpha_y|\leq2$ and $0\leq|\alpha_t|,|\alpha_s|\leq1$:
\begin{enumerate}[(i)]
\item For $x>y>0$, we have
\begin{eqnarray*}
\sup_{s,t}|R(x,y,\sigma,t;s)| & \leq & 
C \sum_{\nu_\mathrm{out}^+,\nu^+} \rme^{\nu_\mathrm{out}^+(\sigma) x} \rme^{-\nu^+(\sigma)y},
\\ \nonumber
\sup_{s,t}|\partial_{x,y,t,s}^\alpha R(x,y,\sigma,t;s)| & \leq & 
C (|\sigma|+\rme^{-\eta|y|})^{\alpha_y}
\sum_{\nu_\mathrm{out}^+,\nu^+} \rme^{\nu_\mathrm{out}^+(\sigma) x} \rme^{-\nu^+(\sigma)y};
\end{eqnarray*}
\item For $x>0>y$, we have
\begin{eqnarray*}
\sup_{s,t}|R(x,y,\sigma,t;s)| & \leq & 
C \sum_{\nu_\mathrm{in}^-,\, \nu_\mathrm{out}^+} \rme^{\nu_\mathrm{out}^+(\sigma)x}\rme^{-\nu_\mathrm{in}^-(\sigma)y},
\\ \nonumber
\sup_{s,t}|\partial_{x,y,t,s}^\alpha R(x,y,\sigma,t;s)| & \leq & 
C (|\sigma|+\rme^{-\eta|y|})^{\alpha_y}
\sum_{\nu_\mathrm{in}^-,\, \nu_\mathrm{out}^+} \rme^{\nu_\mathrm{out}^+(\sigma)x}\rme^{-\nu_\mathrm{in}^-(\sigma)y};
\end{eqnarray*}
\item For $0>x>y$, we have
\begin{eqnarray*}
\sup_{s,t}|R(x,y,\sigma,t;s)| & \leq & 
C \sum_{\nu_\mathrm{in}^-,\, \nu^-} \rme^{\nu^-(\sigma) x} \rme^{-\nu_\mathrm{in}^-(\sigma)y},
\\ \nonumber
\sup_{s,t}|\partial_{x,y,t,s}^\alpha R(x,y,\sigma,t;s)| & \leq & 
C (|\sigma|+\rme^{-\eta|y|})^{\alpha_y}
\sum_{\nu_\mathrm{in}^-,\, \nu^-} \rme^{\nu^-(\sigma) x} \rme^{-\nu_\mathrm{in}^-(\sigma)y}.
\end{eqnarray*}
\end{enumerate}
Symmetric bounds hold for $x<y$.
\end{Theorem}

In the remainder of this section, we prove Theorem~\ref{resbounds}. Due to Corollary~\ref{c:g*}, it suffices to extend the two integral terms
\begin{equation}\label{e:ged}
\int_{-\infty}^x \Phi^\mathrm{s}(x,z,\sigma) \mathcal{B}(z) V_\ell(z,y,\sigma;s) \,\rmd z
+ \int_\infty^x \Phi^\mathrm{u}(x,z,\sigma) \mathcal{B}(z) V_\ell(z,y,\sigma;s) \,\rmd z,
\end{equation}
with $\ell=3$, in the representation (\ref{E:G_dich_reln}) of the resolvent kernel $G_*$. Since $\mathcal{B}$ does not depend on $\sigma$, and $V_\ell$ is analytic near $\sigma=0$, we need to extend the exponential dichotomies $\Phi^\mathrm{s}(x,z,\sigma)$ and $\Phi^\mathrm{u}(x,z,\sigma)$ for $x\gtrless z$ with $x,z\in\mathbb{R}$. First, we use spatial dynamics to extend the exponential dichotomies on the half lines $\mathbb{R}^+$ and $\mathbb{R}^-$ analytically across $\sigma=0$. Afterwards, we use the assumptions on the Floquet spectrum to construct a meromorphic extension of the exponential dichotomy on $\mathbb{R}$ and derive pointwise bounds for this extension. Finally, we transfer these bounds to the resolvent kernel $G_*$ by estimating the integrals (\ref{e:ged}).


\subsection{Analytic extension of the exponential dichotomies on $\mathbb{R}^\pm$}

Consider the spatial-dynamical system (\ref{E:sd_evalue}),
\begin{equation}\label{e:evp}
U_x = \mathcal{A}(x,\sigma) U
= \begin{pmatrix} 0 & 1 \\ \partial_t+\sigma+f_{uu}(\bar{u})[\bar{u}_x,\cdot] & f_u(\bar{u}) \end{pmatrix} U,
\end{equation}
on the space $Y_m=H^{m+\frac12}\times H^{m}$ for $m>2$ fixed. For $\Re\sigma>0$, this equation possesses the exponential dichotomies $\Phi^\mathrm{s}(x,y,\sigma)$ and $\Phi^\mathrm{u}(x,y,\sigma)$, which are defined and analytic in $\sigma$ for $x>y$ and $x<y$, respectively, with $x,y\in\mathbb{R}$. Our goal is to construct analytic extensions of these dichotomies separately for $x,y\in\mathbb{R}^+$ and $x,y\in\mathbb{R}^-$ from $\Re\sigma>0$ to a small ball $B_\epsilon(0)$ centered at $\sigma=0$. Throughout this section, $\epsilon$ denotes a positive, and possibly small, constant that we may adjust during the arguments to follow.

First, we consider the asymptotic equations
\begin{equation}\label{e:aed}
U_x
= \begin{pmatrix} 0 & 1 \\ \partial_t + \sigma & f_u(u_\pm) \end{pmatrix} U
=: \mathcal{A}_\pm(\sigma) U.
\end{equation}
Using the fact that the operators $\mathcal{A}_\pm(\sigma)$ leave the $2N$-dimensional subspaces $\Span\{\rme^{\rmi kt}\hat{V}:\;\hat{V}\in\mathbb{C}^{2N}\}\subset Y_m$ invariant for each $k\in\mathbb{Z}$, it was shown in \cite{SandstedeScheel06} that their spectrum is discrete and given by the spatial eigenvalues
\begin{equation}\label{E:spatial_evalues}
\frac{a_j^\pm}{2}+\frac{1}{2}\sqrt{[a_j^\pm]^2+4(\sigma+\rmi k)}, \qquad
\frac{a_j^\pm}{2}-\frac{1}{2}\sqrt{[a_j^\pm]^2+4(\sigma+\rmi k)}, \qquad
j=1,\dots,N, \quad k\in\mathbb{Z},
\end{equation}
where the $a_j^\pm$ are the nonzero, real, distinct eigenvalues of $f_u(u_\pm)$ guaranteed by Hypothesis~\ref{H1}. Furthermore, there is an $\eta>0$ so that these eigenvalues have distance $3\eta$ from the imaginary axis, uniformly in $\sigma\in B_\epsilon(0)$, except for the $N$ spatial eigenvalues
\[
\nu_j^\pm(\sigma) = -\frac{\sigma}{a_j^\pm} + \frac{2\sigma^2}{[a_j^\pm]^3} + \mathrm{O}(|\sigma|^3), \qquad
j=1,\ldots,N,
\]
which arise from (\ref{E:spatial_evalues}) when setting $k=0$ and expanding in $\sigma$ near zero. The eigenvectors of $\mathcal{A}_\pm(\sigma)$ associated with the eigenvalues $\nu_j^\pm(\sigma)$ do not depend on $t$ and are given by
\begin{equation*}
\mathcal{V}_j^\pm(\sigma) :=
\begin{pmatrix} 1 \\ \nu_j^\pm(\sigma) \end{pmatrix} r_j^\pm, \qquad
j=1,\ldots,N,
\end{equation*}
where $r_j^\pm$ are the right eigenvectors of $f_u(u_\pm)$ belonging to $a_j^\pm$. Key to our analysis is the fact that these eigenvectors are linearly independent and analytic in $\sigma$ for all $\sigma$ near zero.

\begin{figure}
\centering\includegraphics{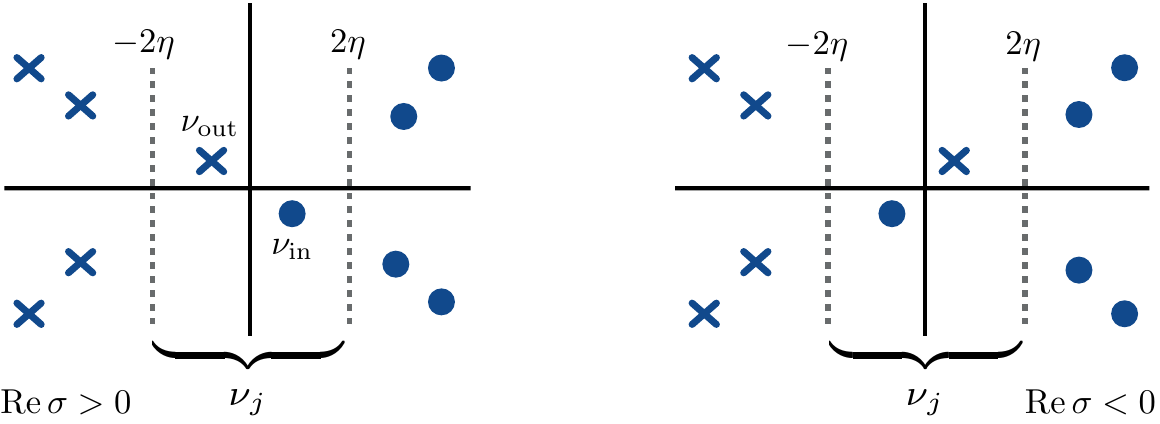}
\caption{The spatial spectrum of $\mathcal{A}_+(\sigma)$ is shown for $\Re\sigma>0$ [left] and $\Re\sigma<0$ [right]. The spatial eigenvalues that are stable for $\Re\sigma>0$ are indicated by crosses, while unstable eigenvalues are shown as bullets. As indicated, the $N$ small eigenvalues $\nu_j^+$ cross through the imaginary axis as $\Re\sigma$ changes sign.}
\label{f:ss}
\end{figure}

For $\Re\sigma>0$, the eigenvalues $\nu_j^\pm(\sigma)$ move off the imaginary axis, and the operators $\mathcal{A}_\pm(\sigma)$ are hyperbolic; see Figure~\ref{f:ss}. In fact, since $a_\mathrm{out}^-<0<a_\mathrm{in}^-$ and $a_\mathrm{in}^+<0<a_\mathrm{out}^+$, we see that
\begin{equation*}
\Re\nu_\mathrm{out}^-(\sigma)>0
\quad\mbox{and}\quad
\Re\nu_\mathrm{out}^+(\sigma)<0
\quad\mbox{for}\quad
\Re\sigma>0,
\end{equation*}
so that these spatial eigenvalues contribute respectively to the unstable eigenspace $\tilde{E}^\mathrm{u}_-(\sigma)$ of $\mathcal{A}_-(\sigma)$ and the stable eigenspace $\tilde{E}^\mathrm{s}_+(\sigma)$ of $\mathcal{A}_+(\sigma)$ when $\Re\sigma>0$. Thus, we define the subspaces
\[
\mathcal{R}_\mathrm{out}^\pm(\sigma) = \Span\{\mathcal{V}_j^\pm(\sigma):\; a_j^\pm\gtrless0\}, \qquad
\mathcal{R}_\mathrm{in}^\pm(\sigma) = \Span\{\mathcal{V}_j^\pm(\sigma):\; a_j^\pm\lessgtr0\}
\]
of outgoing and incoming modes, respectively, which are analytic in $\sigma\in B_\epsilon(0)$. Similarly, we can define the spectral subspaces $\tilde{E}^\mathrm{ss}_\pm(\sigma)$ and $\tilde{E}^\mathrm{uu}_\pm(\sigma)$ belonging to $\mathcal{A}_\pm(\sigma)$ that correspond to the stable and unstable eigenvalues with real part less than $-3\eta$ and larger than $3\eta$, respectively. Using these definitions, it follows from the above discussion and the results in \cite{SandstedeScheel06} that the decompositions
\[
\underbrace{\tilde{E}^\mathrm{uu}_-(\sigma)\oplus\mathcal{R}_\mathrm{out}^-(\sigma)}_{=\tilde{E}^\mathrm{u}_-(\sigma)\mbox{ for }\Re\sigma>0}
\oplus \tilde{E}^\mathrm{ss}_-(\sigma)\oplus\mathcal{R}_\mathrm{in}^-(\sigma) = Y_m,
\qquad
\underbrace{\tilde{E}^\mathrm{ss}_+(\sigma)\oplus\mathcal{R}_\mathrm{out}^+(\sigma)}_{=\tilde{E}^\mathrm{s}_+(\sigma)\mbox{ for }\Re\sigma>0}
\oplus \tilde{E}^\mathrm{uu}_+(\sigma)\oplus\mathcal{R}_\mathrm{in}^+(\sigma) = Y_m
\]
and the associated projections $\tilde{P}^\mathrm{u}_-(\sigma)$ and $\tilde{P}^\mathrm{s}_+(\sigma)$ exist and are analytic for $\sigma\in B_\epsilon(0)$. In particular, the unstable subspace $\tilde{E}^\mathrm{u}_-(\sigma)$ of $\mathcal{A}_-$, the stable subspace $\tilde{E}^\mathrm{s}_+(\sigma)$ of $\mathcal{A}_+$, and their spectral complements can be extended analytically from $\Re\sigma>0$ to the ball $B_\epsilon(0)$.

\begin{Lemma}\label{l:aed}
There are positive constants $C$ and $\epsilon$ so that (\ref{e:evp}) has an exponential dichotomy $\Phi^\mathrm{s,u}_+(x,y,\sigma)$ on $\mathbb{R}^+$ that is analytic in $\sigma\in B_\epsilon(0)$ and satisfies
\begin{eqnarray}
\|\Phi^\mathrm{s}_+(x,y,\sigma)\|_{L(Y_m)} & \leq & C \sum_{\nu_\mathrm{out}^+} \rme^{\nu_\mathrm{out}^+(\sigma)(x-y)}, \qquad x>y\geq0
\label{ed:bdd} \\ \nonumber
\|\Phi^\mathrm{u}_+(x,y,\sigma)\|_{L(Y_m)} & \leq & C \sum_{\nu_\mathrm{in}^+} \rme^{\nu_\mathrm{in}^+(\sigma)(x-y)}, \qquad y>x\geq0
\end{eqnarray}
and
\begin{eqnarray*}
\|\partial_x\Phi^\mathrm{s}_+(x,y,\sigma)\|_{L(Y_m,Y_{m-\frac12})} & \leq & C \left(|\sigma|+\rme^{-\eta|x-y|}\right) \sum_{\nu_\mathrm{out}^+} \rme^{\nu_\mathrm{out}^+(\sigma)(x-y)}, \qquad x>y\geq0 \\ \nonumber
\|\partial_x\Phi^\mathrm{u}_+(x,y,\sigma)\|_{L(Y_m,Y_{m-\frac12})} & \leq & C \left(|\sigma|+\rme^{-\eta|x-y|}\right) \sum_{\nu_\mathrm{in}^+} \rme^{\nu_\mathrm{in}^+(\sigma)(x-y)}, \qquad y>x\geq0
\end{eqnarray*}
for $\sigma\in B_\epsilon(0)$. Furthermore, the associated projection $P^\mathrm{s}_+(x,\sigma):=\Phi^\mathrm{s}_+(x,x,\sigma)$ on $\mathbb{R}^+$ satisfies $P^\mathrm{s}_+(x,\sigma)\to\tilde{P}^\mathrm{s}_+(\sigma)$ as $x\to\infty$. The same statement with symmetric bounds holds for dichotomies on $\mathbb{R}^-$.
\end{Lemma}

This result can be viewed as an infinite-dimensional version of the Gap Lemma, which was established in finite dimensions in \cite{GardnerZumbrun98,KapitulaSandstede98}.

\begin{Proof}
Since the distances of the strong stable and strong unstable spectrum of $\mathcal{A}_+(\sigma)$ to the imaginary axis are larger than $3\eta$, we know that the shock profile $\bar{u}(x,t)$ approaches its limit $u_+$ exponentially with rate $3\eta$ in the $C^1$-norm as $x\to\infty$.

We shall use the following result: If there are constants $\kappa^\mathrm{s}<\kappa^\mathrm{u}$ such that the spectrum of the asymptotic operator $\mathcal{A}_+(\sigma)$ is the union of two spectral sets defined by eigenvalues with real part respectively less than $\kappa^\mathrm{s}$ and larger than $\kappa^\mathrm{u}$, uniformly in $\sigma\in B_\epsilon(0)$, then (\ref{e:evp}) has exponential dichotomies on $\mathbb{R}^+$ with rates $\kappa^\mathrm{s}$ and $\kappa^\mathrm{u}$ as outlined in Definition~\ref{d:ed}, and these dichotomies can be chosen so that they are analytic in $\sigma$. Furthermore, as $x\to\infty$, the associated $x$-dependent projections converge exponentially with rate $\min\{3\eta,|\kappa^\mathrm{u}-\kappa^\mathrm{s}|\}$ to the spectral projections of $\mathcal{A}_+(\sigma)$ associated with the two spectral sets. This claim follows from \cite[Theorem~1]{PeterhofSandstedeScheel97} upon using exponential weights, and we refer to \cite{SandstedeScheel08morse} for further details.

Choosing $\kappa^\mathrm{s}=-2\eta<-\eta=\kappa^\mathrm{u}$, we find analytic dichotomies $\Phi^\mathrm{ss}_+(x,y,\sigma)$ and $\Phi^\mathrm{cu}_+(x,y,\sigma)$ corresponding to solutions that decay with rate at least $-2\eta$ as $x$ increases and solutions that grow not faster than with rate $\eta$ in backward time. Similarly, picking $\kappa^\mathrm{s}=\eta<2\eta=\kappa^\mathrm{u}$, we obtain dichotomies $\Phi^\mathrm{cs}_+(x,y,\sigma)$ and $\Phi^\mathrm{uu}_+(x,y,\sigma)$ that correspond to solutions which grow with rate at most $\eta$ as $x$ increases and solutions which decay with rate at least $-2\eta$ in backward time. The difference between these dichotomies is whether we subsume the $N$ center directions that belong to the small spatial eigenvalues $\nu_j^+(\sigma)$ into the unstable or the stable part of the spectrum. Using the convergence of the associated projections $P^\mathrm{cs}_+(x,\sigma)$ and $P^\mathrm{cu}_+(x,\sigma)$, we see that the subspace $E^\mathrm{c}_+(y,\sigma):=\Rg P^\mathrm{cs}_+(y,\sigma)\cap\Rg P^\mathrm{cu}_+(y,\sigma)$ has dimension $N$ for all $y\geq0$ and all $\sigma$. Furthermore, solutions $U(x)$ with initial data $U(y)$ in $E^\mathrm{c}_+(y,\sigma)$ exist for all $x\geq0$ with $U(x)\in E^\mathrm{c}_+(x,\sigma)$, since we can use $\Phi^\mathrm{cs}_+(x,y,\sigma)$ to evolve for $x>y$ and $\Phi^\mathrm{cu}_+(x,y,\sigma)$ to evolve for $x<y$. Thus, we successfully isolated the $N$-dimensional center directions from their infinite-dimensional stable and unstable counterparts on which we already have exponential dichotomies $\Phi^\mathrm{ss}_+(x,y,\sigma)$ and $\Phi^\mathrm{uu}_+(x,y,\sigma)$ that are analytic in $\sigma\in B_\epsilon(0)$.

Next, we decompose the $N$-dimensional center space $E^\mathrm{c}_+(y,\sigma)$ into two complementary subspaces which are composed of solutions that converge with uniform rate $2\eta$ to $\mathcal{R}_\mathrm{out}^+(\sigma)$ or to $\mathcal{R}_\mathrm{in}^+(\sigma)$, respectively, as $x\to\infty$. We proceed as in \cite[\S4.3]{Sandstede02}. First, we write (\ref{e:evp}) as
\begin{equation}\label{e:evpB}
U_x = [\mathcal{A}_+(\sigma) + \mathcal{B}_+(x)] U, \qquad
\mathcal{B}_+(x) = \begin{pmatrix} 0 & 0 \\ f_{uu}(\bar{u})[\bar{u}_x,\cdot] & f_u(\bar{u}) - f_u(u_+) \end{pmatrix},
\end{equation}
where $\|\mathcal{B}_+(x)\|_{L(Y_m)}\leq C\rme^{-3\eta|x|}$ for $x\geq0$. Pick an index $j$ so that $a_j^+>0$ is an outgoing characteristic. We seek a solution $U(x)$ of (\ref{e:evpB}) of the form
\begin{equation}\label{e:gapsoln}
U(x) = \rme^{\nu_j^+(\sigma)(x-L)} \mathcal{V}_j^+(\sigma) + V(x),
\end{equation}
where $L>0$ is some large constant and we require that $|V(x)|_{Y_m}\leq C\rme^{-2\eta x}$ as $x\to\infty$. To construct $V(x)$, we consider the integral equation
\begin{eqnarray}\label{e:gapie}
V(x) & = &
\int_\infty^x \rme^{\mathcal{A}_+(\sigma)\tilde{P}^\mathrm{cu}_+(\sigma)(x-z)} \tilde{P}^\mathrm{cu}_+(\sigma) \mathcal{B}_+(z) \left[ V(z) + \rme^{\nu_j^+(\sigma)(z-L)} \mathcal{V}_j^+(\sigma) \right]\,\rmd z
\\ \nonumber & &
+ \int_L^x \rme^{\mathcal{A}_+(\sigma)\tilde{P}^\mathrm{s}_+(\sigma)(x-z)} \tilde{P}^\mathrm{s}_+(\sigma) \mathcal{B}_+(z) \left[ V(z) + \rme^{\nu_j^+(\sigma)(z-L)} \mathcal{V}_j^+(\sigma) \right]\,\rmd z
\end{eqnarray}
for $x\geq L$. Upon fixing a sufficiently large $L$, it was shown in \cite[\S4.3 and (4.12)]{Sandstede02} that (\ref{e:gapie}) has a unique solution $V_j(x,\sigma)$ for $x\geq L$ that grows with rate at most $\eta$ as $x \to \infty$. Once this is established, one can show that this solution depends analytically on $\sigma$ and, in fact, converges exponentially with rate $2\eta$ to zero as $x\to\infty$, because the first integral term becomes zero in this limit. In order to construct this solution for all $x \geq 0$, we now need to flow it backward from $x = L$ to $x = 0$. However, we cannot necessarily do so in this infinite dimensional setting because we do not know whether the initial data $V_j(L,\sigma)$ lies in $E^\mathrm{c}_+(L,\sigma)$. However, it is easy to see that
\[
E^\mathrm{c}_\mathrm{out}(L,\sigma) := \left[ \Rg P^\mathrm{ss}_+(L,\sigma) \oplus \Span\{ V_j(L,\sigma):\; a_j^+>0 \} \right] \cap E^\mathrm{c}_+(L,\sigma)
\]
has the same dimension as $\mathcal{R}_\mathrm{out}^+(\sigma)$, possibly after making $L$ larger.

Proceeding in the same fashion for initial data in $\mathcal{R}_\mathrm{in}^+(\sigma)$, we can construct an analytic complement $E^\mathrm{c}_\mathrm{in}(L,\sigma)$ of $E^\mathrm{c}_\mathrm{out}(L,\sigma)$ in $E^\mathrm{c}_+(L,\sigma)$. Since we can evolve initial data in $E^\mathrm{c}_+(L,\sigma)$ for all $x\geq0$, we can define an analytic and invariant decomposition in $E^\mathrm{c}_+(x,\sigma)$ for each $x\geq0$. Adding $\Phi^\mathrm{ss}_+(x,y,\sigma)$ and $\Phi^\mathrm{uu}_+(x,y,\sigma)$ to the center evolutions that we just constructed defines an analytic extension of the exponential dichotomy on $\Phi^\mathrm{s}_+(x,y,\sigma)$ and $\Phi^\mathrm{u}_+(x,y,\sigma)$ into $B_\epsilon(0)$. The bounds stated in Lemma~\ref{l:aed} are a consequence of the ansatz (\ref{e:gapsoln}) and the exponential bounds for $V_j(x,\sigma)$.
\end{Proof}


\subsection{Meromorphic extension of the exponential dichotomy on $\mathbb{R}$}

We define $E^\mathrm{s}_+(\sigma)$ and $E^\mathrm{u}_-(\sigma)$ to be the ranges of the projections $P^\mathrm{s}_+(0,\sigma):=\Phi^\mathrm{s}_+(0,0,\sigma)$ and $P^\mathrm{u}_-(0,\sigma):=\Phi^\mathrm{u}_-(0,0,\sigma)$, respectively. It follows from \cite[Theorem~2]{PeterhofSandstedeScheel97} that the exponential dichotomies $\Phi^\mathrm{s}_+(x,y,\sigma)$ and $\Phi^\mathrm{u}_-(x,y,\sigma)$, which we defined in Lemma~\ref{e:aed} separately on $\mathbb{R}^+$ and $\mathbb{R}^-$, fit together at $x=y=0$ to produce an exponential dichotomy on $\mathbb{R}$ if and only if $E^\mathrm{s}_+(\sigma)\oplus E^\mathrm{u}_-(\sigma)=Y_m$. Thus, if this equation were true for all $\sigma$ near zero, then (\ref{e:evp}) would admit an analytic exponential dichotomy on $\mathbb{R}$ for all such $\sigma$, which could be constructed explicitly through the analytic projection onto $E^\mathrm{s}_+(\sigma)$ with null space $E^\mathrm{u}_-(\sigma)$; see \cite[(3.20)]{PeterhofSandstedeScheel97} and (\ref{e:med}) below. As we shall see in Lemma~\ref{l:fredholm} below, the direct sum decomposition of $Y_m$ through stable and unstable subspaces fails at $\sigma=0$, due to the presence of the embedded spatial and temporal translation eigenmodes. Therefore, we instead show that the projection onto $E^\mathrm{s}_+(\sigma)$ with null space $E^\mathrm{u}_-(\sigma)$ has a meromorphic extension in $\sigma$, with a pole at $\sigma=0$, which we can use to construct a meromorphic exponential dichotomy on $\mathbb{R}$.

Consider (\ref{e:evp}),
\begin{equation}\label{n:evp}
U_x = \mathcal{A}(x,\sigma) U
= \begin{pmatrix} 0 & 1 \\ \partial_t+\sigma+f_{uu}(\bar{u})[\bar{u}_x,\cdot] & f_u(\bar{u}) \end{pmatrix} U,
\end{equation}
and its formally transposed equation
\begin{equation}\label{n:adj}
W_x = -\mathcal{A}(x,\sigma)^T W
= - \begin{pmatrix} 0 & -\partial_t+\sigma+f_{uu}^T(\bar{u})[\bar{u}_x,\cdot] \\ 1 & f_u^T(\bar{u}) \end{pmatrix} W,
\end{equation}
taken with respect to the real inner product in $X=L^2(S^1)\times L^2(S^1)$. A calculation \cite{SandstedeScheel01,SandstedeScheel06} shows that
\begin{equation}\label{e:sp}
\frac{\rmd}{\rmd x} \langle W(x),U(x) \rangle_{X} = 0 \qquad \forall x
\end{equation}
for solutions $U(x)$ of (\ref{n:evp}) and $W(x)$ of (\ref{n:adj}).

We now set $\sigma=0$ and consider the resulting equations
\begin{equation}\label{n:evp0}
U_x = \mathcal{A}(x,0) U
= \begin{pmatrix} 0 & 1 \\ \partial_t+f_{uu}(\bar{u})[\bar{u}_x,\cdot] & f_u(\bar{u}) \end{pmatrix} U
\end{equation}
and
\begin{equation}\label{n:adj0}
W_x = -\mathcal{A}(x,0)^T W
= - \begin{pmatrix} 0 & -\partial_t+f_{uu}^T(\bar{u})[\bar{u}_x,\cdot] \\ 1 & f_u^T(\bar{u}) \end{pmatrix} W.
\end{equation}
Equation (\ref{n:evp0}) admits the two linearly independent solutions
\begin{equation*}
\overline{U}_1(x) = \partial_x \begin{pmatrix} \bar{u} \\ \bar{u}_x \end{pmatrix}, \qquad
\overline{U}_2(x) = \partial_t \begin{pmatrix} \bar{u} \\ \bar{u}_x \end{pmatrix},
\end{equation*}
which are defined for $x\in\mathbb{R}$ and decay exponentially to zero as $x\to\pm\infty$. Next, consider the adjoint equation (\ref{n:adj0}), which admits the solutions
\begin{equation}\label{sol:adj0}
W(x) = \begin{pmatrix} - f_u^T(\bar{u}) w \\ w \end{pmatrix},
\end{equation}
where $w\in\mathbb{R}^N$ is arbitrary. These solutions can be used to define $N$ bounded and linearly independent solutions by substituting the $N$ basis vectors $e_j$ in $\mathbb{R}^N$ for $w$. As shown in \cite{SandstedeScheel06}, these solutions originate from the smooth functional
\[
\mathcal{E}: \quad
Y_m \longrightarrow \mathbb{R}^n, \quad
\begin{pmatrix} u \\ v \end{pmatrix} \longmapsto
\int_0^{2\pi} [v- f(u)]\,\rmd t,
\]
which is conserved under the evolution of the system
\[
\begin{pmatrix} u_x \\ v_x \end{pmatrix} =
\begin{pmatrix} v \\ \partial_t u + f_u(u)v \end{pmatrix}
\]
on $Y_m$ that time-periodic shock profiles with period $2\pi$ satisfy; see (\ref{main}).

As outlined in \S\ref{introduction}, the condition (\ref{majdacondition}) in Hypothesis~(S3) implies that there is a unique nonzero vector $\psi_1\in\mathbb{R}^N$, up to scalar multiples, that is perpendicular to the outgoing eigenvectors of $f_u(u_\pm)$ so that $\psi_1\perp[R^+_\mathrm{out}\oplus R_\mathrm{out}^-]$. We define
\[
\Psi_1(x) = \begin{pmatrix} -f_u^T(\bar{u}) \psi_1 \\ \psi_1 \end{pmatrix}
\]
to be the associated solution of (\ref{n:adj0}). Hypothesis~(S4) implies that a second solution of the adjoint equation is given by
\[
\Psi_2(x) = \begin{pmatrix} -\partial_t\psi_2 - f_u^T(\bar{u})\psi_2 \\ \psi_2 \end{pmatrix},
\]
where $\psi_2(x,t)$ appears in (S4). Assumption (\ref{E:transverse_connection}) implies that $\Psi_1$ and $\Psi_2$ are linearly independent. Recall the definition $E_\pm^j(\sigma):=\Rg P^\pm_j(0,\sigma)$ for $j=\mathrm{s,u}$. 

\begin{Lemma}\label{l:fredholm}
The linear mapping
\[
\iota(\sigma):\quad E^\mathrm{s}_+(\sigma) \times E^\mathrm{u}_-(\sigma) \longrightarrow Y_m, \quad (V^\mathrm{s},V^\mathrm{u}) \longmapsto V^\mathrm{s}-V^\mathrm{u}
\]
is Fredholm with index zero for all $\sigma\in B_\epsilon(0)$. Furthermore, we have
\[
E^\mathrm{s}_+(0) \cap E^\mathrm{u}_-(0) = \Span\{\overline{U}_1(0),\overline{U}_2(0)\}, \qquad
[E^\mathrm{s}_+(0) + E^\mathrm{u}_-(0)]^\perp = \Span\{\Psi_1(0),\Psi_2(0)\},
\]
and $E^\mathrm{s}_+(\sigma)\oplus E^\mathrm{u}_-(\sigma)=Y_m$ for all $\sigma\in B_\epsilon(0)\setminus\{0\}$.
\end{Lemma}

\begin{Proof}
Spectral stability of the shock profile together with \cite[\S4]{SandstedeScheel00} implies that $\iota(\sigma)$ is invertible, and therefore Fredholm with index zero, for $\Re\sigma>0$. Furthermore, \cite[Corollary~1]{PeterhofSandstedeScheel97} implies that the nullspace of $\iota(0)$ is finite-dimensional, while \cite[Comment on p.~273]{PeterhofSandstedeScheel97} and \cite[Lemma~6.1 and \S6.2]{SandstedeScheel01} show that the range of $\iota(0)$ is closed and has finite codimension. Thus, $\iota(\sigma)$ is Fredholm with index zero for $\sigma=0$ and hence for all $\sigma\in B_\epsilon(0)$, possibly after making $\epsilon$ smaller, as the set of Fredholm operators of a given index is open.

Next, it is clear that $\overline{U}_1(0)$ and $\overline{U}_2(0)$ lie in $E^\mathrm{s}_+(0) \cap E^\mathrm{u}_-(0)$, and Hypothesis~(S2) implies that this space does not contain any other initial data that lead to nontrivial localized solutions of (\ref{n:evp0}). Hence, the proof of Lemma~\ref{l:aed} implies that any other element in this intersection corresponds to a solution $V(x)$ of (\ref{n:evp0}) for which $V(x)$ converges to a nonzero element of $\mathcal{R}_\mathrm{out}^+$ as $x\to\infty$ or of $\mathcal{R}_\mathrm{out}^-$ as $x\to-\infty$. Suppose the former case occurs with $V(x)\to V_\infty^+=(r_+,0)^T\in\mathcal{R}_\mathrm{out}^+\setminus\{0\}$ as $x\to\infty$. Since $\mathcal{R}_\mathrm{out}^-\cap \mathcal{R}_\mathrm{out}^+=\{0\}$ by (S3) and because $\mathcal{R}_\mathrm{out}^\pm$ is invariant under $f_u(u_\pm)$, we can pick a $\psi_+\in L_\mathrm{in}^-$ so that $\langle\psi_+,f_u(u_+)r_+\rangle=1$ and $\langle\psi_+,f_u(u_-)r_-\rangle=0$ for all $r_-\in R_\mathrm{out}^-$. Define the associated solution $\Psi(x)$ of (\ref{n:adj0}) through (\ref{sol:adj0}) and observe that $\langle\Psi(x),V(x)\rangle_X=-1$ for all $x$ by construction and (\ref{e:sp}). However, $V(x)$ converges to $R_\mathrm{out}^-$ (or to zero) as $x\to-\infty$, and we reach a contradiction to our choice of $\psi_+$. This proves our claim about $E^\mathrm{s}_+(0) \cap E^\mathrm{u}_-(0)$. A similar argument shows that $\Psi_1(0)$ and $\Psi_2(0)$ are perpendicular to the range of $\iota(0)$ and therefore span the complement of the range as claimed.

We can regard $\iota(\sigma)$ as being analytic in $\sigma$ by viewing the equivalent operator
\[
\iota(\sigma):\quad E^\mathrm{s}_+(0) \times E^\mathrm{u}_-(0) \longrightarrow Y_m, \quad (V^\mathrm{s},V^\mathrm{u}) \longmapsto P^\mathrm{s}_+(0,\sigma) V^\mathrm{s} - P^\mathrm{u}_-(0,\sigma) V^\mathrm{u}.
\]
Thus, $\iota(\sigma)$ has a nontrivial nullspace either for all $\sigma$ or else only for a discrete set of $\sigma$ in its region of analyticity. For $\sigma>0$, a nontrivial nullspace of $\iota(\sigma)$ corresponds to a Floquet eigenvalue, and Hypothesis~(S1) precludes their existence. Hence, we conclude that $\iota(\sigma)$ is invertible for all $\sigma\in B_\epsilon(0)\setminus\{0\}$, possibly after making $\epsilon$ smaller.
\end{Proof}

Lemma~\ref{l:fredholm} implies that there are closed subspaces $E^\mathrm{s}_0$, $E^\mathrm{u}_0$, $E^\mathrm{pt}_0$, and $E^\psi_0$ of $Y_m$ with
\begin{equation}\label{e:ds0}
\underbrace{E^\mathrm{s}_0 \oplus E^\mathrm{pt}_0}_{=E^\mathrm{s}_+(0)} \oplus
\underbrace{E^\mathrm{u}_0 \oplus E^\psi_0}_{=E^\mathrm{u}_+(0)} = Y_m,
\end{equation}
where
\[
E^\mathrm{pt}_0 = E^\mathrm{s}_+(0) \cap E^\mathrm{u}_-(0) = \Span\{\overline{U}_1(0),\overline{U}_2(0)\}, \qquad
E^\psi_0 = [E^\mathrm{s}_+(0) + E^\mathrm{u}_-(0)]^\perp = \Span\{\Psi_1(0),\Psi_2(0)\}.
\]
Note that $E^\mathrm{u}_0$ is not uniquely determined, and we shall use this freedom below in Lemma~\ref{l:med} to make a specific choice that simplifies the estimates. We define $\mathcal{P}$ to be the projection onto $E^\mathrm{u}_0 \oplus E^\psi_0$ with null space $E^\mathrm{s}_0 \oplus E^\mathrm{pt}_0$.

\begin{Lemma}\label{l:med}
For each $\sigma \in B_\epsilon(0) \setminus \{0\}$, there is a unique mapping $h_+(\sigma): E^\mathrm{u}_+(\sigma) \to E^\mathrm{s}_+(\sigma)$ so that $E^\mathrm{u}_-(\sigma)=\graph h_+(\sigma)$. Furthermore, we can choose $E^\mathrm{u}_0$ subject to (\ref{e:ds0}) so that $h_+(\sigma)$ can be written as $(1 - \mathcal{P}) \tilde{h}_+(\sigma)\mathcal{P}$, where
\[
\tilde{h}_+(\sigma): E^\mathrm{u}_0 \oplus E^\psi_0 \to E^\mathrm{s}_0 \oplus E^\mathrm{pt}_0, \qquad
\tilde{h}_+(\sigma) = \tilde{h}^+_\mathrm{a}(\sigma) + \tilde{h}^+_\mathrm{p}(\sigma),
\]
and $\tilde{h}^+_\mathrm{a}(\sigma)$ is analytic for $\sigma \in B_\epsilon(0)$, while $ \tilde{h}^+_\mathrm{p}(\sigma)$ is given by
\[
\tilde{h}^+_\mathrm{p}(\sigma)(V^\mathrm{u},V^\psi) = \frac{1}{\sigma} (0,\tilde{M}_0 V^\psi),
\]
where $\tilde{M}_0: E^\psi_0\to E^\mathrm{pt}_0$ is invertible and has the matrix representation
\[
\tilde{M}_0 = \begin{pmatrix}\displaystyle
\langle \psi_1,[\bar{u}] \rangle_{\mathbb{R}^N} & 0 \\ \displaystyle
\int_\mathbb{R} \langle \psi_2, \bar{u}_x \rangle_{L^2(S^1)}\,\rmd x & \displaystyle
\int_\mathbb{R} \langle \psi_2, \bar{u}_t \rangle_{L^2(S^1)}\,\rmd x
\end{pmatrix}^{-1}
\]
with respect to the bases $\{\overline{U}_1(0),\overline{U}_2(0)\}$ and $\{\Psi_1(0),\Psi_2(0)\}$. Similarly, for each $\sigma \in B_\epsilon(0) \setminus \{0\}$, there is a unique mapping $h_-(\sigma): E^\mathrm{s}_-(\sigma) \to E^\mathrm{u}_-(\sigma)$ so that $E^\mathrm{s}_+(\sigma)=\graph h_-(\sigma)$. This mapping has a meromorphic representation analogous to the one given above for $h_+(\sigma)$, but now involving the matrix $-\tilde{M}_0$.
\end{Lemma}

In particular, $h_+(\sigma)$ is meromorphic on $B_\epsilon(0)$ with a simple pole at $\sigma=0$.

\begin{Proof}
Our proof mimics Lyapunov--Schmidt reduction. We use the coordinates
\begin{equation*}
(V^\mathrm{s}, V^\mathrm{pt},V^\mathrm{u},V^\psi) \in E^\mathrm{s}_0 \oplus E^\mathrm{pt}_0 \oplus E^\mathrm{u}_0 \oplus E^\psi_0
\end{equation*}
and indicate the range of operators by the appropriate superscript: the mapping $g^{\mathrm{u}\psi}(\sigma)$, for instance, maps into $E^\mathrm{u}_0 \oplus E^\psi_0$.

Lemma~\ref{l:fredholm} implies that we can write $E_-^\mathrm{u}(\sigma)$ uniquely as a graph over $E^\mathrm{u}_0 \oplus E^\mathrm{pt}_0$ with values in $E^\mathrm{s}_0 \oplus E^\psi_0$. Thus, there are unique analytic mappings $h_j^\mathrm{s}(\sigma)$ and $h_j^\psi(\sigma)$ with
\[
E_-^\mathrm{u}(\sigma):\quad
V = V^\mathrm{u} + V^\mathrm{pt} + h_1^\mathrm{s}(\sigma) V^\mathrm{u} + h_2^\mathrm{s}(\sigma) V^\mathrm{pt} + h_1^\psi(\sigma) V^\mathrm{u} + h_2^\psi(\sigma) V^\mathrm{pt}.
\]
Let $\sigma=0$. Setting $V^\mathrm{u}=0$, we find that
\[
V^\mathrm{pt} + h_2^\mathrm{s}(0) V^\mathrm{pt} + h_2^\psi(0) V^\mathrm{pt} \in E^\mathrm{u}_-(0) \qquad \forall\; V^\mathrm{pt}\in E_0^\mathrm{pt}.
\]
Since $E_0^\mathrm{pt}\subset E_-^\mathrm{u}(0)$, we conclude that $h_2^\psi(0)=h_2^\mathrm{s}(0)=0$. Next, set $V^\mathrm{pt}=0$, then
\[
V^\mathrm{u} + h_1^\mathrm{s}(0) V^\mathrm{u} + h_1^\psi(0) V^\mathrm{u} \in E^\mathrm{u}_-(0) \qquad \forall\; V^\mathrm{u}\in E_0^\mathrm{u}.
\]
The only requirement for $E_0^\mathrm{u}$ is that $E^\mathrm{u}_0\oplus E^\psi_0=E_+^\mathrm{u}(0)$. Upon replacing $E^\mathrm{u}_0$ by $\graph h_1^\psi(0)\subset E_+^\mathrm{u}(0)$, we can therefore assume that $h_1^\psi(0)=0$. Thus, the above discussion shows that we can write $E_-^\mathrm{u}(\sigma)$ as
\[
E_-^\mathrm{u}(\sigma): \quad V = V^\mathrm{u} + V^\mathrm{pt} + \sigma h^\psi(\sigma)(V^\mathrm{u} + V^\mathrm{pt}) + h^\mathrm{s}(\sigma) V^\mathrm{u} + \sigma h^\mathrm{s}(\sigma) V^\mathrm{pt},
\]
where all mappings in the above expression are analytic in $\sigma$.

Since the subspaces $E^\mathrm{s}_+(\sigma)$ and $E^\mathrm{u}_+(\sigma)$ are analytic, we also have
\begin{eqnarray}\label{e:p}
E_+^\mathrm{s}(\sigma): & \quad & V = V^\mathrm{s} + V^\mathrm{pt} + \sigma g^{\mathrm{u}\psi}(\sigma) (V^\mathrm{s} + V^\mathrm{pt}) \\ \nonumber
E_+^\mathrm{u}(\sigma): & \quad & V = V^\mathrm{u} + V^\psi + \sigma g^{\mathrm{s,pt}}(\sigma) (V^\mathrm{u} + V^\psi)
\end{eqnarray}
for unique mappings $g^{\mathrm{u}\psi}$ and $g^{\mathrm{s,pt}}$ that are analytic in $\sigma$.

We need to write $E_-^\mathrm{u}(\sigma)$ as a graph over $E_+^\mathrm{u}(\sigma)$ with values in $E_+^\mathrm{s}(\sigma)$. Thus, consider
\begin{eqnarray}\label{E:med_1}
\lefteqn{\tilde{V}^\mathrm{u} + \tilde{V}^\mathrm{pt} + \sigma h^\psi(\sigma)(\tilde{V}^\mathrm{u}+\tilde{V}^\mathrm{pt}) + h^\mathrm{s}(\sigma) \tilde{V}^\mathrm{u} + \sigma h^\mathrm{s}(\sigma) \tilde{V}^\mathrm{pt}} \\ & = & \nonumber
\left[V^\mathrm{u} + V^\psi + \sigma g^{\mathrm{s,pt}}(\sigma) (V^\mathrm{u}+V^\psi)\right] + \left[V^\mathrm{s} + V^\mathrm{pt} + \sigma g^{\mathrm{u}\psi}(\sigma) (V^\mathrm{s}+V^\mathrm{pt})\right],
\end{eqnarray}
where we need to express $(V^\mathrm{s},V^\mathrm{pt})$ in terms of  $(V^\mathrm{u},V^\psi)$
so that (\ref{E:med_1}) is true. Upon writing (\ref{E:med_1}) in components, we see that
\[
\tilde{V}^\mathrm{u} = V^\mathrm{u} + \sigma g^\mathrm{u}(\sigma)(V^\mathrm{s} + V^\mathrm{pt}), \qquad
\tilde{V}^\mathrm{pt} = V^\mathrm{pt} + \sigma g^\mathrm{pt}(\sigma)(V^\mathrm{u} + V^\psi).
\]
Substituting these expressions into the stable component of (\ref{E:med_1}), we obtain
\[
h^\mathrm{s}(\sigma)[V^\mathrm{u} + \sigma g^\mathrm{u}(\sigma)(V^\mathrm{s} + V^\mathrm{pt})] + \sigma h^\mathrm{s}(\sigma)[V^\mathrm{u} + V^\psi + \sigma g^\mathrm{s,pt}(\sigma)(V^\mathrm{u} + V^\psi)] = V^\mathrm{s} + \sigma g^\mathrm{s}(\sigma)(V^\mathrm{u} + V^\psi),
\]
which we can solve for $V^\mathrm{s}$ so that
\begin{eqnarray}
V^\mathrm{s} & = &
\left[\mathrm{id}_{E^\mathrm{s}_0}-\sigma h^\mathrm{s}(\sigma)g^\mathrm{u}(\sigma)\right]^{-1} 
\Big(
h^\mathrm{s}(\sigma)[V^\mathrm{u} + \sigma g^\mathrm{u}(\sigma)V^\mathrm{pt}]
\nonumber \\ & & \nonumber
+ \sigma h^\mathrm{s}(\sigma)[V^\mathrm{u} + V^\psi + \sigma g^\mathrm{s,pt}(\sigma)(V^\mathrm{u} + V^\psi)]
- \sigma g^\mathrm{s}(\sigma)(V^\mathrm{u} + V^\psi)
\Big) \\ \label{E:med_2}  & =: &
h_1^\mathrm{s}(\sigma) V^\mathrm{u} + \sigma h_2^\mathrm{s}(\sigma)(V^\mathrm{pt} + V^\psi),
\end{eqnarray}
where $h_j^\mathrm{s}(\sigma)$ is analytic in $\sigma$. Finally, the $E^\psi_0$-component of equation (\ref{E:med_1}) is given by
\begin{eqnarray*}
\lefteqn{ \hspace*{-2.5cm} \sigma h^\psi(\sigma)\Big( V^\mathrm{u} + \sigma g^\mathrm{u}(\sigma)(h_1^\mathrm{s}(\sigma) V^\mathrm{u} + \sigma h_2^\mathrm{s}(\sigma)(V^\mathrm{pt} + V^\psi) + V^\mathrm{pt}) + V^\mathrm{pt} + \sigma g^\mathrm{pt}(\sigma)(V^\mathrm{u} + V^\psi) \Big) } \\ & = &
V^\psi + \sigma g^\psi(\sigma)(h_1^\mathrm{s}(\sigma) V^\mathrm{u} + \sigma h_2^\mathrm{s}(\sigma)(V^\mathrm{pt} + V^\psi) + V^\mathrm{pt}),
\end{eqnarray*}
which is of the form
\begin{equation*}
\left[\mathrm{id}_{E^\psi_0}+\sigma h_1^\psi(\sigma)\right] V^\psi = \sigma \underbrace{\left[h^\psi(\sigma) - g^\psi(\sigma) + \sigma h_2^\psi(\sigma)\right]}_{=: M(\sigma)} V^\mathrm{pt} + \sigma h_3^\psi(\sigma) V^\mathrm{u},
\end{equation*}
where all mappings are analytic in $\sigma$. Assume for the moment that
\[
M(0)=h^\psi(0)-g^\psi(0):\quad E^\mathrm{pt}_0\longrightarrow E^\psi_0
\]
is invertible. We then have
\begin{eqnarray*}
V^\mathrm{pt} & = &
\frac{1}{\sigma} M(\sigma)^{-1}[\mathrm{id}_{E^\psi_0}+\sigma h_2^\psi(\sigma)] V^\psi - M(\sigma)^{-1}h_3^\psi(\sigma) V^\mathrm{u} \\ & =: &
\frac{1}{\sigma} \tilde{M}(\sigma) V^\psi - h_1^\mathrm{pt}(\sigma) V^\mathrm{u},
\end{eqnarray*}
where $\tilde{M}(0)=[h^\psi(0)-g^\psi(0)]^{-1}$, and $\tilde{M}(\sigma)$ and $h_1^\mathrm{pt}(\sigma)$ are analytic in $\sigma$.
Using also (\ref{E:med_2}), we see that
\[
(V^\mathrm{s},V^\mathrm{pt}) = \tilde{h}_\mathrm{a}^+(\sigma)(V^\mathrm{u},V^\psi) + \frac{1}{\sigma}(0,\tilde{M}(0)V^\psi)
\]
as claimed.

It remains to prove that $M(0)=h^\psi(0)-g^\psi(0)$ is invertible and to derive an expression for its matrix representation. Hence, we need to find expressions for the $E^\psi_0$-components of the graphs of $E_+^\mathrm{s}(\sigma)$ and $E_-^\mathrm{u}(\sigma)$ over $E_0^\mathrm{pt}$. We start with $E_+^\mathrm{s}(\sigma)=\Rg\Phi^\mathrm{s}_+(0,0,\sigma)$. Using the definition (\ref{e:p}) of $g^\psi(\sigma)$ and recalling that $X=L^2(S^1)\times L^2(S^1)$, we see that
\begin{equation}\label{e:der}
\langle \Psi_i(0),g^\psi(0) \overline{U}_j(0) \rangle_X = \langle \Psi_i(0), \mathrm{D}_\sigma\Phi^\mathrm{s}_+(0,0,\sigma)|_{\sigma=0} \overline{U}_j(0) \rangle_X,
\end{equation}
since the $E^\psi_0$-components of the projection $P^\mathrm{s}_+(0,\sigma)$ and the parametrization (\ref{e:p}) differ only at order $\mathrm{O}(\sigma^2)$. To derive an expression for $\mathrm{D}_\sigma \Phi_+^\mathrm{s}(0,0,0)\overline{U}_j(0)$, we recall that $\Phi_+^\mathrm{s}(x,0,\sigma)\overline{U}_j(0)$ satisfies equation (\ref{e:evp}), which can be written as
\[
U_x = \left[\mathcal{A}(x,0) + \sigma \begin{pmatrix} 0 & 0 \\ 1 & 0 \end{pmatrix} \right] U.
\]
Thus, $V(x)=\mathrm{D}_\sigma \Phi_+^\mathrm{s}(x,0,0)\overline{U}_j(0)$ is a bounded solution of
\[
V_x = \mathcal{A}(x,0)V + \begin{pmatrix} 0 & 0 \\ 1 & 0 \end{pmatrix} \overline{U}_j(x),
\]
where we used the fact that $\Phi_+^\mathrm{s}(x,0,0)\overline{U}_j(0) = \overline{U}_j(x)$. Using the construction of $\Phi^\mathrm{s}_+(x,y,\sigma)$ in Lemma~\ref{l:aed}, we know that there are analytic coefficients $\alpha_{jk}(\sigma)$ so that
\begin{equation}\label{E:med_4}
\Phi_+^\mathrm{s}(x,0,\sigma)\overline{U}_j(0) = \Phi_+^\mathrm{ss}(x,0,\sigma)\overline{U}_j(0) + \sum_{k=1}^{p-1} \sigma \alpha_{jk}(\sigma) V_k(x,\sigma),
\end{equation}
where the factor $\sigma$ arises since $\overline{U}_j(x)=\Phi_+^\mathrm{ss}(x,0,0)\overline{U}_j(0)$ decays exponentially with rate at least $3\eta$. Thus, $V(x)$ is of the form
\[
V(x) = \Phi_+^\mathrm{s}(x,0,0) V_0^+ + \underbrace{\int_0^x \Phi_+^\mathrm{ss}(x,z,0) \begin{pmatrix} 0 & 0 \\ 1 & 0 \end{pmatrix} \overline{U}_j(z)\, \rmd z + \int_\infty^x \Phi_+^\mathrm{cu}(x,z,0) \begin{pmatrix} 0 & 0 \\ 1 & 0 \end{pmatrix} \overline{U}_j(z)\, \rmd z}_{\text{particular solution converges to }0\text{ as }x\to\infty} 
\]
for an appropriate $V_0^+$: indeed, (\ref{E:med_4}) shows that $V(x)\to\mathcal{R}^+_\mathrm{out}$ as $x\to\infty$, as does the particular solution in the above expression; thus, the only other contribution to $V(x)$ can come from $\Phi_+^\mathrm{s}(x,0,0)$. Hence, 
\[
V(0) = \Phi_+^\mathrm{s}(0,0,0) V_0^+ + \int_\infty^0 \Phi_+^\mathrm{cu}(0,z,0) \begin{pmatrix} 0 & 0 \\ 1 & 0 \end{pmatrix} \overline{U}_j(z)\,\rmd z
\]
and, setting $u_1:=\bar{u}_x$ and $u_2:=\bar{u}_t$, we obtain
\begin{eqnarray*}
\langle \Psi_i(0),g^\psi(0) \overline{U}_j(0) \rangle_X & \stackrel{(\ref{e:der})}{=} &
\langle \Psi_i(0),V(0) \rangle_X \\ & = &
\int_\infty^0 \left\langle \Psi_i(0), \Phi_+^\mathrm{cu}(0,z,0) \begin{pmatrix} 0 & 0 \\ 1 & 0 \end{pmatrix} \overline{U}_j(z)\right\rangle_X \rmd z \\ & = &
\int_\infty^0 \left\langle \Phi_+^\mathrm{cu}(0,z,0)^T \Psi_i(0), \begin{pmatrix} 0 & 0 \\ 1 & 0 \end{pmatrix} \overline{U}_j(z)\right\rangle_X \rmd z \\ & = &
\int_\infty^0 \left\langle \Psi_i(z),\begin{pmatrix} 0 & 0 \\ 1 & 0 \end{pmatrix} \overline{U}_j(z) \right\rangle_X \rmd z
\\ & = &
\int_\infty^0 \langle \psi_i(z),u_j(z) \rangle_{L^2(S^1)}\,\rmd z,
\end{eqnarray*}
where we have used the fact that $\Phi_+^\mathrm{cu}(0,z,0)^T\Psi_i(0)=\Psi_i(z)$; see \cite[Lemma~5.1]{SandstedeScheel01}. Proceeding in the same fashion for $E_-^\mathrm{u}(\sigma)$, we find that
\[
\langle \Psi_i(0),h^\psi(0) \overline{U}_j(0) \rangle_X = \int_{-\infty}^0 \langle \psi_i(z),u_j(z) \rangle_{L^2(S^1)}\,\rmd z,
\]
and therefore
\[
M(0) = h^\psi(0)-g^\psi(0) = \left( \int_{-\infty}^\infty \langle \psi_i(z),u_j(z) \rangle_{L^2(S^1)}\, \rmd z \right)_{ij}.
\]
Furthermore, since $\psi_1(x,t) = \psi_1 \in \mathbb{R}^n$, we have
\[
\int_\mathbb{R} \int_0^{2\pi} \langle \psi_1(x,t), \bar{u}_t(x,t) \rangle_{\mathbb{R}^N}\, \rmd t\, \rmd x = \int_\mathbb{R} \left\langle \psi_1, \int_0^{2\pi}\bar{u}_t(x,t)\, \rmd t \right\rangle_{\mathbb{R}^N}\, \rmd x = 0,
\]
which proves that $M(0)$ is lower triangular. Hypotheses~(S3) and (S4) imply that the diagonal entries of $M(0)$ are nonzero, so that $M(0)$ is invertible. 

This completes the proof for $h_+(\sigma)$. The proof for $h_-(\sigma)$ is analogous, but we need to be careful with the signs: The integral representations of the solutions $V(x)$ given above stay the same, but the roles of $h^\psi$ and $g^\psi$ are reversed. Thus, the matrix that appears in the representation of $h_\mathrm{p}^-(\sigma)$ is $-\tilde{M}_0$, and not $\tilde{M}_0$.
\end{Proof}

We define
\begin{equation}\label{E:tilde_dich}
\begin{array}{rclcl}
\tilde{P}_+^\mathrm{s}(x,\sigma) & := & P_+^\mathrm{s}(x,\sigma) - \Phi_+^\mathrm{s}(x,0,\sigma) h_+(\sigma) \Phi_+^\mathrm{u}(0,x,\sigma) & & x\geq0 \\[1ex]
\tilde{\Phi}_+^\mathrm{s}(x,y,\sigma) & := & \Phi_+^\mathrm{s}(x,y,\sigma) \tilde{P}_+^\mathrm{s}(y, \sigma) & & x\geq y\geq0 \\[1ex]
\tilde{\Phi}_+^\mathrm{u}(x,y,\sigma) & := & (1-\tilde{P}_+^\mathrm{s}(x,\sigma)) \Phi_+^\mathrm{u}(x,y,\sigma) & & y\geq x\geq0
\end{array}
\end{equation}
and similarly for $x,y \leq 0$. As in \cite[(3.20)]{PeterhofSandstedeScheel97}, equation (\ref{E:tilde_dich}) defines an exponential dichotomy on $\mathbb{R}^+$ with projection $\tilde{P}_+^\mathrm{s}(x,\sigma)$, since $\Rg P_+^\mathrm{s}(x,\sigma)=\Rg\tilde{P}_+^\mathrm{s}(x,\sigma)$ and $(1-\tilde{P}_+^\mathrm{s}(x,\sigma))(1-P_+^\mathrm{s}(x,\sigma))=1-\tilde{P}_+^\mathrm{s}(x,\sigma)$ for all $\sigma$. Similarly, $\tilde{P}_-^\mathrm{u}(x,\sigma)$ is the projection for the exponential dichotomy $\tilde{\Phi}_-^\mathrm{u}(x,y,\sigma)$ and $\tilde{\Phi}_-^\mathrm{s}(x,y,\sigma)$ on $\mathbb{R}^-$. By construction
\begin{equation}\label{E:tilde_proj}
\tilde{P}_+^\mathrm{s}(0,\sigma) = 1-\tilde{P}_-^\mathrm{u}(0,\sigma) \qquad \forall\;\sigma\neq0,
\end{equation}
and the Laurent series of these two operators coincide at $\sigma=0$, since the contribution of the pole at $\sigma=0$ is, in both cases, given by the matrix $\tilde{M}_0$. Hence, we can define a meromorphic exponential dichotomy on $\mathbb{R}$ via
\begin{equation}\label{e:med}
\Phi^\mathrm{s}(x,y, \sigma) = \left\{\begin{array}{lcl}
\tilde{\Phi}^\mathrm{s}_+(x,y,\sigma) & & x>y\geq0 \\[1ex]
\tilde{\Phi}_+^\mathrm{s}(x,0,\sigma)\tilde{\Phi}_-^\mathrm{s}(0,y, \sigma) & & x\geq0>y \\[1ex]
\tilde{\Phi}^\mathrm{s}_-(x,y,\sigma) & & 0>x>y 
\end{array}\right.
\end{equation}
for $x>y$, and an analogous expression for $\Phi^\mathrm{u}(x,y,\sigma)$ for $x<y$. Note that, if we fix $x>0$ and let $y\to0$, we obtain from the first two equations in (\ref{e:med}) the two expressions
\[
\Phi^\mathrm{s}(x,0,\sigma) = \begin{cases}
\tilde{\Phi}^\mathrm{s}_+(x,0,\sigma) \\
\tilde{\Phi}_+^\mathrm{s}(x,0,\sigma) \tilde{\Phi}_-^\mathrm{s}(0,0,\sigma)
= \tilde{\Phi}_+^\mathrm{s}(x,0,\sigma) (1-\tilde{P}_-^\mathrm{u}(0,\sigma))
= \tilde{\Phi}_+^\mathrm{s}(x,0,\sigma),
\end{cases}
\]
which coincide due to (\ref{E:tilde_proj}). This completes the meromorphic extension of the exponential dichotomies on $\mathbb{R}$ for $\sigma\in B_\epsilon(0)$.

It remains to obtain pointwise bounds for these operators. It suffices to consider the case $x>y$, as the case $x<y$ is completely analogous. Recall that $h_\pm(\sigma)=h_\mathrm{a}^\pm(\sigma)+h_\mathrm{p}^\pm(\sigma)$, where $h_\mathrm{a}^\pm(\sigma)$ is analytic in $\sigma$, while
\[
h_\mathrm{p}^\pm(\sigma): \quad Y_m \longrightarrow Y_m, \quad
V \longmapsto \frac{\pm1}{\sigma} \sum_{i,j=1}^2 m_{ij}^0 \langle\Psi_j(0),V\rangle_X \overline{U}_i(0),
\]
where $m_{ij}^0$ denote the entries of $\tilde{M}_0$.

We have
\[
\begin{array}{rclcl}
\Phi^\mathrm{s}(x,y,\sigma) & = &
\underbrace{\Phi_+^\mathrm{s}(x,y, \sigma) - \Phi_+^\mathrm{s}(x,0,\sigma)h^+_\mathrm{a}(\sigma)\Phi_+^\mathrm{u}(0,y,\sigma)}_{\text{(a)}}
- \underbrace{\Phi_+^\mathrm{s}(x,0,\sigma)h^+_\mathrm{p}(\sigma)\Phi_+^\mathrm{u}(0,y,\sigma)}_{\text{(i)}}
&& x>y\geq0 \\[2.5em]
\Phi^\mathrm{s}(x,y, \sigma) & = &
\underbrace{\Phi_+^\mathrm{s}(x,0, \sigma)\Phi_-^\mathrm{s}(0,y, \sigma) - \Phi_+^\mathrm{s}(x,0,\sigma)h^+_\mathrm{a}(\sigma)\Phi_+^\mathrm{u}(0,0,\sigma)\Phi_-^\mathrm{s}(0,y,\sigma)}_{\text{(b)}} && x>0>y \\[2.5em] && \qquad
- \underbrace{\Phi_+^\mathrm{s}(x,0,\sigma)h^+_\mathrm{p}(\sigma)\Phi_+^\mathrm{u}(0,0,\sigma)\Phi_-^\mathrm{s}(0,y,\sigma)}_{\text{(ii)}}
\\[2.5em]
\Phi^\mathrm{s}(x,y,\sigma) & = &
\underbrace{\Phi_-^\mathrm{s}(x,y, \sigma) + \Phi_-^\mathrm{u}(x,0,\sigma)h^-_\mathrm{a}(\sigma)\Phi_-^\mathrm{s}(0,y,\sigma)}_{\text{(c)}} + \underbrace{\Phi_-^\mathrm{u}(x,0,\sigma)h^-_\mathrm{p}(\sigma)\Phi_-^\mathrm{s}(0,y,\sigma)}_{\text{(iii)}}
&& 0>x>y. \\[2.5em] 
\end{array}
\]
We now use the bounds (\ref{ed:bdd}) for the exponential dichotomies on $\mathbb{R}^\pm$ that we derived in Lemma~\ref{l:aed}. First, we estimate the terms given by (a)-(c). For case~(a) with $x>y\geq0$, we find
\begin{eqnarray*}
\|\Phi^\mathrm{s}(x,y,\sigma)\|_{L(Y_m)} & = &
\|\Phi_+^\mathrm{s}(x,y,\sigma) - \Phi_+^\mathrm{s}(x,0,\sigma) h_\mathrm{a}^+(\sigma) \Phi_+^\mathrm{u}(0,y,\sigma)\|_{L(Y_m)}
\nonumber \\  & \leq &
K \sum_{\nu^+_\mathrm{out},\nu^+_\mathrm{in}} \left[ \rme^{\nu^+_\mathrm{out}(\sigma)(x-y)} + \rme^{\nu^+_\mathrm{out}(\sigma)x} \rme^{-\nu^+_\mathrm{in}(\sigma)y}\right].
\end{eqnarray*}
Case~(b) for $x\geq0>y$ gives
\begin{eqnarray*}
\|\Phi^\mathrm{s}(x,y,\sigma)\|_{L(Y_m)} & = &
\|\Phi_+^\mathrm{s}(x,0, \sigma)\Phi_-^\mathrm{s}(0,y, \sigma) - \Phi_+^\mathrm{s}(x,0,\sigma)h^+_\mathrm{a}(\sigma)\Phi_+^\mathrm{u}(0,0,\sigma)\Phi_-^\mathrm{s}(0,y,\sigma)\|_{L(Y_m)}
\nonumber \\ & \leq &
K \sum_{\nu^+_\mathrm{out},\nu^-_\mathrm{in}} \rme^{\nu^+_\mathrm{out}(\sigma)x} \rme^{-\nu^-_\mathrm{in}(\sigma)y}.
\end{eqnarray*}
Finally, case~(c) with $0>x>y$ gives
\begin{eqnarray*}
\|\Phi^\mathrm{s}(x,y,\sigma)\|_{L(Y_m)} & = &
\|\Phi_-^\mathrm{s}(x,y, \sigma) + \Phi_-^\mathrm{u}(x,0,\sigma)h^-_\mathrm{a}(\sigma)\Phi_-^\mathrm{s}(0,y,\sigma)\|_{L(Y_m)}
\nonumber \\ & \leq &
K \sum_{\nu^-_\mathrm{out},\nu^-_\mathrm{in}} \left[ \rme^{\nu^-_\mathrm{in}(\sigma)(x-y)} + \rme^{\nu^-_\mathrm{out}(\sigma)x} \rme^{-\nu^-_\mathrm{in}(\sigma)y}\right].
\end{eqnarray*}
Next, we consider the meromorphic terms (i)-(iii). For case~(i) with $x>y\geq0$, we obtain
\begin{eqnarray*}
\lefteqn{\Phi_+^\mathrm{s}(x,0,\sigma) h^+_\mathrm{p}(\sigma) \Phi_+^\mathrm{u}(0,y,\sigma)}
\nonumber \\ \nonumber & = &
\frac{1}{\sigma} \Phi_+^\mathrm{s}(x,0,\sigma) \sum_{i,j} m_{ij}^0 \overline{U}_i(0) \langle\Psi_j(0), \Phi_+^\mathrm{u}(0,y,\sigma) \cdot \rangle_X
\\ \nonumber & = &
\frac{1}{\sigma} \Phi_+^\mathrm{s}(x,0,\sigma) \left[ P_+^\mathrm{ss}(0,0)-P_+^\mathrm{ss}(0,\sigma)+P_+^\mathrm{ss}(0,\sigma) \right] \sum_{i,j} m_{ij}^0 \overline{U}_i(0) \langle\underbrace{\Phi_+^\mathrm{u}(0,y,\sigma)^T \Psi_j(0)}_{=:\Psi_j(y,\sigma)},\cdot\rangle_X
\\ \nonumber & = &
\frac{1}{\sigma} \Phi_+^\mathrm{s}(x,0,\sigma) [P_+^\mathrm{ss}(0,\sigma)+\mathrm{O}(\sigma)] \sum_{i,j} m_{ij}^0 \overline{U}_i(0) \langle \Psi_j(y,\sigma), \cdot \rangle_X
\\ \nonumber & = &
\frac{1}{\sigma} \Phi_+^\mathrm{s}(x,0,\sigma) P_+^\mathrm{ss}(0,\sigma) \sum_{i,j} m_{ij}^0 \overline{U}_i(0) \langle \Psi_j(y,\sigma), \cdot \rangle_X
+ \underbrace{\Phi_+^\mathrm{s}(x,0,\sigma) \mathrm{O}(1) \sum_{i,j} m_{ij}^0 \overline{U}_i(0) \langle \Psi_j(y,\sigma), \cdot \rangle_X}_{=:\text{ rest}}
\\ \nonumber & = &
\sum_{i,j} \left[ \frac{1}{\sigma} \Phi_+^\mathrm{s}(x,0,0) P_+^\mathrm{ss}(0,0) m_{ij}^0 \overline{U}_i(0) + \mathrm{O}(\rme^{-3\eta x})\right] \langle \Psi_j(y,\sigma),\cdot \rangle_X + \mbox{rest}
\\ & = &
\sum_{i,j} \left[ \frac{1}{\sigma} m_{ij}^0 \overline{U}_i(x) + \mathrm{O}(\rme^{-3\eta x})\right] \langle \Psi_j(y,\sigma),\cdot \rangle_X + \mbox{rest}.
\end{eqnarray*}
Note that the terms denoted by the Landau symbol $\mathrm{O}$ are all analytic in $\sigma$. Equation (\ref{ed:bdd}) shows that $\Psi_j(y,\sigma)=\Phi^\mathrm{u}_+(0,y,\sigma)^T\Psi_j(0)$ satisfies the bound
\begin{equation}\label{e:psi>0}
|\Psi_j(y,\sigma)| \leq K \sum_{\nu_\mathrm{in}^+} \rme^{-\nu_\mathrm{in}^+(\sigma)y}, \qquad
y\geq0.
\end{equation}
We remark that, as a consequence of \cite[Lemma~5.1]{SandstedeScheel01}, $\Phi^\mathrm{u}_+(0,y,\sigma)^T$ is the stable evolution of the adjoint equation, and $\Psi_j(y,\sigma)=\Phi^\mathrm{u}_+(0,y,\sigma)^T\Psi_j(0)$ therefore satisfies the adjoint equation.

The term~(ii) for $x\geq0>y$ can be treated similar to (i), and we obtain
\begin{equation*}
\Phi_+^\mathrm{s}(x,0,\sigma) h^+_\mathrm{p}(\sigma) \Phi_-^\mathrm{s}(0,y,\sigma)
= \sum_{i,j} \left[ \frac{1}{\sigma} m_{ij}^0 \overline{U}_j(x) + \mathrm{O}(\rme^{-3\eta x})\right]\langle \Psi_j(y,\sigma),\cdot\rangle_X + \mbox{rest},
\end{equation*}
where
\begin{equation}\label{e:psi<0}
|\Psi_j(y,\sigma)| \leq K \sum_{\nu_\mathrm{in}^-} \rme^{-\nu_\mathrm{in}^-(\sigma)y}, \qquad
y\leq0.
\end{equation}
Finally, case~(iii) for $0>x>y$ gives
\begin{equation*}
\Phi_-^\mathrm{u}(x,0,\sigma) h^-_\mathrm{p}(\sigma) \Phi_-^\mathrm{s}(0,y,\sigma)
= \sum_{i,j} \left[ \frac{1}{\sigma} m_{ij}^0 \overline{U}_j(x) + \mathrm{O}(\rme^{-3\eta|x|})\right]\langle \Psi_j(y,\sigma),\cdot\rangle_X + \mbox{rest},
\end{equation*}
where $\Psi_j(y,\sigma)$ obeys the bound (\ref{e:psi<0}). We summarize our findings in the following lemma.

\begin{Lemma}\label{l:ed}
There exists an $\epsilon>0$ so that (\ref{n:evp}) has meromorphic exponential dichotomies $\Phi^\mathrm{s}(x,y,\sigma)$ and $\Phi^\mathrm{u}(x,y,\sigma)$, defined respectively for $x>y$ and $x<y$, for $\sigma\in B_\epsilon(0)$ such that
\begin{equation}\label{E:lem4_1}
\Phi^k(x,y,\sigma) = \frac{\mp1}{\sigma} \sum_{i,j} m_{ij}^0 \overline{U}_i(x) \langle\Psi_j(y,\sigma),\cdot\rangle_X + \tilde{\Phi}^k(x,y,\sigma), \qquad k=\mathrm{s,u},
\end{equation}
(the minus and plus signs are for $\Phi^\mathrm{s}$ and $\Phi^\mathrm{u}$, respectively), where $\Psi_j(y,\sigma)$ and $\tilde{\Phi}^k(x,y,\sigma)$ are analytic in $\sigma$. Furthermore, $\Psi(y,\sigma)$ obeys the bounds (\ref{e:psi>0}) and (\ref{e:psi<0}), while we have the bounds
\begin{equation}\label{E:lem4_2}
\|\tilde{\Phi}^s(x,y,\sigma)\|_{L(Y_m)} \leq K \left\{\begin{array}{lcl}
\sum_{\nu^+_\mathrm{out},\nu^+_\mathrm{in}} \left[ \rme^{\nu^+_\mathrm{out}(\sigma)(x-y)} + \rme^{\nu^+_\mathrm{out}(\sigma)x} \rme^{-\nu^+_\mathrm{in}(\sigma)y}\right]
&& 0\leq y<x \\[1em]
\sum_{\nu^-_\mathrm{in},\nu^+_\mathrm{out}} \rme^{\nu^+_\mathrm{out}(\sigma)x} \rme^{-\nu^-_\mathrm{in}(\sigma)y}
&& y<0\leq x \\[1em]
\sum_{\nu^-_\mathrm{out},\nu^-_\mathrm{in}} \left[ \rme^{\nu^-_\mathrm{in}(\sigma)(x-y)} + \rme^{\nu^-_\mathrm{out}(\sigma)x} \rme^{-\nu^-_\mathrm{in}(\sigma)y}\right]
&& y<x<0
\end{array}\right.
\end{equation}
for $\tilde{\Phi}^\mathrm{s}(x,y,\sigma)$, with analogous bounds for $\tilde{\Phi}^\mathrm{u}(x,y,\sigma)$. Derivatives of $\tilde{\Phi}^k(x,y,\sigma)$ and $\Psi_j(y,\sigma)$ with respect to $x$ and $y$ satisfy the same estimates multiplied by $[|\sigma|+\rme^{-\eta|x|}]$ and $[|\sigma|+\rme^{-\eta|y|}]$, respectively, as operators in $L(Y_m,Y_{m-\frac12})$.
\end{Lemma}


\subsection{Pointwise bounds for the resolvent kernel}

We now use the asymptotic formulas for the exponential dichotomy given in Lemma~\ref{l:ed} to complete the proof of Theorem~\ref{resbounds}. Lemma~\ref{l:ed} shows that the meromorphic extension of the exponential dichotomy can be written as
\[
\Phi^k(x,y,\sigma) = \tilde{E}_1(x,y,\sigma) + \tilde{E}_2(x,y,\sigma) + \tilde{\Phi}^k(x,y,\sigma), \qquad
k=\mathrm{s,u},
\]
where the $\tilde{E}_{1,2}$ terms come from the meromorphic components of (\ref{E:lem4_1}). Thus, we may write the representation (\ref{E:G_dich_reln}) of $G_*$ as
\begin{eqnarray}
G_*(x,y,\sigma;s) & = &
P_1 \int_{-\infty}^x \left[\tilde{E}_1 + \tilde{E}_2 + \tilde{\Phi}^\mathrm{s}\right](x,z,\sigma) \mathcal{B}(z) V_\ell(z,y,\sigma;s) \,\rmd z
\label{E:G_dich_reln_2} \\ &&
+ P_1 \int_{\infty}^x \left[\tilde{E}_1 + \tilde{E}_2 + \tilde{\Phi}^\mathrm{u}\right](x,z,\sigma) \mathcal{B}(z) V_\ell(z,y,\sigma;s) \,\rmd z, \nonumber
\end{eqnarray}
where $\mathcal{B}$ is defined in (\ref{e:B}). We shall prove that the terms that involve $\tilde{\Phi}^\mathrm{s,u}$ may all be included in the term $R(x,y,\sigma;s)$ given in Theorem~\ref{resbounds}, while the terms involving $\tilde{E}_{1,2}$ lead directly to the terms $E_{1,2}(x,y,\sigma;s)$ of Theorem~\ref{resbounds}. Throughout, we will often suppress any dependence on $s$ for notational convenience.

We now consider (\ref{E:G_dich_reln_2}) and begin by estimating the integrals that involve $\tilde{\Phi}^\mathrm{s,u}$. Each of these integrals can be estimated in the same fashion, and we therefore give details only for the integral that contains $\tilde{\Phi}^\mathrm{s}$ in the case where $x\geq0>y$. For $x\geq0>y$, equation (\ref{E:lem4_2}) and Lemma~\ref{l:bs2} imply that
\begin{eqnarray}\label{e:phibdds}
\lefteqn{\left|\int_{-\infty}^x \tilde{\Phi}^\mathrm{s}(x,z,\sigma) \mathcal{B}(z) V_\ell(z,y,\sigma) \,\rmd z \right|_{Y_m}} \\ \nonumber
& \leq &
\int_{-\infty}^x \|\tilde{\Phi}^\mathrm{s}(x,z,\sigma)\|_{L(Y_m)} \rme^{-\eta|z-y|}\rme^{\eta|z-y|}|\mathcal{B}(z) V_\ell(z,y,\sigma)|_{Y_m} \,\rmd z
\\ \nonumber & \leq &
C \left( \int_{-\infty}^x \|\tilde{\Phi}^\mathrm{s}(x,z,\sigma)\|_{L(Y_m)}^2 \rme^{-2\eta|z-y|} \,\rmd z \right)^{1/2}
\left( \int_{-\infty}^x \rme^{2\eta|z-y|}|\mathcal{B}(z) V_\ell(z,y,\sigma)|_{Y_m}^2 \,\rmd z\right)^{1/2}
\\ \nonumber & \leq &
C  \left[ \sum_{\nu^-_\mathrm{in},\nu^+_\mathrm{out}} \rme^{\nu^+_\mathrm{out}(\sigma)x} \rme^{-\nu^-_\mathrm{in}(\sigma)y} \right]
\left|\rme^{\eta|\cdot-y|}\mathcal{B}(\cdot) V_\ell(\cdot,y,\sigma)\right|_{L^2(\mathbb{R},Y_m)}
\\ \nonumber & \leq &
C  \left[ \sum_{\nu^-_\mathrm{in},\nu^+_\mathrm{out}} \rme^{\nu^+_\mathrm{out}(\sigma)x} \rme^{-\nu^-_\mathrm{in}(\sigma)y} \right].
\end{eqnarray}
Derivatives of this integral with respect to $y$ can be estimated as in (\ref{e:lest}) below, using again (\ref{E:lem4_2}) and Lemma~\ref{l:bs2}. Derivatives with respect to $(x,t,s)$ lead to the same estimates as above, since they can be accounted for by estimating the derivatives of the dichotomies in $L(Y_m,Y_{m-\frac12})$ instead of $L(Y_m)$. In summary, the estimates for the terms $\tilde{\Phi}^\mathrm{s,u}$ in the exponential dichotomy transfer to the resolvent kernel and are captured by the term $R(x,y,\sigma;s)$ in Theorem~\ref{resbounds}.

Next, we consider the integrals in (\ref{E:G_dich_reln_2}) that involve the terms $\tilde{E}_j$. We focus on $\tilde{E}_1$ as the term $\tilde{E}_2$ can be treated similarly. Equation (\ref{E:lem4_1}) implies that the two integrals in (\ref{E:G_dich_reln_2}) that involve $\tilde{E}_1$ can be combined as follows:
\begin{eqnarray}
&&
P_1\left[ \int_{-\infty}^x \tilde{E}_1(x,z,\sigma) \mathcal{B}(z) V_\ell(z,y,\sigma;s) \,\rmd z \right](t) +
P_1\left[ \int_{\infty}^x \tilde{E}_1(x,z,\sigma) \mathcal{B}(z) V_\ell(z,y,\sigma;s) \,\rmd z \right](t) \nonumber \\ \label{e:E1} && \qquad\qquad =
-\frac{1}{\sigma} \bar{u}_x(x,t) \int_{-\infty}^\infty \sum_j m_{1j}^0 \left\langle\Psi_j(z,\sigma),\mathcal{B}(z) V_\ell(z,y,\sigma;s) \right\rangle_X \,\rmd z.
\end{eqnarray}
Assume first that $y\leq0$. Using (\ref{e:gapsoln}) and the definition of $\Psi_j(z,\sigma)$, we have
\begin{equation}\label{n:psi}
\Psi_j(z,\sigma) = \sum_{\nu_\mathrm{in}^-} c_{j,\mathrm{in}}(z,\sigma) \rme^{-\nu_\mathrm{in}^-(\sigma)z}, \qquad |\partial_z c_{j,\mathrm{in}}(z,\sigma)| \leq C\rme^{-\eta|z|}
\end{equation}
for $z\leq0$. Thus, we may write
\begin{eqnarray}\label{e:efe}
\lefteqn{ \int_{-\infty}^\infty \sum_j m_{1j}^0 \left\langle\Psi_j(z,\sigma),\mathcal{B}(z) V_\ell(z,y,\sigma;s) \right\rangle_X \,\rmd z } \\ \nonumber & = &
\sum_j m_{1j}^0 \left( \sum_{\nu_\mathrm{in}^-} \int_{-\infty}^0 \rme^{-\nu_\mathrm{in}^-(\sigma)z} \left\langle c_{j,\mathrm{in}}(z,\sigma),\mathcal{B}(z) V_\ell(z,y,\sigma;s) \right\rangle_X \,\rmd z + \int_0^\infty \left\langle\Psi_j(z,\sigma),\mathcal{B}(z) V_\ell(z,y,\sigma;s) \right\rangle_X \,\rmd z \right)
\\ \nonumber & = &
\sum_{\nu_\mathrm{in}^-} \rme^{-\nu_\mathrm{in}^-(\sigma)y} \sum_j m_{1j}^0 \left( \int_{-\infty}^0  \rme^{-\nu_\mathrm{in}^-(\sigma)(z-y)} \left\langle c_{j,\mathrm{in}}(z,\sigma),\mathcal{B}(z) V_\ell(z,y,\sigma;s) \right\rangle_X \,\rmd z
\right. \\ \nonumber && \left.
+ \frac{1}{p} \int_0^\infty \rme^{\nu_\mathrm{in}^-(\sigma)y} \left\langle\Psi_j(z,\sigma),\mathcal{B}(z) V_\ell(z,y,\sigma;s) \right\rangle_X \,\rmd z \right)
\\ \nonumber & =: &
\sum_{\nu_\mathrm{in}^-} \rme^{-\nu_\mathrm{in}^-(\sigma)y} l_{1,\mathrm{in}}^-(y,\sigma)^T,
\end{eqnarray}
where $l_{1,\mathrm{in}}^-(y,\sigma)\in\mathbb{R}^N$ for each fixed $(y,\sigma)$, so that (\ref{e:E1}) becomes
\[
-\frac{1}{\sigma} \bar{u}_x(x,t) \sum_{\nu_\mathrm{in}^-} \rme^{-\nu_\mathrm{in}^-(\sigma)y} l_{1,\mathrm{in}}^-(y,\sigma)^T.
\]
We now claim that there are positive constants $C$ and $\eta$ so that
\begin{equation}\label{e:lbdds}
|l_{1,\mathrm{in}}(y,\sigma)|\leq C, \qquad
|\partial_y l_{1,\mathrm{in}}(y,\sigma)|\leq C\left[|\sigma|+\rme^{-\eta|y|}\right]
\end{equation}
for all $y\leq0$ and all $\sigma$ with $\Re\sigma\geq-\epsilon$. Furthermore, the same bounds are then true for $s$-derivatives of $l_{1,\mathrm{in}}(y,\sigma)$, since these are equivalent to $t$-derivatives, which are taken care of by the regularity of functions in $Y_m$. This then yields the expressions (\ref{E1lambda}) in Theorem~\ref{resbounds} and the bounds (\ref{ljkbounds}) in Theorem~\ref{greenbounds}, since we can expand the analytic term $l_{1,\mathrm{in}}(y,\sigma)$ in $\sigma$ and subsume the $\sigma$-dependent part into $R(x,y,\sigma;s)$.

The first estimate in (\ref{e:lbdds}) follows from (\ref{e:efe}) and Lemma~\ref{l:bs2} as in (\ref{e:phibdds}) above, since $|\rme^{-\nu_\mathrm{in}^-(\sigma)(z-y)}|\leq C\rme^{\frac{\eta}{2}|z-y|}$. It remains to show that differentiation with respect to $y$ leads to exponential decay with respect to $y$. We estimate here only the most dangerous component in (\ref{e:efe}): Lemma~\ref{l:bs2} and (\ref{n:psi}) give
\begin{eqnarray}\label{e:lest}
\lefteqn{ \partial_y \int_{-\infty}^0  \rme^{-\nu_\mathrm{in}^-(\sigma)(z-y)} \left\langle c_{j,\mathrm{in}}(z,\sigma),\mathcal{B}(z) V_\ell(z,y,\sigma;s) \right\rangle_X \,\rmd z }
\\ \nonumber & = &
\int_{-\infty}^0  \rme^{-\nu_\mathrm{in}^-(\sigma)(z-y)} \left\langle c_{j,\mathrm{in}}(z,\sigma),\mathcal{B}(z) \partial_y V_\ell(z,y,\sigma;s) \right\rangle_X \,\rmd z
+ \mathrm{O}(|\sigma|)
\\ \nonumber & = &
\mathrm{O}\left( \int_{-\infty}^0 \left|\rme^{\nu_\mathrm{in}^-(\sigma)(z-y)}\right| \rme^{-\eta|z|} \rme^{-\eta|z-y|} \,\rmd z \right)
+ \mathrm{O}(|\sigma|+\rme^{-\eta|y|})
\\ \nonumber & = &
\mathrm{O}(|\sigma|+\rme^{-\frac{\eta}{4}|y|})
\end{eqnarray}
and renaming $\eta$ gives the second estimate in (\ref{e:lbdds}).


\section{Estimates of the Green's distribution and linear stability}\label{lin}

With the bounds on the resolvent kernel $G_*$ determined in Theorem~\ref{resbounds}, it is now possible to obtain the desired pointwise bounds on the Green's distribution $\mathcal{G}_*$ by a simplified version of the analysis of \cite{ZumbrunHoward98,MasciaZumbrun03}.


\subsection{Bounds for the Green's distribution}

\begin{figure}[t]
\centering\includegraphics{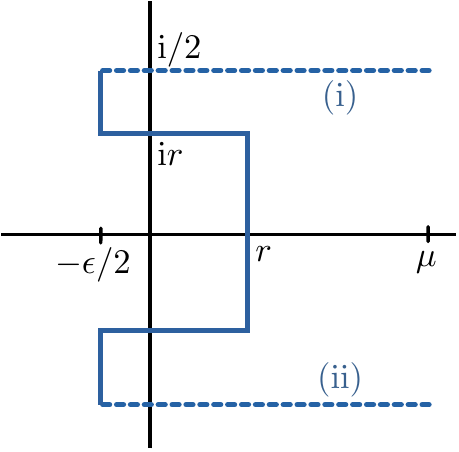}
\caption{Plotted is the contour $\tilde{\Gamma}$ (solid and dashed lines) that is used to calculate the Green's distribution $\mathcal{G}_*$ from the resolvent kernel $G_*$. The contributions from integrating along pieces (i) and (ii) cancel, yielding the contour $\Gamma$ (solid line). }
\label{F:contour}
\end{figure}

In this section, we prove Theorem~\ref{greenbounds}. Corollary~\ref{c:g*} states that the Green's distribution $\mathcal{G}_*$ is given as the contour integral
\begin{equation}\label{e:gdci}
\mathcal{G}_*(x,t;y,s) = \frac{1}{2\pi\rmi} \int_{\mu - \rmi/2}^{\mu + \rmi/2} \rme^{\sigma t} 
[G_*(x,y,\sigma;s)](t) \,\rmd\sigma
\end{equation}
that involves the resolvent kernel $G_*$. Applying Cauchy's integral theorem and Theorem~\ref{resbounds}, which, in particular, states that $G_*$ may be meromorphically extended into the region $\{\sigma\in\mathbb{C}:\;\Re\sigma>-\epsilon\}$, we first observe that (\ref{e:gdci}) may be modified to
\[
\mathcal{G}_*(x,t;y,s) = \frac{1}{2\pi\rmi} \oint_{\tilde{\Gamma}} \rme^{\sigma t} [G_*(x,y,\sigma;s)](t) \,\rmd\sigma,
\]
where $\tilde{\Gamma}$ is defined in Figure~\ref{F:contour} for some small constant $r>0$. The fact that the key relation \eqref{keyper} persists under analytic extension implies that the integration along the top and bottom pieces, $[\mu -\rmi/2, -\epsilon/2 -\rmi/2]$ and $[-\epsilon/2 +\rmi/2, \mu +\rmi/2]$, cancel; see Figure~\ref{F:contour}. We therefore have
\[
\mathcal{G}_*(x,t;y,s) = \frac{1}{2\pi\rmi} \oint_{\Gamma} \rme^{\sigma t} [G_*(x,y,\sigma;s)](t) \,\rmd\sigma,
\]
where 
\[
\Gamma := \left[-\frac\epsilon2 -\frac\rmi2,-\frac\epsilon2 -\rmi r\right] \cup 
\left[-\frac\epsilon2 -\rmi r, r-\rmi r\right] \cup 
[r -\rmi r, r+\rmi r] \cup 
\left[r+\rmi r, -\frac\epsilon2 +\rmi r\right] \cup
\left[-\frac\epsilon2 +\rmi r, -\frac\epsilon2 +\frac\rmi2\right].
\]
This is a representation on a contour that corresponds exactly to the low-frequency part of the contour used to begin the arguments of \cite[\S8]{ZumbrunHoward98} and \cite[\S7]{MasciaZumbrun03}, and we may therefore move the contours for the individual meromorphic pieces $E_1$, $E_2$ and $R$ of $G_*$ exactly as in these references. Since our resolvent bounds for $G_*$ as well as the initial contour $\Gamma$ are the same as for the low-frequency estimates in the time-independent case treated there, we obtain the same bounds for $\mathcal{G}_*$ that were obtained in \cite{ZumbrunHoward98,MasciaZumbrun03} for the entire Green's function in the time-independent case. Note that the new term $E_2$ that arises in $G_*$ for the time-dependent case has exactly the same form as the term $E_1$, except for the factor $\partial_t \bar u$ in place of $\partial_x \bar u$, so it may again be treated in exactly the same way as before. We omit the contour estimates and refer the reader instead to \cite[\S8]{ZumbrunHoward98} and \cite[\S7]{MasciaZumbrun03} for details. This establishes the estimates \eqref{E1}--\eqref{ljkbounds} of Theorem~\ref{greenbounds}.

Finally, we recall that the singular parts $\mathcal{G}_0$, $\mathcal{G}_1$, and $\mathcal{G}_2$ of the Green's function have already been estimated in Theorem~\ref{t:gd}. This completes the proof of Theorem \ref{greenbounds}.
\qed

We remark that the contour estimates carried out in the present time-dependent case are in fact somewhat simpler than those carried out in \cite{ZumbrunHoward98,MasciaZumbrun03} for the time-independent case. The reason is that difficulties associated with high frequencies have been subsumed into the iterative parametrix-type construction of the resolvent kernel and were dealt with in Lemma~\ref{parametrixlem} --- an illustration of conservation of difficulty. Though we will not use this, we remark that the bounds we obtain for $G_*$ are nonsingular as $t\to0^+$ and therefore somewhat better than those that hold for the low-frequency part of the time-independent Green's function. We shall use this observation later in the proof of Proposition~\ref{parambds}.


\subsection{Proof of Theorem \ref{stabcrit}}

From the bounds of Theorem~\ref{greenbounds}, we obtain the equivalence of spectral and linearized stability as follows.

As in the time-independent case, the pointwise bounds of Theorem~\ref{greenbounds} yield linearized orbital stability in $L^1\cap H^\mathrm{s}$, for any $s$, by standard $L^1\to L^p$ convolution bounds; see \cite[\S 9]{ZumbrunHoward98} and \cite[\S8]{MasciaZumbrun03}. This gives the sufficiency of spectral stability. Necessity follows from the theory of effective spectral projections, and the spectral expansions of the pointwise Green's function in terms of the effective spectral projections for $(x,y)$ restricted to a bounded domain, that was developed in \cite[\S9]{ZumbrunHoward98} and \cite[\S8]{MasciaZumbrun03}. We omit the details as they are essentially the same as in the time-independent case of those papers.

This concludes the proof of Theorem~\ref{stabcrit}.\qed


\section{Nonlinear stability}\label{nonlinsec}

We establish nonlinear stability by a modified version of the arguments used in \cite{RaoofiZumbrun07,HowardZumbrun06} to prove stability of general-type shock waves in the time-independent case. First, we will separate out the motion of perturbations along spatial and temporal translates of the shock profile from the rest of the time evolution by seeking solutions $u(x,t)$ in the form
\[
u(x,t) = \bar{u}(x-q(t),t-\tau(t)) + v(x,t), \qquad
\zeta(t) := (q,\tau)(t).
\]
The decomposition of the Green's distribution into the terms $\mathcal{E}_j$, which capture the movement of perturbations along the shock profile, and the remainder $\tilde{\mathcal{G}}$, which describes the evolution along the characteristics, allows us to seek $(\zeta,v)$ as fixed points of an appropriate set of integral equations that are similar to equations (\ref{nl:q}) and (\ref{nl:v}). An iteration scheme is then set up to construct solutions to these equations. To prove temporal decay, the anticipated spatio-temporal decay properties of $(\zeta,v)$ are encoded in an appropriate set of template functions, which are used to control the norms of iterates.


\subsection{Fixed-point iteration scheme}

We now introduce the fixed-point iteration scheme by which we shall simultaneously construct and estimate the solution of the perturbed shock problem. We will use the notation
\[
\zeta = (q, \tau), \quad
\zeta_* = (q_*, \tau_*), \quad
\bar{u}^\zeta(x,t) = \bar{u}(x-q, t-\tau), \quad
\frac{d\bar{u}^\zeta}{d\zeta}(x,t) = -(\bar{u}_x,\bar{u}_t)(x-q, t-\tau)
\]
and
\begin{eqnarray*}
\frac{\partial\bar{u}^\zeta}{\partial\zeta}\Big|_{\zeta_*} \dot{\zeta} & = & - (\bar{u}_x(x-q_*, t-\tau_*), \bar{u}_t(x-q_*, t-\tau_*)) \cdot (\dot{q}(t), \dot{\tau}(t)) \\
& = & -\bar{u}_x(x-q_*, t-\tau_*)\dot{q}(t) - \bar{u}_t(x-q_*, t-\tau_*)\dot{\tau}(t).
\end{eqnarray*}
We fix $\zeta=\zeta_*$ and consider the linearization about $u^{\zeta_*}$. In particular,
\[
\pi(y,s,t) = (\pi_1,\pi_2)(y,s,t)
\]
will depend on $\zeta^*$ through the profile $u^{\zeta_*}$ about which we linearize, though we will suppress this dependence in our notation. For later use, we note the important fact
\begin{equation}\label{efacts}
\pi(y,s,s)\equiv 0.
\end{equation}
Our starting point, similar to \cite{HowardZumbrun06,RaoofiZumbrun07}, is the observation that, if $u$ solves (\ref{main}), then the perturbation $v$ defined via
\[
u(x,t) = \bar{u}^{\zeta_*+\zeta(t)}(x,t) + v(x,t)
\]
satisfies
\begin{equation}\label{mixedperteq}
v_t-\mathcal{L}^{\zeta_*} v = [Q^{\zeta_*}(\zeta,v)]_x
+ \frac{\partial \bar{u}^\zeta}{\partial \zeta}\Big|_{\zeta_*}\dot \zeta (t) 
+ [R^{\zeta_*}(\zeta, v)]_x + S^{\zeta_*}(\zeta, \zeta_t),
\end{equation}
where
\begin{eqnarray*}
\mathcal{L}^{\zeta_*} v & := & v_{xx} - [f_u(\bar{u}^{\zeta_*})v]_x \\
Q^{\zeta_*}(\zeta, v) & := &
\Big(f(\bar{u}^{\zeta_*+\zeta})+ f_u(\bar{u}^{\zeta_*+\zeta})v
-f(\bar{u}^{\zeta_*+\zeta}+v) \Big) =\mathrm{O}(|v|^2) \\
R^{\zeta_*}(\zeta, v) & := &
\Big(f_u(\bar{u}^{\zeta_*}(x,t))-f_u(\bar{u}^{\zeta_*+\zeta}(x,t))\Big)v
=\mathrm{O}( \rme^{-\eta|x|}|\zeta| |v| ) \\
S^{\zeta_*}(\zeta, \dot{\zeta}) & := &
\left( \frac{\partial\bar{u}^\zeta}{\partial \zeta}\Big|_{\zeta_*+\zeta(t)}-\frac{\partial \bar{u}^\zeta}{\partial \zeta}\Big|_{\zeta_*}
\right) \dot \zeta
= \mathrm{O}( \rme^{-\eta|x|}|\dot \zeta| |\zeta| ).
\end{eqnarray*}
The above asymptotics hold as long as $|v|$ remains bounded,
by Taylor's Theorem together with exponential decay of $\partial_\zeta \bar{u}^\zeta$ in space, which follows from Hypothesis~\ref{H1}. Note that, in the above, we have suppressed the time-dependence of $\zeta =\zeta(t)$ in the superscripts for notational convenience.

We can now describe the iteration scheme for $(\zeta_*,\zeta,v)$. For given $\zeta_*^{n-1}$ and $\zeta^{n-1}(\cdot)$, let
\begin{equation}\label{inidat}
v_0^{n-1}(x) := u_0(x)-\bar{u}^{\zeta_*^{n-1}}(x,0)
\end{equation}
and define $v^n$ to be the solution of
\begin{eqnarray}
v^n(x,t) & = & \int_{-\infty}^{\infty} \tilde{\mathcal{G}}^{n-1}(x, t; y,0) v^{n-1}_0(y)\, \rmd y + \int_0^t\int_{-\infty}^{\infty}  \tilde{\mathcal{G}}^{n-1} (x, t; y,s)  S^{\zeta_*^{n-1}}(\zeta^{n-1},\dot\zeta^{n-1})(y,s)\,\rmd y\,\rmd s \nonumber \\ \label{un} & &
-\int_0^t\int_{-\infty}^{\infty} \tilde{\mathcal{G}}^{n-1}_y(x,
t; y,s) \left[ Q^{\zeta_*^{n-1}}(\zeta^{n-1},v^{n})+R^{\zeta_*^{n-1}}(\zeta^{n-1},v^{n})\right](y,s)\,\rmd y\,\rmd s,
\end{eqnarray}
where $ \tilde{\mathcal{G}}^{n-1}$ is the decaying part of the Green's distribution
$\mathcal{G}^{n-1} = \mathcal{E}_1^{n-1} + \mathcal{E}_2^{n-1} + \mathcal{\tilde{G}}^{n-1}$ of the linearized equation around $\bar{u}^{\zeta_*^{n-1}}$. Further, set
\begin{eqnarray}\label{zetan}
\zeta^n(t) & := & -\int^\infty_{-\infty} (\pi^{n-1}(y,0,t)-\pi^{n-1}(y,0,\infty)) v_0^{n-1}(y)\,\rmd y \\ \nonumber &&
-\int^t_0\int^{\infty}_{-\infty} (\pi^{n-1}(y,s,t)-\pi^{n-1}(y,s,\infty))
S^{\zeta_*^{n-1}}(\zeta^{n-1},\dot\zeta^{n-1})(y,s)\,\rmd y\,\rmd s \\ \nonumber &&
+\int^t_0\int^{\infty}_{-\infty} (\pi^{n-1}_{y}(y,s,t)-\pi^{n-1}_y(y,s,\infty))\left[ Q^{\zeta_*^{n-1}}(\zeta^{n-1},v^{n})+ R^{\zeta_*^{n-1}}(\zeta^{n-1}, v^{n})\right](y,s)\,\rmd y\,\rmd s \\ \nonumber &&
+\int_t^{\infty}\int^{\infty}_{-\infty} \pi^{n-1}(y,s,\infty)
S^{\zeta_*^{n-1}}(\zeta^{n-1},\dot\zeta^{n-1}) (y,s)\,\rmd y\,\rmd s, \\ \nonumber &&
-\int_t^{\infty}\int^{\infty}_{-\infty} \pi^{n-1}_y(y,s,\infty) \left[ Q^{\zeta_*^{n-1}}(\zeta^{n-1},v^{n})+R^{\zeta_*^{n-1}}
(\zeta^{n-1}, v^{n})\right] (y,s)\,\rmd y\,\rmd s
\end{eqnarray}
and
\begin{eqnarray}\label{zetasn}
\zeta_*^n & := & \zeta_*^{n-1}
-\int_{-\infty}^{\infty}\pi^{n-1}(y,0,\infty) v^{{n-1}}_0(y)\,\rmd y \\ \nonumber &&
-\int_0^{\infty}\int_{-\infty}^{\infty} \pi^{n-1}(y,s, \infty) S^{\zeta_*^{n-1}}(\zeta^{n-1},\dot\zeta^{n-1})(y,s)\,\rmd y\,\rmd s \\ \nonumber &&
+\int_0^{\infty}\int_{-\infty}^{\infty}\pi_y^{n-1}(y,s,
\infty) \left[ Q^{\zeta_*^{n-1}}(\zeta^{n-1},v^n)+R^{\zeta_*^{n-1}}(\zeta^{n-1},
v^n)\right](y,s)\,\rmd y\,\rmd s,
\end{eqnarray}
where $\pi^{n-1} = (\pi_1^{n-1}, \pi_2^{n-1})$ is associated with $\mathcal{E}_1^{n-1}$ and $\mathcal{E}_2^{n-1}$ of the linearized operator at $\bar{u}^{\zeta_*^{n-1}}$. Formally, we can then define an associated iteration map $\mathcal{T}$ by
\[
(\zeta^n, \zeta_*^n)=
\mathcal{T} (\zeta^{n-1}, \zeta_*^{n-1}).
\]
We now clarify the relation between fixed points of $\mathcal{T}$ and solutions of the viscous system of conservation laws.

\begin{Lemma}\label{soln}
Under \eqref{un}--\eqref{zetasn}, the function
$u^n := \bar{u}^{\zeta_*^{n-1}+\zeta^{n-1}} +v^n$
satisfies the equation
\begin{equation}\label{unequation}
u^n_t + f(u^n)_x - u^n_{xx} =
\frac{\partial \bar{u}^\zeta}{\partial \zeta}\Big|_{\zeta_*^{n-1}}
\big(\dot\zeta^n(t) -\dot\zeta^{n-1}(t)\big)
\end{equation}
with initial data
$u^n(\cdot, 0)=
u_0(\cdot) + (\bar{u}^{\zeta_*^{n-1}+\zeta^{n-1}(0)}(\cdot, 0) -\bar{u}^{\zeta_*^{n-1}}(\cdot, 0))$, where
\begin{eqnarray}\dot{q}^n(t) & = &
-\int^\infty_{-\infty} \partial_t \pi^{n-1}_1 (y,0,t)
v_0^{n-1}(y)\,\rmd y  - \int^t_0\int^{\infty}_{-\infty} \partial_t \pi^{n-1}_1 (y,s,t)
S^{\zeta_*^{n-1}}(\zeta^{n-1},\dot\zeta^{n-1})(y,s)\,\rmd y\,\rmd s
\nonumber \\ \label{dotdeltan} &&
+\int^t_0\int^{\infty}_{-\infty} \partial_t \partial_y \pi^{n-1}_1(y,s,t)
\left[Q^{\zeta_*^{n-1}}(\zeta^{n-1},v^n)+ R^{\zeta_*^{n-1}}(\zeta^{n-1}, v^n)\right](y,s)\,\rmd y\,\rmd s
\end{eqnarray}
and
\begin{eqnarray}
\dot \tau^n (t) & = &
-\int^\infty_{-\infty} \partial_t \pi^{n-1}_2 (y,0,t)
v_0^{n-1}(y)\,\rmd y  - \int^t_0\int^{\infty}_{-\infty} \partial_t \pi^{n-1}_2 (y,s,t)
S^{\zeta_*^{n-1}}(\zeta^{n-1},\dot\zeta^{n-1})(y,s)\,\rmd y\,\rmd s
\nonumber \\ \label{dottaun} &&
+\int^t_0\int^{\infty}_{-\infty} \partial_t \partial_y \pi^{n-1}_2(y,s,t)
\left[Q^{\zeta_*^{n-1}}(\zeta^{n-1},v^n)+ R^{\zeta_*^{n-1}}(\zeta^{n-1}, v^n)\right](y,s)\,\rmd y\,\rmd s.
\end{eqnarray}
In particular, $u^n$ satisfies
\eqref{main}
with initial data $u_0$ if and only if
\begin{equation*}
(\zeta^n, \zeta_*^n)= (\zeta^{n-1}, \zeta_*^{n-1}),
\end{equation*}
that is, if and only if $(\zeta^n, \zeta_*^n)$ is a fixed point of $\mathcal{T}$, in which case also $\zeta(0)=\zeta(\infty)=0$.
\end{Lemma}

\begin{Proof}
Equations \eqref{dotdeltan} and \eqref{dottaun} follow by differentiation of
\eqref{zetan}, recalling property \eqref{efacts} and the related fact that $\int_\mathbb{R} \pi_y(y,t,t)[Q,R](y)\,\rmd y=0$ for functions $Q$ and $R$ of the type considered here; see \cite[\S4.2.4]{Zumbrun04} for more details. From \eqref{zetan}, we find that $\zeta^n(\infty)=0$ and
\begin{eqnarray*}
\zeta^n(t)-\zeta^n(0) & = &-\int_{-\infty}^{\infty}\pi^{n-1}(y,0,t) v^{{n-1}}_0(y)\,\rmd y
-\int_0^t\int_{-\infty}^{\infty} \pi^{n-1}(y,s, t) S^{\zeta_*^{n-1}}
(\zeta^{n-1},\dot\zeta^{n-1})(y,s)\,\rmd y\,\rmd s \\ &&
+\int_0^t\int_{-\infty}^{\infty}\pi^{n-1}_y(y,s, t)
\left[Q^{\zeta_*^{n-1}}(\zeta^{n-1}, v^n)+R^{\zeta_*^{n-1}}(\zeta^{n-1},
v^n)\right](y,s)\,\rmd y\,\rmd s.
\end{eqnarray*}
Setting $t=\infty$ in this equation and
comparing with \eqref{zetasn}, we find that
$\zeta^n(0)=\zeta^n(\infty)=0$
if and only if $\zeta^n_*=\zeta^{n-1}_*$.

From (\ref{un}) and (\ref{zetan}) we conclude that
\begin{eqnarray*}
\lefteqn{
v^n(x,t) = \int_{-\infty}^{\infty}
\mathcal{G}^{n-1}(x, t; y,0) v^{{n-1}}_0(y)\,\rmd y
+\int_0^t\int_{-\infty}^{\infty} \mathcal{G}^{n-1}(x, t; y,s) S^{\zeta_*^{n-1}}(\zeta^{n-1},\dot\zeta^{n-1})(y,s)\,\rmd y\,\rmd s} \\ &&
+\int_0^t\int_{-\infty}^{\infty} \mathcal{G}^{n-1}(x, t;
y,s) \left[ Q^{\zeta_*^{n-1}}(\zeta^{n-1}, v^n)+R^{\zeta_*^{n-1}}(\zeta^{n-1},
v^n)\right]_y(y,s)\,\rmd y\,\rmd s
+  \frac{\partial \bar{u}^\zeta}{\partial \zeta}\Big|_{\zeta_*^{n-1}} (\zeta^n(t)-\zeta^n(0))
\end{eqnarray*}
and thus, by Duhamel's Principle,
\[
v^n(t)-\mathcal{L}^{\zeta_*^{n-1}}v^n =Q^{\zeta_*^{n-1}}(\zeta^{n-1}, v^n)_x
+R^{\zeta_*^{n-1}}(\zeta^{n-1}, v^n)_x +
S^{\zeta_*^{n-1}}(\zeta^{n-1},\dot\zeta^{n-1}) +
\frac{\partial \bar{u}^\zeta}{\partial \zeta}\Big|_{\zeta_*^{n-1}}
\dot \zeta^n(t).
\]
Setting $u^n = v^n +\bar{u}^{\zeta_*^{n-1}+\zeta^{n-1}}$,
we then obtain \eqref{unequation}
by a straightforward calculation comparing with \eqref{mixedperteq},
with the claimed initial data
\[
u^n(\cdot, 0)=
\bar{u}^{\zeta_*^{n-1}+\zeta^{n-1}(0)}(\cdot,0)+u^{n-1}_0(\cdot)
= u_0(\cdot) + (\bar{u}^{\zeta_*^{n-1}+\zeta^{n-1}(0)} -\bar{u}^{\zeta_*^{n-1}})(\cdot,0).
\]
Note that the right-hand side of the above equation is equal to $u_0$
if and only if $\zeta^{n-1}(0)=0$: if $\zeta^n\equiv \zeta^{n-1}$ is a fixed point, this is true if and only if $\zeta^n(0)=0$ or, equivalently, $\zeta^n_*=\zeta^{n-1}_*$.
\end{Proof}

\begin{Remark}\label{inirmk}
In \eqref{unequation}--\eqref{dottaun}, the values of $(v^n, \dot \zeta^n)$ at time $T$ depend only on the values for $0\le t\le T$, and not on future times. By \eqref{un}, we have, evidently,
\begin{equation}\label{ineq}
v^n(\cdot, 0)=v^{n-1}_0.
\end{equation}
\end{Remark}


\subsection{Auxiliary convolution and energy estimates}\label{auxests}

In this section, we shall collect various estimates that we will need to construct fixed points of the iteration scheme explained above. We begin by presenting two corollaries of Theorem~\ref{greenbounds}.

\begin{Corollary}\label{eboundscor}
Under the assumptions of Theorem~\ref{nonlin}, and in the notation of Theorem~\ref{greenbounds}, the following holds for $j=1,2$ and $y\leq0$:
\begin{eqnarray*}
|\pi_j(y,s,t)| & \leq &
C\sum_{a_\mathrm{in}^-}
\left(\errfn\left(\frac{y+a_\mathrm{in}^-(t-s)}{\sqrt{4(t-s)}}\right)
-\errfn\left(\frac{y-a_\mathrm{in}^-(t-s)}{\sqrt{4(t-s)}}\right)\right) \nonumber \\ \nonumber
|\pi_j (y,s,t) - \pi_j (y,\infty,s)| & \leq &
C \errfn\left(\frac{|y|-a(t-s)}{M\sqrt{(t-s)}}\right) \qquad \mbox{ for some } a>0 \\ \nonumber
|\partial_t \pi_j(y,s,t)| & \leq &
C (t-s)^{-1/2} \sum_{a_\mathrm{in}^-} \rme^{-|y+a_\mathrm{in}^-(t-s)|^2/M(t-s)} \\
|\partial_y \pi_j(y,s,t)| & \leq &
C (t-s)^{-1/2} \sum_{a_\mathrm{in}^-} \rme^{-|y+a_\mathrm{in}^-(t-s)|^2/M(t-s)} \\ \nonumber &&
+\; C\rme^{-\eta|y|}\sum_{a_\mathrm{in}^-}
\left(\errfn\left(\frac{y+a_\mathrm{in}^-(t-s)}{\sqrt{4(t-s)}}\right)
-\errfn\left(\frac{y-a_\mathrm{in}^-(t-s)}{\sqrt{4(t-s)}}\right)\right) \\ \nonumber
|\partial_y \pi_j (y,s,t) - \partial_y \pi_j(y,\infty,s)| & \leq &
C (t-s)^{-1/2} \sum_{a_\mathrm{in}^-} \rme^{-|y+a_\mathrm{in}^-(t-s)|^2/M(t-s)} \\ \nonumber
|\partial_{yt} \pi_j(y,s,t)| & \leq &
C\left( (t-s)^{-1}+(t-s)^{-1/2}\rme^{-\eta|y|} \right) \sum_{a_\mathrm{in}^-} \rme^{-|y+a_\mathrm{in}^-(t-s)|^2/M(t-s)}.
\end{eqnarray*}
Symmetric estimates are true for $y\ge 0$.
\end{Corollary}

\begin{Proof}
This follows from a straightforward calculation using \eqref{e_1} and \eqref{ljkbounds};
see \cite{MasciaZumbrun03,Zumbrun00} and \cite[Remark~7]{HowardZumbrun06}.
\end{Proof}

In the next corollary, we state a useful property of $\pi=(\pi_1,\pi_2)$.

\begin{Corollary}[\cite{HowardZumbrun06}]
Under the assumptions of Theorem~\ref{nonlin}, and in the notation of Theorem~\ref{greenbounds},
\begin{equation}\label{Ifacteq}
\int^\infty_{-\infty} \pi(y,s,\infty) \frac{\partial\bar{u}^{\zeta}(\zeta_*)}{\partial \zeta}(y,s)\,\rmd y = \mathrm{id}_{\mathbb{R}^2}.
\end{equation}
\end{Corollary}

\begin{Proof}
This follows from the fact that both $\bar{u}^{\zeta}_x$ and $\bar{u}^{\zeta}_t$ are, for any fixed $\zeta$, stationary solutions of the linearized equations. Hence
\[
\int_{-\infty}^{\infty}\mathcal{G}(x,t;y,s)
\bar{u}_x^\zeta(y,s)\,\rmd y\equiv
\bar{u}_x^\zeta(x,t), \qquad \int_{-\infty}^{\infty}\mathcal{G}(x,t;y,s)
\bar{u}_t^\zeta(y,s)\,\rmd y\equiv
\bar{u}_t^\zeta(x,t).
\]
Because $\bar{u}_x$ and $\bar{u}_t$ are linearly independent and, under the assumption of spectral stability, $\mathcal{E}_1$ and $\mathcal{E}_2$ represent the only nondecaying parts of $\mathcal{G}(x,t;y,s)$, we have
\[
\lim_{t\to \infty} \int_{-\infty}^{\infty} \pi_1(y,s,t)\bar{u}_t^\zeta(y,s)dy = \lim_{t\to \infty} \int_{-\infty}^{\infty} \pi_2(y,s,t)\bar{u}_x^\zeta(y,s)dy = 0, 
\]
which leads to \eqref{Ifacteq}. 
\end{Proof}

Next, we investigate the dependence of $\pi$ on $\zeta_*$.

\begin{Proposition}[Parameter-dependent bounds]\label{parambds}
Under the assumptions of Theorem~\ref{nonlin}, and in the notation of Theorem~\ref{greenbounds}, there exists a constant $C$ such that
\begin{eqnarray}
|\partial_{\zeta_*}\pi| & \sim & C |\pi| \nonumber \\ \nonumber
|\partial_{\zeta_*}\partial_t\pi| & \sim & C |\partial_t\pi| \\ \nonumber
|\partial_{\zeta_*}(\pi(y,s,t)-\pi(y,s,\infty))| & \sim & C |(\pi(y,s,t)-\pi(y,s,\infty))|
\\ \label{ederivbds}
|\partial_{\zeta_*}\partial_y\pi| & \sim & C |\partial_y\pi| \\ \nonumber
|\partial_{\zeta_*}(\pi_y(y,s,t)-\pi_y(y,s,\infty))| & \sim & C|(\pi_y(y,s,t)-\pi_y(y,s,\infty)| \\ \nonumber
|\partial_{\zeta_*}\partial_t\partial_y\pi| & \sim & C |\partial_t\partial_y\pi|
\end{eqnarray}
and
\begin{equation}\label{Gderivbds}
|\partial_{\zeta_*}\partial_{y}^\alpha \tilde{\mathcal{G}}(x,t;y,s)| \sim
C|\partial_{y}^\alpha \tilde{\mathcal{G}}(x,t;y,s)|
\end{equation}
for $0\leq|\alpha|\leq1$ and $y\leq0$, and symmetrically for $y\geq0$, where by $\sim$ we mean that the left-hand side obeys the same bounds as given for the right-hand side in Theorem~\ref{greenbounds} and Corollary~\ref{eboundscor} above.
\end{Proposition}

\begin{Proof}
Evidently, $\partial_{\zeta_*}=(\partial_x+\partial_y, \partial_t+\partial_s)$. The bounds \eqref{ederivbds} for $\pi_1$ and $\pi_2$ follow by direct calculation, together with the estimate
\[
\partial_{\zeta_*} \partial_y l_{j,\mathrm{in}}^\pm(y,s) = \mathrm{O}(\rme^{-\eta|y|}).
\]
To obtain the estimates \eqref{Gderivbds}, we first indicate how to estimate $s$ and $t$ derivatives of $\tilde{\mathcal{G}}$, which we have not discussed so far. From the decomposition $\tilde{\mathcal{G}}=\sum_{j=0}^\ell\mathcal{G}_j+\mathcal{G}_*-\mathcal{E}$
together with the bounds on $\mathcal{G}_j$ stated in Lemma~\ref{parametrixlem}, it is sufficient to estimate the $t$- and $s$-derivatives of
\begin{equation}\label{srep}
\mathcal{G}_*(x,t;y,s) - \mathcal{E} =
\frac{1}{2\pi\rmi} \oint_{\Gamma} \rme^{\sigma(t-s)} [R(x,y,\sigma;s)](t) \,\rmd\sigma,
\end{equation}
where we have written the integral in terms of the unshifted time-coordinates, with initial time $t=s$, to capture the dependence on $s$---everywhere else in the paper, we shifted to $t=0$, yielding the factor $\rme^{\sigma t}$ in place of $\rme^{\sigma(t-s)}$. In \eqref{srep}, we see that $t$- and $s$-derivatives yield a factor $\pm \sigma$, where they fall on $\rme^{\sigma(t-s)}$, and otherwise contribute a $t$- or $s$-derivative of $R$. The additional factor $\sigma$ yields a contribution to the inverse Laplace transform that is smaller by a factor $(1+(t-s))^{-\frac12}$ than our estimate for $|\mathcal{G}_*-\mathcal{E}|$. Likewise, derivatives that fall on $R$ are harmless, since we established pointwise for those in Theorem~\ref{resbounds}. Similarly, using again the results from \S\ref{S:res_decomp_est}, the effect of additional spatial derivatives on \eqref{srep} is a factor of order $|\sigma|+\rme^{-\eta|y|}$ for $y$-derivatives, while $x$-derivatives do not change the bounds we obtain. Thus, for either temporal or spatial derivatives, the issue reduces to the estimates on $\mathcal{G}_j$ already carried out.
\end{Proof}

We shall also make use of the following technical lemmas proved in \cite{HowardZumbrun06}.

\begin{Lemma}[Linear estimates I]\label{iniconvolutions}
Under the assumptions of Theorem~\ref{nonlin}, there exists a constant $C$ such that
\begin{eqnarray*}
\int_{-\infty}^{\infty}|\tilde{\mathcal{G}}(x,t;y,0)|(1+|y|)^{-3/2}\,\rmd y & \le &
C(\theta_\mathrm{gauss}+\theta_\mathrm{inner}+\theta_\mathrm{outer})(x,t) \\
\int_{-\infty}^{\infty}|\pi_t(y,0,t)|(1+|y|)^{-3/2}\,\rmd y & \leq & C(1+t)^{-3/2} \\
\int_{-\infty}^{\infty}|\pi(y,0,t)|(1+|y|)^{-3/2}\,\rmd y & \leq & C \\
\int^{\infty}_{-\infty} |\pi(y,0,t)-\pi(y,0,\infty)| (1+|y|)^{-3/2}\,\rmd y & \leq &
C(1+t)^{-1/2}
\end{eqnarray*}
for $0\le t<\infty$, where $\tilde{\mathcal{G}}$ and $\pi=(\pi_1,\pi_2)^T$ are defined in Theorem~\ref{greenbounds}.
\end{Lemma}

\begin{Proof}
These estimates can be established exactly as in \cite[Lemma~3]{HowardZumbrun06}, using Corollary~\ref{eboundscor} and the explicit bounds for $\tilde{\mathcal{G}}$.
\end{Proof}

\begin{Lemma}[Nonlinear estimates I]\label{convolutions}
Let
\[
\Theta(y,s) := (1+s)^{1/2}s^{-1/2}(\theta_\mathrm{gauss} + \theta_\mathrm{inner}+\theta_\mathrm{outer})^2(y,s) + (1+s)^{-1} (\theta_\mathrm{gauss}+\theta_\mathrm{inner}+\theta_\mathrm{outer})(y,s).
\]
Under the assumptions of Theorem~\ref{nonlin}, there is a constant $C$ such that
\begin{eqnarray*}
\int_0^t\int_{-\infty}^{\infty}|\tilde{\mathcal{G}}_y(x,t;y,s)|\Theta(y,s)\,\rmd y \,\rmd s & \leq &
C(\theta_\mathrm{gauss}+\theta_\mathrm{inner}+\theta_\mathrm{outer})(x,t) \\
\int_0^t\int_{-\infty}^{\infty}|\pi_{yt}(y,s,t)|\Theta(y,s)\,\rmd y \,\rmd s & \leq & C(1+t)^{-1} \\
\int_t^{\infty} \int_{-\infty}^{\infty}|\pi_y(y,s,\infty)|\Theta(y,s)\,\rmd y & \leq & C(1+t)^{-1/2} \\
\int_0^t\int_{-\infty}^{\infty} |\pi_y(y,s,t)- \pi_y(y,s,\infty)| \Theta(y,s)\,\rmd y\,\rmd s & \leq & C(1+t)^{-1/2}
\end{eqnarray*}
for $0\le t<\infty$, where $\tilde{\mathcal{G}}$ and $\pi=(\pi_1,\pi_2)^T$ are defined as in Theorem~\ref{greenbounds}.
\end{Lemma}

\begin{Proof}
The estimates can be established as in \cite[Lemma~4]{HowardZumbrun06}, again using Corollary~\ref{eboundscor} and the explicit bounds for $\tilde{\mathcal{G}}$.
\end{Proof}

\begin{Lemma}[Nonlinear estimates II]\label{occonvolutions}
Let
\begin{eqnarray*}
\Phi_1(y,s) & := & \rme^{-\eta|y|}(1+s)^{-1/2}(\theta_\mathrm{gauss}+\theta_\mathrm{inner}+\theta_\mathrm{outer})(y,s)
\;\leq\; C\rme^{-\eta|y|/2}(1+s)^{-3/2} \\
\Phi_2(y,s) & := & \rme^{-\eta|y|}(1+s)^{-3/2}.
\end{eqnarray*}
Under the assumptions of Theorem~\ref{nonlin}, there is a constant $C$ such that
\begin{eqnarray*}
\int_0^t\int_{-\infty}^{\infty} |\tilde{\mathcal{G}}_y(x,t;y,s)| \Phi_1(y,s) \,\rmd y \,\rmd s
& \leq & C(\theta_\mathrm{gauss}+\theta_\mathrm{inner}+\theta_\mathrm{outer})(x,t) \\
\int_0^t\int_{-\infty}^{\infty} |\pi_{yt}(y,s,t)| \Phi_1(y,s) \,\rmd y \,\rmd s
& \leq & C(1+t)^{-1} \\
\int_0^{t}\int_{-\infty}^{\infty} |\pi_{y}(y,s,t)-\pi_y(y,s,\infty)| \Phi_1(y,s) \,\rmd y \,\rmd s  
& \leq & C(1+t)^{-1/2} \\
\int_t^{\infty}\int_{-\infty}^{\infty} |\pi_{y}(y,s,\infty)| \Phi_1(y,s) \,\rmd y \,\rmd s
& \leq & C(1+t)^{-1/2}
\end{eqnarray*}
and
\begin{eqnarray*}
\int_0^t\int_{-\infty}^{\infty} |\tilde G(x,t;y,s)| \Phi_2(y,s) \,\rmd y \,\rmd s
& \leq & C(\theta_\mathrm{gauss}+\theta_\mathrm{inner}+\theta_\mathrm{outer})(x,t) \\
\int_0^t\int_{-\infty}^{\infty} |\pi_{t}(y,s,t)| \Phi_2(y,s) \,\rmd y \,\rmd s
& \leq & C(1+t)^{-3/2} \\
\int_0^t\int_{-\infty}^{\infty} |\pi(y,s,t)-\pi(y,s,\infty)| \Phi_2(y,s) \,\rmd y \,\rmd s
& \leq & C(1+t)^{-3/2},\\
\int_t^{\infty}\int_{-\infty}^{\infty} |\pi(y,s,\infty)| \Phi_2(y,s) \,\rmd y \,\rmd s
& \leq & C(1+t)^{-1/2}
\end{eqnarray*}
for $0\le t<\infty$, where $\tilde G$ and $\pi$ are defined as in Theorem~\ref{greenbounds}.
\end{Lemma}

\begin{Proof}
This follow as in \cite[Lemma~5]{HowardZumbrun06} and \cite[Lemmas~4.2-4.3]{RaoofiZumbrun07}.
\end{Proof}

We shall also use the following auxiliary energy estimate, adapted essentially unchanged from \cite{MasciaZumbrun04_largeamp, Zumbrun04, Raoofi05, RaoofiZumbrun07}. Let $u$ be a solution of 
\begin{equation}\label{viscous}
u_t + f(u)_x - u_{xx} = \frac{\partial\bar{u}}{\partial\zeta}(x)\Big|_{\zeta_*} \gamma(t)
\end{equation}
for some given function $\gamma(t)$, and define the function $v$ via
\begin{equation}\label{pert}
u = \bar{u}^{\zeta_*+\zeta(t)} + v.
\end{equation}

\begin{Lemma}[\cite{MasciaZumbrun04_largeamp, Zumbrun04, Raoofi05}] \label{aux}
Under the hypotheses of Theorem~\ref{nonlin}, let $u_0\in H^3$, and suppose that, for $0\le t\le T$, the suprema of $ |\dot\zeta|$ and $|\gamma|$ and the $H^3$ norm of $v$, defined by \eqref{viscous} and \eqref{pert}, each remain bounded by a sufficiently small constant. There are then constants $\theta_{1,2}>0$ so that
\begin{equation}\label{Ebounds}
\|v(t)\|_{H^3}^2 \leq C \rme^{-\theta_1 t} \|v(0)\|^2_{H^3} + C \int_0^t \rme^{-\theta_2(t-\tau)} (|v|_{L^2}^2 + |\dot{\zeta}|^2 + |\gamma|^2)(\tau)\,\rmd\tau
\end{equation}
for all $0\leq t\leq T$.
\end{Lemma}

\begin{Proof}
This follows by parabolic energy estimates similar to, but in fact much simpler, than those used to treat the case of partially elliptic viscosity in \cite{MasciaZumbrun04_largeamp, Zumbrun04, Raoofi05}. Specifically, we write the perturbation equation for $v$ as
\begin{equation}\label{vperturteq}
v_t + \left[ \int_0^1 f_u(\bar u^{\zeta_*+\zeta(t)}+\tau v)\,\rmd\tau\, v\right]_x - v_{xx} =
\frac{\partial\bar{u}(x)}{\partial\zeta}\Big|_{\zeta_*} (\gamma(t)+\dot{\zeta}(t)).
\end{equation}
Observing that $\partial_x^\alpha(\partial\bar{u}/\partial\zeta)|_{\zeta_*}(x)=\mathrm{O}(\rme^{-\eta|x|})$ is bounded in the $L^1$ norm for $|\alpha|\leq3$, we take the $L^2$ inner product in $x$ of $\sum_{j=0}^3\partial_x^{2j}v$ against \eqref{vperturteq}, integrate by parts and rearrange the resulting terms to arrive at the inequality
\[
\partial_t \|v\|_{H^3}^2(t) \leq -\theta \|\partial_x^4 v\|_{L^2}^2 +
C \left(\|v\|_{H^3}^2 + |\gamma(t)|^2 + |\dot{\zeta}(t)|^2\right),
\]
where $\theta>0$, for some sufficiently large $C>0$, so long as $\|v\|_{H^3}$ remains bounded. Using the Sobolev interpolation
\[
\|v\|_{H^3}^2 \leq \tilde{C}^{-1} \|\partial_x^4 v\|_{L^2}^2 + \tilde{C} \| v\|_{L^2}^2
\]
for $\tilde{C}>0$ sufficiently large, we obtain 
\[
\partial_t \|v\|_{H^3}^2(t) \leq -\tilde{\theta} \|v\|_{H^3}^2 + C \left(\|v\|_{L^2}^2 + |\gamma(t)|^2 + |\dot{\zeta}(t)|^2\right),
\]
from which \eqref{Ebounds} follows by Gronwall's inequality.
\end{Proof}

\begin{Lemma}[\cite{RaoofiZumbrun07}] \label{weightest}
Under the hypotheses of Theorem~\ref{nonlin}, let $M_0:=\|(1+|x|^2)^{3/4}v_0(x)\|_{H^3}<\infty$, and suppose that, for $0\le t\le T$, the suprema of $ |\dot\zeta|$ and $|\gamma|$
and the $H^{3}$ norm of $v$, determined by \eqref{viscous} and \eqref{pert}, each remain bounded by some constant $C>0$. Then there exists some $M=M(C)>0$ such that, for all $0\le t\le T$, 
\begin{equation}\label{weightedEbounds}
\left\|(1+|x|^2)^{3/4}v(x,t)\right\|_{H^3}^2 \le M\rme^{Mt} \left( M_0 +\int_0^t (|\dot\zeta|^2 + |\gamma|^2)(\tau)\, \rmd\tau \right).
\end{equation}
\end{Lemma}

\begin{Proof}
This follows by standard Friedrichs symmetrizer estimates carried out in the weighted $H^3$ norm. Specifically, making the coordinate change $v=(1+|x|^2){-3/4}w$, we obtain from \eqref{vperturteq} the modified equation
\begin{equation}\label{zperturteq}
w_t + \left[ \int_0^1 f_u(\bar u^{\zeta_*+\zeta(t)}+\tau v)\,\rmd\tau\, w\right]_x - w_{xx} =
\frac{\partial\bar{u}(x)}{\partial\zeta}|_{\zeta_*} (\gamma(t) + \dot{\zeta}(t))
\end{equation}
plus lower-order commutator terms of order $\mathrm{O}(|w|+|w||v|+|w_x||v|)$, which are bounded by $M(|w|+|w_x|)$ by assumption, and similarly in the equations for $\partial_x^j w$ for $j=1,\dots,3$. Likewise, $\int f_u\,\rmd\tau$ and its derivatives up to order two remain uniformly bounded by Sobolev embedding and the assumed bound on $\|v\|_{H^3}$. Performing the same energy estimates on \eqref{zperturteq} as carried out on \eqref{vperturteq} in the proof of Lemma~\ref{aux}, we readily obtain the result by Gronwall's inequality (indeed, somewhat better, thanks to parabolic smoothing). We refer to \cite[Lemma~5.2]{RaoofiZumbrun07} for further details in the general partially parabolic case.
\end{Proof}

\begin{Remark}\label{zetareg}
Using Sobolev embeddings and equation \eqref{viscous}, we see that Lemma~\ref{weightest} immediately implies that, if $\|(1+|x|^2)^{3/4}v_0(x)\|_{H^3}<\infty$ and if $\|v(\cdot,t)\|_{H^3}$, $|\dot\zeta(t)|$ and $|\gamma(t)|$ are uniformly bounded on $0\le t\le T$, then
\[
|(1+|x|^2)^{3/4}v(x,t)| \quad\mbox{and}\quad
|(1+|x|^2)^{3/4}v_t(x,t)|
\]
are uniformly bounded on $0\le t\le T$ as well.
\end{Remark}


\subsection{Proof of Theorem~\ref{nonlin}}\label{proof}

We are almost ready to prove the main theorem. We will need the following definitions and lemmas. Define the norm
\begin{equation}\label{Beta1norm}
|\zeta|_{\mathcal{B}_1} := |\zeta(t)(1+t)^{1/2}|_{L^\infty_t} + |\dot \zeta(t)(1+t)|_{L^\infty_t}
\end{equation}
and the Banach space $\mathcal{B}_1:=\{\zeta:\mathbb{R}\to\mathbb{R}^2:\; |\zeta|_{\mathcal{B}_1}<\infty\}$. We also define
\begin{equation}
|v|_{\mathcal{B}_2} := \left| \frac{1}{(\theta_\mathrm{gauss}+\theta_\mathrm{inner}+\theta_\mathrm{outer})(x,t)} v(x,t)\right|_{L^{\infty}_{x,t}}
\label{Beta2norm}
\end{equation}
and the Banach space $\mathcal{B}_2:=\{v:\mathbb{R}^2\to\mathbb{R}^n:\; |v|_{\mathcal{B}_2}<\infty\}$. The next lemma gives local existence of the integral equations \eqref{unequation}--\eqref{dottaun} for $(v^n,\dot{\zeta}^n)$.

\begin{Lemma}[$H^3$ local theory] \label{local}
Under the hypotheses of Theorem~\ref{nonlin},
let
\[
M_1:=\|v_0(x)\|_{H^3} + |\zeta^{n-1}|_{\mathcal{B}_1} + |\zeta_*^{n-1}|
<\infty,
\]
where $|\cdot|_{\mathcal{B}_1}$ is defined in (\ref{Beta1norm}).
There exists some $T=T(M_1)>0$ sufficiently small
and $C=C(M_1,T)>0$ sufficiently large such that,
on $0\le t\le T$, there exists a unique solution
\[
(v^n,\dot \zeta^n)\in L^\infty_t(H^3_x) \times C^{0}_t
\]
of \eqref{unequation}--\eqref{dottaun} that satisfies
\[
\|v^n(t)\|_{H^3},\, |\dot \zeta^n(t)|\le CM_1.
\]
\end{Lemma}

\begin{Proof}
Short-time existence, uniqueness, and stability follow by (unweighted) energy estimates in $v^n$ similar to \eqref{weightedEbounds} combined with more straightforward estimates on $\dot \zeta^n$ carried out directly from integral equations \eqref{dotdeltan} and \eqref{dottaun}, using a standard (bounded high norm, contractive low norm) contraction mapping argument like those described in \cite[\S4.2.4]{Zumbrun04} and \cite[Proof of Proposition~1.6 and Exercise~1.9]{Zumbrun04_lecturenotes}. We omit the details.
\end{Proof}

\begin{Remark}
A crucial point is that equations \eqref{unequation}--\eqref{dottaun}
depend only on values of $(v^n, \dot\zeta^n)$ on the range $t\in [0,T]$;
see Remark~\ref{inirmk}.
\end{Remark}

The next result allows us to extend solutions provided they stay bounded in an appropriate sense.

\begin{Lemma}\label{continuation}
Under the hypotheses of Theorem~\ref{nonlin}, assume that
\[
(v^n, \dot\zeta^n)\in L^\infty_t(H^3_x)\times L^{\infty}_t
\]
satisfy \eqref{unequation}--\eqref{dottaun} on $0\le t\le T$,
and define 
\begin{equation}\label{zeta}
\beta(t) :=\sup_{x\in\mathbb{R},\, s\in[0,t]}
\left(\frac{|v^n(x,s)|}{(\theta_\mathrm{gauss}+\theta_\mathrm{inner}+\theta_\mathrm{outer})(x,s)}
+|\dot\zeta^n(s)|(1+s) \right).
\end{equation}
If $\beta(T)$, $\|v^{n-1}_0\|_{H^3}$, and
$|\zeta^{n-1}|_{\mathcal{B}_1}$
are bounded by some sufficiently small $\beta_0>0$, then,
for some $\epsilon>0$, the solution $(v^n, \dot \zeta^n)$, and thus $\beta$,
extends to $0\le t\le T+\epsilon$, and $\beta$ is bounded and continuous on $0\le t\le T+\epsilon$.
\end{Lemma}

\begin{Proof}
By \eqref{ineq}, we have
\[
\|v^{n-1}_0\|_{H^3}= \|v^{n}(\cdot, 0)\|_{H^3},
\]
and Lemma~\ref{aux} and smallness of both $\beta(T)$ and $|\zeta^{n-1}|_{\mathcal{B}_1}$ therefore imply boundedness (and smallness) of $\|v^n(t)\|_{H^3}$ and $|\dot \zeta^n(t)|$ on $0\le t\le T$. By Lemma~\ref{local}, applied with initial data $(v^n(T),\zeta^n(T))$, this implies existence and boundedness of $v^n$ in $H^3$ and $\dot{\zeta}^n$ in $\mathbb{R}^2$ on $0\le t\le T+\epsilon$ for some $\epsilon>0$, and thus, by Remark~\ref{zetareg}, boundedness and continuity of $\beta$ on $0\le t\le T+\epsilon$.
\end{Proof}

\begin{Lemma}\label{uexistence}
For $M>0$ sufficiently large, $M_0:=\|(1+|\cdot|^2)^{3/4}(u_0(\cdot)- \bar{u}(\cdot, 0))\|_{H^3}$
sufficiently small, and $|\zeta^{n-1}|_{\mathcal{B}_1}+ M|\zeta_*^{n-1}|\le 2CM_0$ for some constant $C>0$, there exists a $C_1>0$ sufficiently large so that solutions $(v^n, \zeta^n, \zeta_*^n)$ of \eqref{un}--\eqref{zetasn} exist for all $t\ge 0$ and satisfy 
\begin{equation}\label{H4unests}
\|v^n\|_{H^3} \le C_1M_0
\end{equation}
and
\begin{equation}\label{unest}
|v^n|_{\mathcal{B}_2}+|\zeta^{n}|_{\mathcal{B}_1}+ M|\zeta_*^n|\le 2CM_0.
\end{equation}
\end{Lemma}

\begin{Proof}
Define $\beta$ as in \eqref{zeta}.
We claim that, if we can show that 
\begin{equation}\label{claim}
\beta(t)\le CM_0 + C_*(M_0+\beta(t))^2 
\end{equation}
for some fixed $C$ and $C_*>0$, for all time $t$ such that the
solution $(v^n, \dot\zeta^n)\in L^\infty_t(H^3_x)\times
C^{0}_t$ of \eqref{unequation}--\eqref{dottaun} exists and
\begin{equation}\label{finalbd} 
\beta(t)\le \frac32 CM_0, 
\end{equation}
then we can conclude that the solution $(v^n, \dot \zeta^n)$ in fact exists
and satisfies \eqref{finalbd} for all $t\ge 0$, provided
\[
M_0< \frac{2}{5C_*(3C+4)}
\]
is sufficiently small. To see this, first note that, by \eqref{ineq}, \eqref{inidat},
and the exponential decay that follows from Hypothesis~\ref{H1}, we have
\begin{eqnarray}
\|v^{n}(\cdot, 0)\|_{H^3}
& = & \|v^{n-1}_0\|_{H^3} \nonumber \\ \nonumber 
& = & \|u_0-\bar{u}^{\zeta_*^{n-1}}\|_{H^3} \\ \nonumber
& \leq & \|u_0-\bar{u}\|_{H^3} + \| \bar{u}^{\zeta_*^{n-1}}- \bar{u}\|_{H^3} \\ \nonumber
& \leq & M_0 + c_1|\zeta_*^{n-1}| \\ \label{inismall}
& \leq & M_0 \left(1 + \frac{2c_1C}{M}\right),
\end{eqnarray}
which is small when $M_0$ is sufficiently small.
Defining $T$ to be the maximum time up to which a solution $(v^n, \dot \zeta^n)$
exists and $\beta(t) \le \beta_0$ is sufficiently small, 
we can use the assumed bounds on $\zeta^{n-1}$ and $\dot \zeta^{n-1}$ to apply Lemma~\ref{continuation} and conclude that $(v^n, \dot \zeta^n)$ exists up to time $T+\epsilon$ for some $\epsilon>0$
and that $\beta$ remains bounded and continuous up to $T+\epsilon$ as well.
Observing that \eqref{claim} together with
$M_0< 2/(5C_*(3C+4))$ implies that $\beta(t)<\frac32 CM_0$ whenever $\beta(t) \le\frac32 CM_0$,
we find by continuity that $\beta(t)\le\frac32 CM_0$ up to $t=T+\epsilon$
as claimed. We may now repeat this process infinitely many times to conclude that the solution exists for all $t \geq 0$ and satisfies (\ref{finalbd}).

We now show how (\ref{finalbd}) can be used to prove the lemma.
By the definition of $\beta$, (\ref{finalbd}) implies
\[
|v^n|_{\mathcal{B}_2}+ |\dot\zeta^n(t)(1+t)|_{L^\infty_t}\le \frac32 CM_0.
\]
This proves part of equation (\ref{unest}). To complete it, we must show that
\begin{equation}\label{penbd}
|\zeta^n(t)|\le \frac{CM_0}{4} (1+t)^{-1/2}
\end{equation}
and
\begin{equation}\label{lastbd}
|\zeta_*^n|\le \frac{CM_0}{4M}.
\end{equation}
Establishing \eqref{penbd} and \eqref{lastbd}
also proves that
the integral equations for $\zeta^n$ and $\zeta_*^n$ converge, and so we would also obtain,
by Lemma~\ref{soln} and the fact
that $(v^n,\dot \zeta^n)$ satisfies \eqref{unequation}--\eqref{dottaun}
for all $t\ge 0$, that $(v^n, \zeta^n, \zeta_*^n)$
satisfies \eqref{un}--\eqref{zetasn} as claimed.
Finally, recalling \eqref{unequation} and applying Lemma~\ref{aux}
with $\gamma:=\dot \zeta^{n}- \dot \zeta^{n-1}$, we obtain
\eqref{H4unests} so long as \eqref{claim} remains valid,
controlling $\|v^n\|_{H^3}$ by integrating the right-hand side of
\eqref{Ebounds} and using \eqref{finalbd}, the definition of
$\beta$, and the assumed bounds on $\dot \zeta^{n-1}$. We shall carry
out this last calculation in detail in equation (\ref{2calc}), in
the course of proving \eqref{claim}. Thus, it remains to prove \eqref{claim}, \eqref{penbd}, and \eqref{lastbd}.

We now establish \eqref{claim} using \eqref{finalbd}. By Lemma~\ref{aux}, 
we have
\begin{eqnarray}
\|v^n(t)\|_{H^3}^2 & \leq & c \|v^n(0)\|^2_{H^3}\rme^{-\theta t}
+ c\int_0^t \rme^{-\theta_2 (t-\tau )}(|v^n|_{L^2}^2+ |\dot\zeta^n|^2
+|\dot \zeta^n- \dot \zeta^{n-1}|^2)(\tau)\, d\tau \nonumber \\ \nonumber
& \leq & c \|v^n(0)\|^2_{H^3}\rme^{-\theta t}
+ c\int_0^t \rme^{-\theta_2 (t-\tau )}(|v^n|_{L^2}^2+ \max\{|\dot\zeta^n|^2, | \dot \zeta^{n-1}|^2 \})(\tau)\, d\tau \\ \nonumber
& \leq & c_2\big(\|v^n(0)\|^2_{H^3}+ \beta(t)^2\big) (1+t)^{-1/2} \\ \nonumber
& \leq & c_2\big(M_0^2(1 + 2c_1C/M)^2 + (3CM_0/2)^2\big) (1+t)^{-1/2} \\ \label{2calc}
& \leq & (C_1M_0)^2 (1+t)^{-1/2},
\end{eqnarray}
for $C_1>0$ sufficiently large and $M_0$ sufficiently small,
by \eqref{inismall}, \eqref{finalbd}, and the definition of $\beta$.
With \eqref{mixedperteq}, \eqref{H4unests}, the assumption that $|\zeta|_{\mathcal{B}_1}\le 2CM_0$,
and the definitions of $\beta$ and $|\cdot|_{\mathcal{B}_1}$,
we obtain readily
\begin{equation}\label{sourcebds}
|Q^{\zeta_*} +R^{\zeta_*}| \le c(\beta^2+4C^2M_0^2)(\Theta+ \Phi_1)
\end{equation}
and
\begin{equation}\label{derivsourcebds}
|S^{\zeta_*}| \le c(\beta^2+4C^2M_0^2) \Phi_2,
\end{equation}
where $\Theta$ and $\Psi_{1,2}$ are defined in Lemmas~\ref{convolutions}--\ref{occonvolutions}. Applying Lemmas~\ref{iniconvolutions}--\ref{occonvolutions} to \eqref{un}, \eqref{dotdeltan}, and \eqref{dottaun}, we thus obtain \eqref{claim} as claimed. Likewise, we obtain \eqref{penbd} from \eqref{zetan} and \eqref{inismall}, using Lemmas~\ref{iniconvolutions}--\ref{occonvolutions} and the definitions of $\beta$ and $|\cdot|_{\mathcal{B}_1}$.

It remains only to establish \eqref{lastbd}. This is more delicate due to the appearance of $M$ in the denominator of the right-hand side and depends on the key fact that the estimate $\zeta_*^n$ of the asymptotic shock location is to linear order insensitive to the initial guess $\zeta^{n-1}$. To see this, decompose the expression \eqref{zetasn} for $\zeta_*^n$ into its linear part
\begin{eqnarray}\label{deslin}
I & := & \zeta_*^{n-1}-\int_{-\infty}^{\infty}\pi^{n-1}(y,0,\infty) v^{{n-1}}_0(y)\,\rmd y
\\ \nonumber & = &
\underbrace{\zeta_*^{n-1}-\int_{-\infty}^{\infty}\pi|_{\zeta_*=0}(y,0,\infty)
(\bar{u} -\bar{u}^{\zeta_*^{n-1}})(y)\,\rmd y}_{=:I_a}
- \underbrace{\int_{-\infty}^{\infty}\pi|_{\zeta_*=0}(y,0,\infty) (u_0 -\bar{u})(y)\,\rmd y}_{=:I_b}
\\ \nonumber &&
- \underbrace{\int_{-\infty}^{\infty}(\pi^{n-1}-\pi|_{\zeta_*=0})(y,0,\infty) v^{n-1}_0(y)\,\rmd y}_{=:I_c}
\end{eqnarray}
and its nonlinear part
\begin{eqnarray}\label{desnonlin}
II & := & -\int_0^{\infty}\int_{-\infty}^{\infty} \pi^{n-1}(y,0, \infty) S^{\zeta_*^{n-1}}(\zeta^{n-1},\dot\zeta^{n-1})(y,s)\,\rmd y\,\rmd s \\ \nonumber &&
-\int_0^{\infty}\int_{-\infty}^{\infty}\pi^{n-1}(y,0,\infty)
\left[Q^{\zeta_*^{n-1}}(u^n, u^n_x)+R^{\zeta_*^{n-1}}(\zeta^{n-1},
u^n)\right]_y(y,s)\,\rmd y\,\rmd s.
\end{eqnarray}
By estimates like the previous ones, we readily obtain
\[
|II|\le 2c(2CM_0)^2,
\]
which is $\ll CM_0/16M$ for $M_0$ sufficiently small.
Likewise, $|I_c|\le c|\zeta_*|M_0$, by \eqref{inismall},
\eqref{ederivbds}, and the Mean Value Theorem,
hence is $\ll CM_0/16M$ for $M_0$ sufficiently small
(recall that we assume $|\zeta_*|\le 2CM_0$),
and
\[
|I_b| \leq c\| u_0-\bar{u}\|_{L^1}
\leq c_2 \|(1+|x|^2)^{3/4}( u_0-\bar{u})\|_{H^3}
\leq c_2M_0,
\]
hence is $\ll CM_0/16M$ for $M>0$ sufficiently large. Finally, Taylor expanding, and recalling~\ref{H1} and
\eqref{Ifacteq}, we obtain
\[
I_a =
\zeta_*^{n-1}-\zeta_*^{n-1}\int_{-\infty}^{\infty}\pi|_{\zeta_*=0}(y,0,\infty)
\frac{\partial\bar{u}^{\zeta_*}(y)}{\partial \zeta_*}\Big|_{\zeta_*=0}\,\rmd y + \mathrm{O}(|\zeta_*|^2)
= \mathrm{O}(|\zeta_*|^2),
\]
which is also $\ll CM_0/16M$ for $M_0$ sufficiently small
(recall that we assume $|\zeta_*|\le 2CM_0$).
Summing these terms up, we obtain \eqref{lastbd} for $M_0$ sufficiently small
and $M>0$ sufficiently large, as claimed.
This completes the proof.
\end{Proof}

\begin{Proof}[ of Theorem~\ref{nonlin}.]
We define
\[
|(v,\zeta)|_*:= |v|_{\mathcal{B}_1}+ M|\zeta|, \qquad
(v,\zeta)\in \mathcal{B}_1\times \mathbb{R}.
\]
Lemma~\ref{uexistence} implies that, for $r>0$ small, $M>0$ sufficiently large, and
\[
M_0=\|(1+|\cdot|^2)^{3/4}(u_0(\cdot)-\bar{u}(\cdot,0))\|_{H^3}
\]
sufficiently small, the mapping $\mathcal{T}=\mathcal{T}(\zeta,\zeta_*)$ is well defined from
\[
B(0,r)\subset \mathcal{B}_1\times \mathbb{R}\to
\mathcal{B}_1\times \mathbb{R}.
\]
To establish the theorem, 
therefore, it suffices to establish that
$\mathcal{T}$ is a contraction on $B(0,r)$ in the norm $|\cdot|_*$.
Indeed, we can then apply Banach's fixed-point theorem to find that
\[
\mathcal{T}(\zeta, \zeta_*)= (\zeta, \zeta_*)
\]
has a unique solution, and Lemma~\ref{soln} shows that the associated function $u$ satisfies \eqref{main} with initial data $u_0$. The stated decay estimates follow from \eqref{H4unests}, \eqref{unest}, and the definition of the norms in the spaces $\mathcal{B}_1$ and $\mathcal{B}_2$ in (\ref{Beta1norm}) and (\ref{Beta2norm}).

To show that $\mathcal{T}$ is a contraction, we need to
establish the Lipschitz bound
\begin{equation}\label{lipbds}
|\mathcal{T}(\zeta, \zeta_*)-\mathcal{T}(\hat\zeta, \hat\zeta_*)|_* \le \alpha |(\zeta,
\zeta_*)-(\hat\zeta, \hat\zeta_*)|_*
\end{equation}
on $B(0,r)$ for some $\alpha<1$ and some sufficiently small $r$. Assume that $(v^n,\zeta^n, \zeta_*^n)$ satisfies \eqref{un}--\eqref{zetasn}
associated with $(\zeta^{n-1},\zeta_*^{n-1})$, while
$(\hat v^n,\hat \zeta^n, \hat \zeta_*^n)$ satisfies \eqref{un}--\eqref{zetasn}
with $(\zeta^{n-1},\zeta_*^{n-1})$ replaced by $(\hat \zeta^{n-1},\hat \zeta_*^{n-1})$. For each $n$, we define the variations
\[
\Delta v^n:=\hat v^n-v^n,
\qquad
\Delta \zeta^n:=\hat \zeta^n-\zeta^n,
\qquad
\Delta \zeta_*^n:=\hat \zeta^n_*-\zeta^n_*,
\]
and, likewise, define $\Delta \tilde{\mathcal{G}}^{n-1}$, $\Delta \pi^{n-1}$, $\Delta S^{n-1}$, $\Delta Q^{n-1}$, and $\Delta R^{n-1}$ in the obvious way. Using equations \eqref{un}--\eqref{zetasn} and the fact that $\hat{f}^n\hat{g}^n - f^n g^n = \Delta f^ng^n - \hat{f}^n\Delta g^n$,
we obtain
\begin{eqnarray*}
\Delta v^n(x,t) & = &
\int_{-\infty}^{\infty}\Delta \tilde{\mathcal{G}}^{n-1}(x, t; y,0) v^{{n-1}}_0(y)\,\rmd y + \int_{-\infty}^{\infty}
\tilde{\hat{\mathcal{G}}}^{n-1}(x, t; y,0) \Delta v^{{n-1}}_0(y)\,\rmd y \\ &&
+\int_0^t\int_{-\infty}^{\infty} 
\Delta \tilde{\mathcal{G}}^{n-1}(x, t; y,s) 
S^{\zeta_*^{n-1}}(\zeta^{n-1},\dot\zeta^{n-1})(y,s)\,\rmd y\,\rmd s \\ &&
+ \int_0^t\int_{-\infty}^{\infty} 
\tilde{\hat{\mathcal{G}}}^{n-1}(x, t; y,s) 
\Delta S^{n-1}(y,s)\,\rmd y\,\rmd s \\ &&
-\int_0^t\int_{-\infty}^{\infty} \Delta \tilde{\mathcal{G}}^{n-1}_y(x,t; y,s)
\left[Q^{\zeta_*^{n-1}}(v^n)+R^{\zeta_*^{n-1}}(\zeta^{n-1},
v^n)\right](y,s)\,\rmd y\,\rmd s \\&&
-\int_0^t\int_{-\infty}^{\infty}\tilde{\hat{\mathcal{G}}}^{n-1}_y(x,t; y,s)
\left[\Delta Q^{n-1}+ \Delta R^{n-1} \right](y,s)\,\rmd y\,\rmd s,
\end{eqnarray*}
\begin{eqnarray}\label{vardotdeltan}
\Delta \dot{q}^n (t) & = &
-\int^\infty_{-\infty}\Delta \partial_t \pi^{n-1}_1(y,0,t)
v_0^{n-1}(y)\,\rmd y -\int^\infty_{-\infty}\partial_t \hat{\pi}^{n-1}_1(y,0,t)
\Delta v_0^{n-1}(y)\,\rmd y \\ \nonumber &&
- \int^t_0\int^{\infty}_{-\infty} \partial_t \Delta \pi^{n-1}_1(y,s,t)
S^{\zeta_*^{n-1}}(\zeta^{n-1},\dot\zeta^{n-1})(y,s)\,\rmd y\,\rmd s \\ \nonumber &&
- \int^t_0\int^{\infty}_{-\infty} \partial_t \hat{\pi}^{n-1}_1(y,s,t)
\Delta S^{n-1}(y,s)\,\rmd y\,\rmd s \\ \nonumber &&
+\int^t_0\int^{\infty}_{-\infty}\partial_t \partial_y \Delta \pi^{n-1}_1(y,s,t)
\left[ Q^{\zeta_*^{n-1}}(v^n)+ R^{\zeta_*^{n-1}}(\zeta^{n-1}, v^n)\right](y,s)\,\rmd y\,\rmd s,
\\ \nonumber &&
+\int^t_0\int^{\infty}_{-\infty} \partial_t \partial_y \hat{\pi}^{n-1}_1(y,s,t)
\left[\Delta Q^{n-1}+ \Delta R^{n-1}\right](y,s)\,\rmd y\,\rmd s,
\end{eqnarray}
\begin{eqnarray*}
\Delta \dot \tau^n (t) & = & -\int^\infty_{-\infty}\Delta \partial_t \pi^{n-1}_2(y,0,t)
v_0^{n-1}(y)\,\rmd y -\int^\infty_{-\infty}\partial_t \hat{\pi}^{n-1}_2(y,0,t)
\Delta v_0^{n-1}(y)\,\rmd y \\ \nonumber &&
- \int^t_0\int^{\infty}_{-\infty} \partial_t \Delta \pi^{n-1}_2(y,s,t)
S^{\zeta_*^{n-1}}(\zeta^{n-1},\dot\zeta^{n-1})(y,s)\,\rmd y\,\rmd s \\ \nonumber &&
- \int^t_0\int^{\infty}_{-\infty} \partial_t \hat{\pi}^{n-1}_2(y,s,t)
\Delta S^{n-1}(y,s)\,\rmd y\,\rmd s \\ \nonumber &&
+\int^t_0\int^{\infty}_{-\infty}\partial_t \partial_y \Delta \pi^{n-1}_2(y,s,t)
\left[ Q^{\zeta_*^{n-1}}(v^n)+ R^{\zeta_*^{n-1}}(\zeta^{n-1}, v^n)\right](y,s)\,\rmd y\,\rmd s \\ \nonumber &&
+\int^t_0\int^{\infty}_{-\infty} \partial_t \partial_y \hat{\pi}^{n-1}_2(y,s,t)
\left[\Delta Q^{n-1}+ \Delta R^{n-1}\right](y,s)\,\rmd y\,\rmd s,
\end{eqnarray*}
\begin{eqnarray}\label{vardeltasn}
\Delta \zeta_*^n & = &
\Delta \zeta_*^{n-1}-\int_{-\infty}^{\infty}\Delta \pi^{n-1}(y,0,\infty) v^{{n-1}}_0(y)\,\rmd y -\int_{-\infty}^{\infty}\hat{\pi}^{n-1}(y,0,\infty) \Delta v^{{n-1}}_0(y)\,\rmd y
\\ \nonumber &&
-\int_0^{\infty}\int_{-\infty}^{\infty} \Delta \pi^{n-1}(y, s, \infty) S^{\zeta_*^{n-1}}(\zeta^{n-1},\dot\zeta^{n-1})(y,s)\,\rmd y\,\rmd s \\ \nonumber &&
-\int_0^{\infty}\int_{-\infty}^{\infty} \hat{\pi}^{n-1}(y,s, \infty) \Delta S^{n-1} (y,s)\,\rmd y\,\rmd s \\ \nonumber &&
+\int_0^{\infty}\int_{-\infty}^{\infty}\Delta \pi_y^{n-1}(y,s, \infty)
\left[(Q^{\zeta_*^{n-1}}(v^n)+R^{\zeta_*^{n-1}}(\zeta^{n-1},
v^n)\right](y,s)\,\rmd y\,\rmd s \\ \nonumber &&
+\int_0^{\infty}\int_{-\infty}^{\infty}\hat{\pi}_y^{n-1}(y,s,\infty)
\left[\Delta Q^{n-1} + \Delta R^{n-1} \right] (y,s)\,\rmd y\,\rmd s,
\end{eqnarray}
and similarly for $\Delta \zeta^n$, where
\begin{eqnarray*}
\Delta v_0^{n-1}(x) & = &
\bar{u}^{\zeta_*^{n-1}}(x,0) -\bar{u}^{\hat \zeta_*^{n-1}}(x,0) \\ & = &
\frac{\partial \bar{u}^\zeta_*}{\partial \zeta_*}\Big|_{\zeta_*=\zeta_*^{n-1}}\Delta \zeta_*^{n-1}+
\mathrm{O}(|\Delta \zeta_*^{n-1}|^2\rme^{-\eta|x|}) \\ & = &
\mathrm{O}(|\Delta \zeta_*^{n-1}\rme^{-\eta|x|}|).
\end{eqnarray*}

Now define
\begin{equation}\label{xidef}
\xi(t) := \sup_{x\in\mathbb{R},\, s\in[0,t]}
\left(\frac{|\Delta v^n(x,s)|}{(\theta_\mathrm{gauss}+\theta_\mathrm{inner}+\theta_\mathrm{outer})(x,s)}
+|\Delta\dot\zeta^n(s)(1+s)| \right).
\end{equation}
Let $\tilde{r}:=|(\Delta\zeta^{n-1}, \Delta\zeta_*^{n-1})|_*=
|\Delta\zeta^{n-1}|_{B^1} + M|\Delta\zeta_*^{n-1}|$ be
sufficiently small. From \eqref{xidef} and smallness of $r$ and $\tilde{r}$ we obtain immediately
\begin{equation}\label{qrsvar}
|\Delta Q^{n-1} +\Delta R^{n-1}| \le C(r\xi(t)+r\tilde{r})(\Theta+ \Phi_1).
\end{equation}
Also, \eqref{sourcebds} and \eqref{derivsourcebds}  hold with $c(\beta^2 + 4C^2M_0^2)$ replaced by $Cr^2$. We use (\ref{ederivbds}) and (\ref{Gderivbds}) of Proposition~\ref{parambds} to obtain
\begin{equation}\label{deltae}
\Delta \pi^{n-1} \sim \pi\Delta\zeta_* \le \tilde{r} \pi,
\end{equation}
and similar appropriate bounds for 
$\Delta \tilde{\mathcal{G}}^{n-1}$ and its derivatives (of course, $\pi$ in (\ref{deltae})
is defined at a point between $\zeta_*^{n-1}$ and
$\hat\zeta_*^{n-1}$).
Next, using Lemmas~\ref{iniconvolutions}--\ref{occonvolutions} in a procedure parallel to the one used in the proof of Lemma~\ref{uexistence}, we obtain
\[
\xi (t) \le C(\tilde{r} + r\xi(t)).
\]
Since $\xi(t)$ is finite and continuous in $t$, this estimate implies that
\[
\xi(t)\le \frac{C\tilde{r}}{1-Cr}
\]
for a constant $C$ that is independent of $r$ and $\tilde{r}$. Now, replacing $\xi$ in (\ref{qrsvar}) with this bound, we substitute the result back into (\ref{vardotdeltan}) and into the similar formula for $\Delta \zeta^n$. Notice that, with the exception of the first two terms, the terms in \eqref{vardotdeltan} are all quadratic in their source term, so giving us small enough bounds. Hence, using once again Lemmas~\ref{iniconvolutions}--\ref{occonvolutions}, we obtain
\begin{equation}\label{vardotdeltabd}
|\Delta\dot\zeta^n|\le \left(CM_0\tilde{r} +\frac{C}{M}\tilde{r}+ Cr\tilde{r}\right)(1+t)^{-1},
\end{equation}
of which the two first terms in the right-hand side come from the first two terms of \eqref{vardotdeltan}. Similarly, we obtain
\begin{equation}\label{vardeltabd}
|\Delta\zeta^n|\le \left(CM_0 +\frac{C}{M}+Cr\right)\tilde{r}(1+t)^{-\frac12}.
\end{equation}
We notice that $(CM_0 +\frac{C}{M}+ Cr)$ can be made arbitrarily small, provided that $M_0$ and $r$ are small enough and $M$ is large enough. Next, we use \eqref{vardeltasn} to bound $\Delta\zeta_*^n$, using basically the same method used in (\ref{deslin})--(\ref{desnonlin}), and therefore obtaining
\[
M|\Delta\zeta_*^n| \le (CM_0+ Cr)\tilde{r}.
\]
This, together with \eqref{vardotdeltabd} and \eqref{vardeltabd}, gives us (\ref{lipbds}) with $\alpha<1$, finishing the proof of the Theorem~\ref{nonlin}.
\end{Proof}


\section{Summary and open problems}\label{s:conclusions}

In this paper, we considered time-periodic viscous Lax shocks $\bar{u}(x,t)$, which converge to constant time-independ\-ent\footnote{Viscous conservation laws, $u_t+f(u)_x=u_{xx}$, do not support genuinely time-periodic homogeneous rest states. Thus, time-periodic shock profiles can only admit time-independent asymptotic rest states} rest states $u_\pm$ as $x\to\pm\infty$. We showed that spectral stability of a time-periodic shock profile implies its nonlinear stability with respect to small initial perturbations that are smooth and sufficiently localized in space. Specifically, we proved that the corresponding solution converges with algebraic rate $1/\sqrt{t}$ in $L^\infty$ to an appropriate space- and time-translate of the shock profile $\bar{u}$.

The nonlinear stability proof followed the same strategy as in the case of stationary viscous shocks. First, the resolvent kernel of the linearization about the shock is extended meromorphically across the imaginary axis, and pointwise bounds are derived for this extension. These bounds together with analyticity make it possible to derive pointwise estimates for the resulting Green's function, using the inverse Laplace transform together with a careful deformation of the integration contour. Finally, a nonlinear iteration scheme that utilizes the anticipated spatio-temporal decay of perturbations closes the argument. The key difference from the case of stationary shocks is the time-periodicity of the underlying operators. To resolve this issue, we use spatial dynamics and exponential dichotomies to construct the analogue of the resolvent kernel and its meromorphic extension. One novel feature of our analysis is that we do not use Evans functions: instead, the geometric properties that ultimately determine the order of the pole of the meromorphic extension at the origin are encoded through Lyapunov--Schmidt reduction. Another advantage of our approach to meromorphic extensions through exponential dichotomies is that it is abstract and coordinate-free. This has several useful implications, which we shall now discuss.

First, as already discussed in \S\ref{introduction}, our results are also true for undercompressive, overcompressive, and mixed-type shock profiles, provided the notion of spectral stability is appropriately adapted and interpreted. Indeed, Lemma~\ref{l:fredholm} essentially reduces the meromorphic extension to finite-dimensional problems that are analogous to the stationary case, and therefore well understood. Similarly, the nonlinear iteration scheme relies primarily on the template functions that describe the anticipated temporal decay and are identical to those of the stationary case.

The fact that the Evans function is not used in our analysis opens up the possibility of extending our nonlinear-stability analysis to semi-discretizations of viscous shocks. Two-sided spatial finite difference approximations of viscous conservation laws lead to lattice differential equations. Their travelling waves satisfy functional differential equations of mixed type that are ill-posed as initial-value problems. More importantly, their stable and unstable eigenspaces are infinite-dimensional, and it is therefore unclear how Evans functions could be defined (we refer to \cite{Lord00,OhSandstede08} for situations where this can be done through Galerkin approximations). It was, however, shown in \cite{HaerterichSandstedeScheel02,MPVL08} that these equations admit exponential dichotomies. Thus, it is feasible that our approach could be used to establish nonlinear stability of semi-discretized shock profiles, thereby extending the analysis carried out in \cite{BenzoniGavageHuotRousset03} for one-sided finite differences.

We end this paper with three open problems. First, though we attempted to prepare the ground for a future analysis that addresses quasilinear, partially parabolic systems, we did not actually consider real viscosity here. We expect, however, that our main strategy should be applicable to such systems as in the related analyses of \cite{RaoofiZumbrun07,TexierZumbrun06_gas}. Another open problem is to find concrete examples of time-periodic Lax shocks. While previous analyses have shown that time-periodic Lax shocks bifurcate at Hopf bifurcations from stationary shocks, no examples are yet known where these bifurcations occur. The last issue is of a technical nature. Our construction of the resolvent kernel relied on an iterative Birman--Schwinger or parametrix-type argument, which is quite involved and does not immediately give expansions or bounds of the resolvent kernel. As discussed in \S\ref{S:greens_functions}, it should be possible to use exponential dichotomies directly to construct the resolvent kernel, but we have so far not been able to make this argument rigorous.

\begin{Acknowledgment}
Margaret Beck was partially supported under NSF grant DMS-0602891. Bj\"orn Sandstede gratefully acknowledges a Royal Society--Wolfson Research Merit Award. Kevin Zumbrun was partially supported under NSF grant DMS-0300487.
\end{Acknowledgment}



\end{document}